\documentclass[a4paper,11pt]{article}

\usepackage[margin=0.8in]{geometry}
\usepackage{adjustbox}
\usepackage{algorithm}
\usepackage{algorithmic}
\usepackage{amsmath,amssymb,amsfonts,amsthm}
\usepackage{array}
\usepackage{blkarray}
\usepackage{bm}
\usepackage{booktabs}
\usepackage[numbers,sort&compress]{natbib}
\usepackage{hyperref}
\hypersetup{
    colorlinks,
    linkcolor={red!50!black},
    citecolor={blue!50!black},
    urlcolor={blue!80!black}
}
\usepackage[nameinlink,capitalise,noabbrev]{cleveref}
\usepackage{diagbox}
\usepackage{dsfont}
\usepackage{enumerate}
\usepackage[inline]{enumitem}
\usepackage{graphicx}
\usepackage{mathtools}
\usepackage{todonotes}
\usepackage{multirow}
\usepackage{nth}
\usepackage{siunitx}
\usepackage{soul}
\usepackage{stmaryrd}
\usepackage{subcaption}
\usepackage{textcomp}
\usepackage{xcolor}
\usepackage{parskip}
\usepackage{url}
\usepackage{lipsum}
\usepackage{tikz}
\usepackage{tikzsymbols}
\usepackage{makecell}
\usepackage{pgfplots}
\pgfplotsset{compat=newest} 
\usetikzlibrary{arrows,shapes}
\usetikzlibrary{decorations.markings}
\usetikzlibrary{decorations.pathmorphing}
\usetikzlibrary{arrows.meta}
\usetikzlibrary{shadows,arrows,positioning,shapes.geometric}
\tikzset{cross/.style={cross out, draw=black, inner sep=0pt, outer sep=0pt},cross/.default={1pt}}
\usepackage{tikzsymbols}
\usepackage{pgfplots}
\pgfplotsset{compat=newest} 
\usetikzlibrary{plotmarks}
\usetikzlibrary{arrows.meta}
\usepgfplotslibrary{patchplots}
\usepackage{grffile}
\usepackage{pdflscape}
\usepackage{fontawesome}

\usepackage{lineno}
\pgfplotsset{every axis/.append style={
            label style={font=\scriptsize},
            tick label style={font=\scriptsize},
            legend style={font=\scriptsize}
        }}
\pgfplotsset{compat=newest}
\pgfplotsset{plot coordinates/math parser=false}
\pgfplotsset{grid style={dotted,gray}}
\pgfplotsset{
    compat=1.11,
    legend image code/.code={
            \draw[mark repeat=2,mark phase=2]
            plot coordinates {
                    (0cm,0cm)
                    (0.15cm,0cm)        
                    (0.3cm,0cm)         
                };%
        }
}
\usepgfplotslibrary{groupplots}
\newlength\figureheight
\newlength\figurewidth
\pgfdeclarelayer{background}
\pgfdeclarelayer{foreground}
\pgfsetlayers{background,main,foreground}
\usetikzlibrary{external}

\usepackage{aliascnt}

\newaliascnt{proposition}{theorem}

\aliascntresetthe{proposition}

\newaliascnt{lemma}{theorem}

\aliascntresetthe{lemma}

\newaliascnt{assumption}{theorem}

\aliascntresetthe{assumption}

\theoremstyle{definition}

\newaliascnt{definition}{theorem}

\aliascntresetthe{definition}

\newaliascnt{remark}{theorem}
\newtheorem{remark}[remark]{Remark}
\aliascntresetthe{remark}

\newaliascnt{example}{theorem}

\aliascntresetthe{example}

\crefname{theorem}{Theorem}{Theorems}
\Crefname{theorem}{Theorem}{Theorems}
\crefname{corollary}{Corollary}{Corollaries}
\Crefname{corollary}{Corollary}{Corollaries}
\crefname{proposition}{Proposition}{Propositions}
\Crefname{proposition}{Proposition}{Propositions}
\crefname{lemma}{Lemma}{Lemmas}
\Crefname{lemma}{Lemma}{Lemmas}
\crefname{assumption}{Assumption}{Assumptions}
\Crefname{assumption}{Assumption}{Assumptions}
\crefname{remark}{Remark}{Remarks}
\Crefname{remark}{Remark}{Remarks}
\crefname{example}{Example}{Examples}
\Crefname{example}{Example}{Examples}
\crefname{definition}{Definition}{Definitions}
\Crefname{definition}{Definition}{Definitions}
\crefname{algorithm}{Algorithm}{Algorithms}
\Crefname{algorithm}{Algorithm}{Algorithms}

\renewcommand{\vec}[1]{\boldsymbol{#1}}
\newcommand{\x}{\vec{x}}

\renewcommand{\d}{\mathrm{d} }
\newcommand{\flowmap}{\vec{\Phi}}
\newcommand{\submap}{\flowmap}

\newcommand{\R}{\ensuremath{\mathbb{R}}} 

\newcommand{\dx}{\ensuremath{{\mathrm d}x}} 
\newcommand{\dv}{\ensuremath{{\mathrm d}v}} 
\newcommand{\dt}{\ensuremath{{\mathrm d}t}} 

\DeclareMathOperator*{\argmin}{arg\,min}

\usepackage{xcolor}

\ifdefined\tikzexternalrealjob
    \usepackage[acronym,nonumberlist]{glossaries}
\else
    \usepackage[acronym,automake,nonumberlist]{glossaries}
\fi

\ifdefined\tikzexternalrealjob\else
    \makeglossaries
\fi

\graphicspath{{figures/}{figures/euler/}{figures/vlasov_poisson/}}

\begin{document}


\title{\textbf{A neural characteristic mapping method}\\[0.5em]
    \large Lagrangian PINNs based on flow maps for transport-dominated problems}
\date{\today}

\author{
    M. Campos Pinto\footnotemark[2]
    \and V. Fournet$^{*,\dagger}$
    \and E. Franck\footnotemark[3]
    \and P. Krah\footnotemark[4]
    \and V. Michel-Dansac\footnotemark[3]
}

\maketitle
{\renewcommand{\thefootnote}{\fnsymbol{footnote}}%
    \footnotetext[1]{Corresponding author.}%
    \footnotetext[2]{Max Planck Institute for Plasma Physics, Boltzmannstr. 2, Garching, 85748, Germany}%
    \footnotetext[3]{Université de Strasbourg, CNRS, Inria, IRMA, F-67000, Strasbourg, France}%
    \footnotetext[4]{CEA, IRFM, F-13108 Saint Paul-lez-Durance, France}}
\setcounter{footnote}{0}


\begin{abstract}
We propose a novel numerical method for transport-dominated problems using Physics-Informed Neural Networks (PINNs).
In such problems, the solution can develop large gradients and fine structures over time.
This is particularly challenging for classical PINNs, which represent the solution itself.
In our approach, rather than the solution, we choose to approximate the flow map of the equation.
The solution is then recovered by pulling back the exact initial condition,
so that its fine structures are produced by the map rather than represented by the network.
Our method is constructed by approximating the flow map of the equation with multiple PINNs,
one per time subinterval, and composing them using the semigroup property of the flow,
so that each network only has to represent a map close to the identity.
To highlight the advantages of our method
compared to classical PINNs or classical numerical schemes,
we present numerical results on the linear advection equation,
the incompressible Euler equation in vorticity formulation,
the Vlasov-Poisson equation,
and finally a drift-kinetic equation.

\end{abstract}



\section{Introduction} 
This work addresses transport-dominated problems in computational physics of the form
\begin{equation*}
    \partial_t \rho
    + \nabla \cdot \big( \vec{a}[\rho] \rho\big) = 0,
\end{equation*}
where $\rho \coloneqq \rho (t,\vec{x})$ is a density and $\vec{a}[\rho] \coloneqq \vec{a}[\rho](t,\vec{x})$ is a velocity field that may depend on $\rho$. Examples include kinetic models of magnetized plasmas~\cite{krall1986principles,GyselaX-semi-Lagrangian2016}, posed in a phase space of up to six dimensions, and the incompressible Euler equations~\cite{MajdaBertozzi2002}, which advect vorticity along the velocity field it generates.

The solution is the initial condition transported along the characteristics~\cite{Golse2013}, so the numerical difficulty lies not in the transport itself but in the increasing complexity of the solution: stretching and folding generate progressively finer scales that grid-based methods must resolve explicitly~\cite{ArberVann2002,YinMercierYadavSchneiderNave2021}.

Neural networks provide a meshfree representation whose complexity depends on the number of trainable parameters rather than on a prescribed grid~\cite{DeVoreHaninPetrova2021,FranckMichel-DansacNavoretVigon2026}.
In physics-informed neural networks (PINNs)~\cite{RaissiPerdikarisKarniadakis2019,Baty2024}, for instance,
these parameters are determined by minimizing the residual of the governing equation,
without requiring reference data.
Applied directly to the transported quantity,
however, a neural representation faces the same difficulty as grid-based methods,
since a single network must represent a field whose spatial complexity increases in time.
This difficulty is related to the slow decay of the Kolmogorov $n$-width of transport-dominated solution manifolds~\cite{GreifUrban2019,Mojgani_R_2023_j-comput-meth-appl-mech-eng_kolmogorov_wlpinnccmcdpde}.

We therefore adopt a Lagrangian viewpoint, which is the foundation of the characteristic mapping method.
Rather than approximating the transported quantity itself, we approximate with a neural network the flow map $\submap_0^t$, which sends a point at time $t$ back to its position at time $0$, and recover the solution by pulling back the initial condition along this map. This formulation shifts the approximation problem from evolving the solution to the characteristic flow, allowing structural properties of the map (invertibility, periodicity, and symplecticity for Hamiltonian velocity fields) to be directly incorporated into the representation. In the following, we briefly review the state of the art in characteristic mapping methods.

\subsubsection*{State of the art on characteristic mapping methods}

Characteristic mapping methods represent the flow map on a coarse grid and decompose it into submaps, allowing resolutions far beyond the underlying grid to be reached through composition. Such methods have been developed for transport equations with a prescribed velocity field, $ \vec{a}[\rho] (t,\vec{x})=\vec{a}(t,\vec x)$~\cite{MercierYinNave2019}, for nonlinear Vlasov-Poisson equations with self-consistent electric fields~\cite{KrahYinBergmannNaveSchneider2024,KrahLinWilhelmBacchiniNaveGrandgirardSchneider2026,KirchhartWilhelm2024,WilhelmKormann2025}, for multi-species Vlasov equations~\cite{WilhelmEifertTorrilhonOrland2025,WilhelmKrahSchneiderBacchiniGrandgirard2026}, and for extensions including BGK collisions~\cite{YinKrahLinNaveSchneider2026}. They have also been applied to the two- and three-dimensional incompressible Euler equations in vorticity form~\cite{YinMercierYadavSchneiderNave2021,YinSchneiderNave2023} and to magnetohydrodynamics~\cite{YinKrahNaveSchneider2026}. In all these approaches, the map is represented on a grid and the accuracy of the resulting solution is ultimately limited by the interpolation used to represent the map.

Replacing this grid representation by a neural network is a more recent development.
Neural networks have been used to approximate particle trajectories of the
six-dimensional Vlasov-Poisson system~\cite{fournet2026nonlinearcontrolvariatedeltaf} or
to represent the iterated flow maps of the numerical flow iteration~\cite{KirchhartWilhelm2024}.
Most closely related to the present work, the authors of~\cite{MontiMorozovskaPoggiShukla2026}
learn a flow map for the one-dimensional Vlasov-Poisson system with a PINN.
Neural flow maps have also been explored outside the PINN setting in computer graphics,
where a flow map and its Jacobian are learned from a simulated velocity field and
used to transport vorticity with low numerical dissipation in~\cite{deng2023fluid}.
In contrast to the present approach,
the map in that setting is fitted to simulation data rather
than determined from the residual of the governing equation.
More broadly, some machine learning approaches for the Vlasov-Poisson
system are reviewed in~\cite{FranckMichel-DansacNavoretVigon2026},
which also proposes a neural extension of the classical semi-Lagrangian method.

\subsubsection*{Scope and aim of the paper}

\begin{figure*}[htp!]
    \centering
    \setlength{\figurewidth}{0.85\linewidth}
    \setlength{\figureheight}{0.3\linewidth}
    {\scriptsize%
    {\pgfplotsset{compat=1.18}
\definecolor{mplcolor0}{rgb}{0.254900,0.411760,0.882350}
\pgfplotsset{colormap={turbo}{rgb255(0.0000)=(48,18,59) rgb255(0.0159)=(53,30,88) rgb255(0.0317)=(57,42,115) rgb255(0.0476)=(61,53,139) rgb255(0.0635)=(64,64,162) rgb255(0.0794)=(66,75,181) rgb255(0.0952)=(68,86,199) rgb255(0.1111)=(70,97,214) rgb255(0.1270)=(70,107,227) rgb255(0.1429)=(71,118,238) rgb255(0.1587)=(70,128,246) rgb255(0.1746)=(69,138,252) rgb255(0.1905)=(66,148,255) rgb255(0.2063)=(61,158,254) rgb255(0.2222)=(55,168,250) rgb255(0.2381)=(47,178,244) rgb255(0.2540)=(39,190,233) rgb255(0.2698)=(32,199,223) rgb255(0.2857)=(27,208,213) rgb255(0.3016)=(24,215,202) rgb255(0.3175)=(24,222,192) rgb255(0.3333)=(26,228,182) rgb255(0.3492)=(32,234,172) rgb255(0.3651)=(42,239,161) rgb255(0.3810)=(53,243,148) rgb255(0.3968)=(67,247,135) rgb255(0.4127)=(82,250,122) rgb255(0.4286)=(97,252,108) rgb255(0.4444)=(113,254,95) rgb255(0.4603)=(128,255,83) rgb255(0.4762)=(143,255,73) rgb255(0.4921)=(156,254,64) rgb255(0.5079)=(169,251,57) rgb255(0.5238)=(180,248,54) rgb255(0.5397)=(190,244,52) rgb255(0.5556)=(200,239,52) rgb255(0.5714)=(210,233,53) rgb255(0.5873)=(219,226,54) rgb255(0.6032)=(227,219,56) rgb255(0.6190)=(235,211,57) rgb255(0.6349)=(241,203,58) rgb255(0.6508)=(246,195,58) rgb255(0.6667)=(250,186,57) rgb255(0.6825)=(252,177,54) rgb255(0.6984)=(254,167,50) rgb255(0.7143)=(254,155,45) rgb255(0.7302)=(254,144,41) rgb255(0.7460)=(252,132,35) rgb255(0.7619)=(249,117,29) rgb255(0.7778)=(245,105,24) rgb255(0.7937)=(241,93,19) rgb255(0.8095)=(236,83,15) rgb255(0.8254)=(231,73,12) rgb255(0.8413)=(225,65,9) rgb255(0.8571)=(218,57,7) rgb255(0.8730)=(210,49,5) rgb255(0.8889)=(202,42,4) rgb255(0.9048)=(193,35,2) rgb255(0.9206)=(183,29,2) rgb255(0.9365)=(172,23,1) rgb255(0.9524)=(161,18,1) rgb255(0.9683)=(149,13,1) rgb255(0.9841)=(136,8,2) rgb255(1.0000)=(122,4,3)}}
\begin{tikzpicture}
\begin{axis}[
  at={(0.0083\figurewidth, 0.5299\figureheight)},
  anchor=south west,
  width=0.1207\figurewidth,
  height=0.3319\figureheight,
  title={submap $\submap_0^{1}$},
  xmin=0, xmax=1,
  ymin=0, ymax=1,
  xtick=\empty,
  ytick=\empty,
  title style={font=\tiny, at={(0.5,1)}, anchor=south, yshift=-5pt},
  scale only axis
]
\addplot graphics [
  xmin=0, xmax=1,
  ymin=0, ymax=1
] {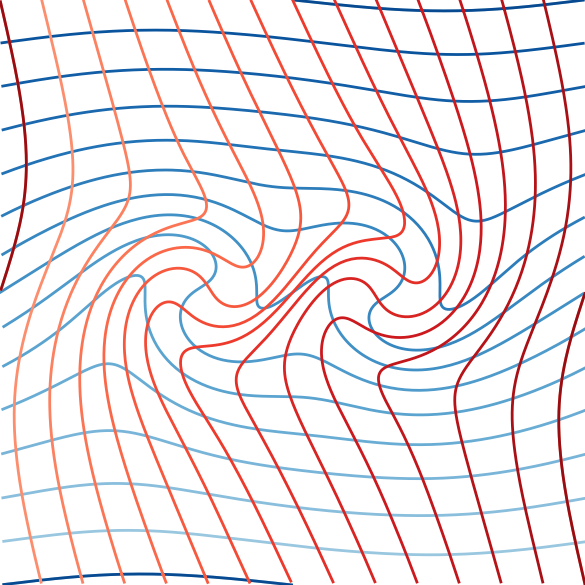};
\end{axis}

\begin{axis}[
  at={(0.1785\figurewidth, 0.5299\figureheight)},
  anchor=south west,
  width=0.1207\figurewidth,
  height=0.3319\figureheight,
  title={submap $\submap_1^{2}$},
  xmin=0, xmax=1,
  ymin=0, ymax=1,
  xtick=\empty,
  ytick=\empty,
  title style={font=\tiny, at={(0.5,1)}, anchor=south, yshift=-5pt},
  scale only axis
]
\addplot graphics [
  xmin=0, xmax=1,
  ymin=0, ymax=1
] {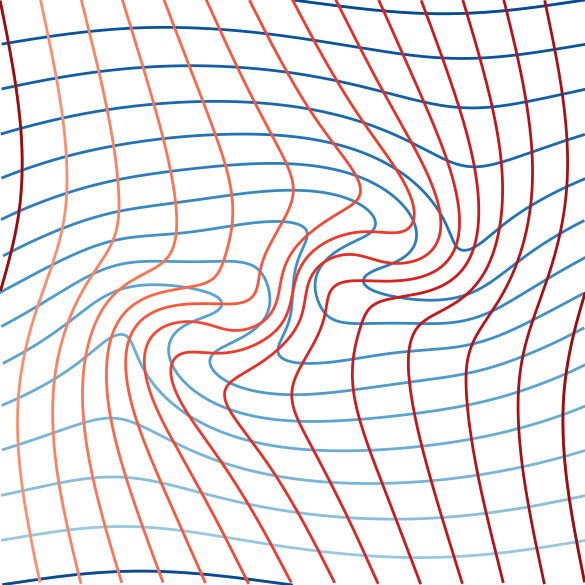};
\end{axis}

\begin{axis}[
  at={(0.0083\figurewidth, 0.0227\figureheight)},
  anchor=south west,
  width=0.1207\figurewidth,
  height=0.3319\figureheight,
  title={submap $\submap_2^{3}$},
  xmin=0, xmax=1,
  ymin=0, ymax=1,
  xtick=\empty,
  ytick=\empty,
  title style={font=\tiny, at={(0.5,1)}, anchor=south, yshift=-5pt},
  scale only axis
]
\addplot graphics [
  xmin=0, xmax=1,
  ymin=0, ymax=1
] {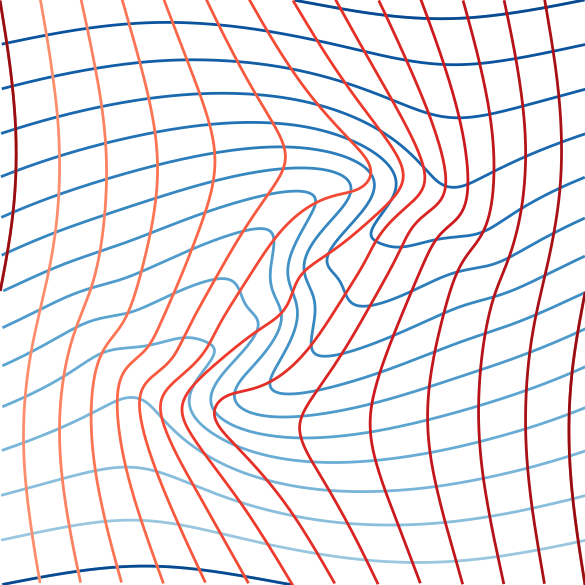};
\end{axis}

\begin{axis}[
  at={(0.1785\figurewidth, 0.0227\figureheight)},
  anchor=south west,
  width=0.1207\figurewidth,
  height=0.3319\figureheight,
  title={submap $\submap_3^{4}$},
  xmin=0, xmax=1,
  ymin=0, ymax=1,
  xtick=\empty,
  ytick=\empty,
  title style={font=\tiny, at={(0.5,1)}, anchor=south, yshift=-5pt},
  scale only axis
]
\addplot graphics [
  xmin=0, xmax=1,
  ymin=0, ymax=1
] {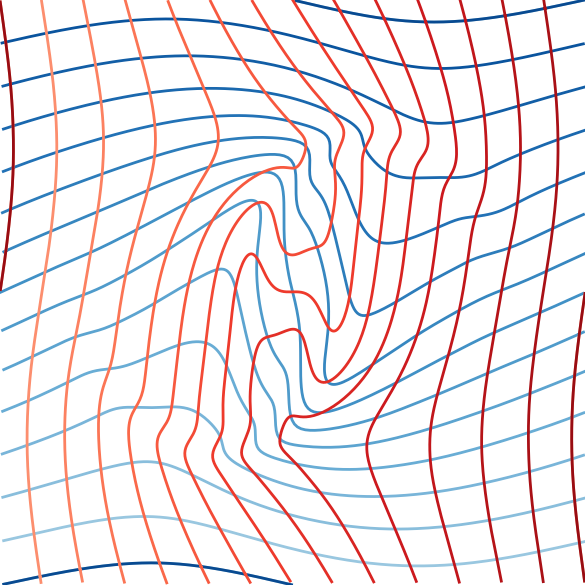};
\end{axis}

\begin{axis}[
  at={(0.3159\figurewidth, 0.0373\figureheight)},
  anchor=south west,
  width=0.2979\figurewidth,
  height=0.8188\figureheight,
  title={Composed map $\submap_{0}^{4}=\submap_{0}^{1}\circ\submap_{1}^{2}\circ\submap_{2}^{3}\circ\submap_{3}^{4}$},
  xmin=0, xmax=1,
  ymin=0, ymax=1,
  xtick=\empty,
  ytick=\empty,
  title style={font=\tiny, at={(0.5,1)}, anchor=south, yshift=-5pt},
  scale only axis
]
\addplot graphics [
  xmin=0, xmax=1,
  ymin=0, ymax=1
] {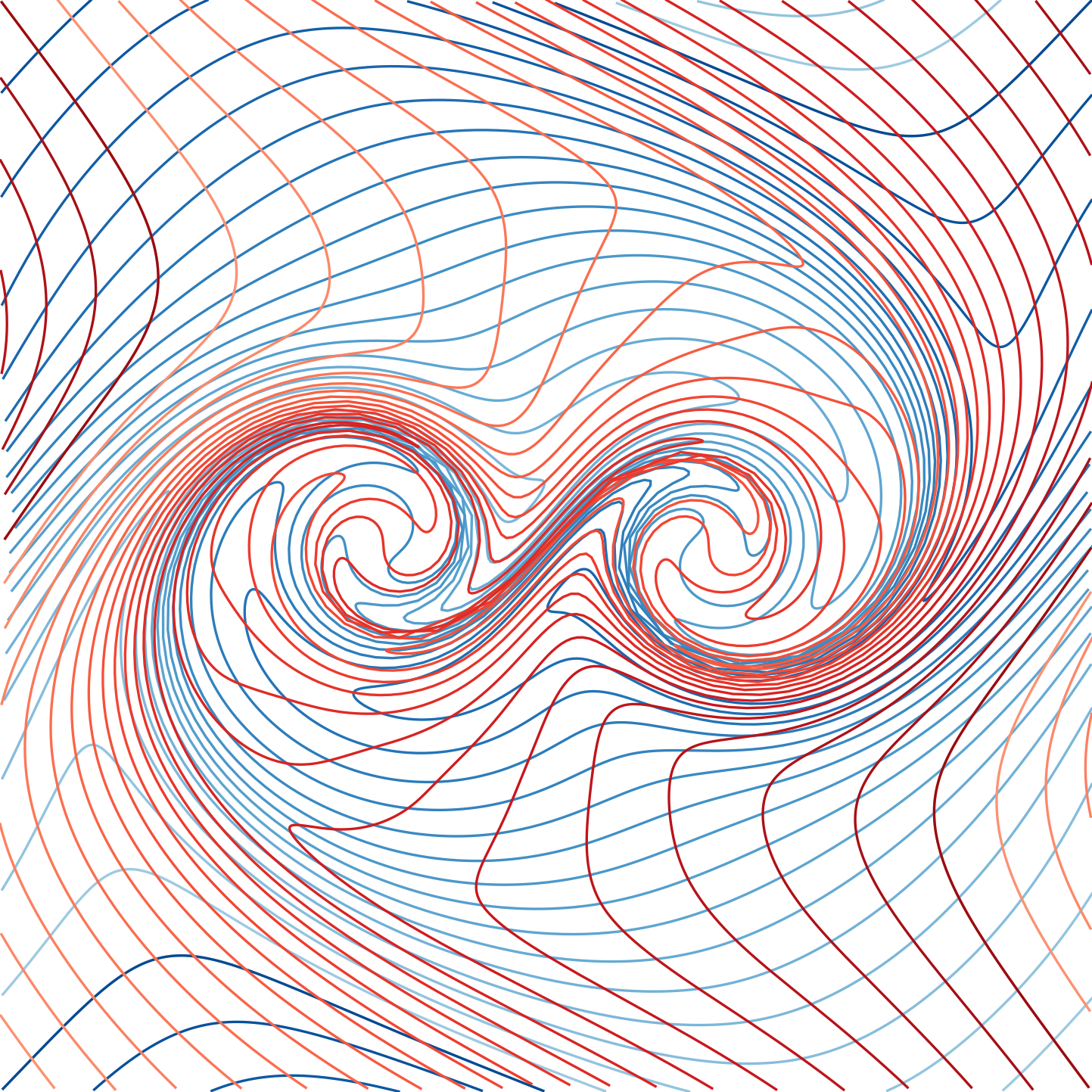};
\end{axis}

\begin{axis}[
  at={(0.6402\figurewidth, 0.0373\figureheight)},
  anchor=south west,
  width=0.2979\figurewidth,
  height=0.8188\figureheight,
  title={Vorticity $\omega(4,\cdot)=\omega_0\circ\submap_{0}^{4}$},
  xmin=0, xmax=0.996094,
  ymin=0, ymax=0.996094,
  xtick=\empty,
  ytick=\empty,
  point meta min=0,
  point meta max=5.69743,
  colormap name=turbo,
  colorbar,
  colorbar style={width=0.0132\figurewidth},
  title style={font=\tiny, at={(0.5,1)}, anchor=south, yshift=-5pt},
  scale only axis
]
\addplot graphics [
  xmin=0, xmax=0.996094,
  ymin=0, ymax=0.996094
] {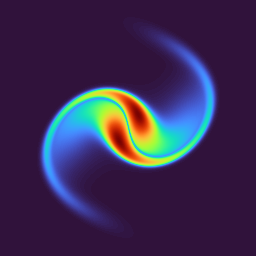};
\end{axis}
\node[anchor=north, text=mplcolor0] at (0.1537\figurewidth, 1.0473\figureheight) {learned neural submaps};
\end{tikzpicture}
}%
}
    \caption{Learning the solution of the two-dimensional Euler equations through its flow. Left: the backward submaps $\submap_{k-1}^{k}$ of the four subintervals $k=1,\dots,4$, each learned by a separate network and each close to the identity, visualized by the image of a Cartesian grid. Middle: their composition, the global backward map at $t=4$. Right: the vorticity $\omega(4,\cdot)=\omega_0\circ\submap_{0}^4$ obtained by pulling the initial condition back through the composed map. The test case is the vortex merger of \cref{sec:euler_vortex_merger}.}
    \label{fig:intro_flow_maps}
\end{figure*}

The aim of this work is to learn transport solutions through their underlying flow maps while combining the characteristic mapping decomposition with neural network representations.
Each subinterval of the evolution is represented by a separate network, so that the individual maps remain close to the identity even when their composition develops strong deformation.
As illustrated in \cref{fig:intro_flow_maps}, the increasing complexity of the solution is thus transferred from the representation of an individual network to the composition of successive maps.

\subsubsection*{Contribution}

The main contributions of this work are:
\begin{itemize}
    \item We develop a software framework inside the scimba library \cite{scimba} for Lagrangian PINN representations that use a composition of multiple submaps.
    \item We show that the submaps can be trained sequentially with a cost per training problem independent of the final simulation time when the velocity does not depend on the solution, and compare different neural network architectures for their representation.
    \item We introduce a Picard iteration for velocity fields coupled through nonlinear equations, and assess the resulting methods on several transport-dominated problems against classical methods.
\end{itemize}

\subsubsection*{Organization of the manuscript}
The remainder of the manuscript is organized as follows. \Cref{sec:nummeth}
introduces the flow map, its semigroup property and the decomposition into
submaps that the rest of the paper relies on. \Cref{sec:pinns} presents the
three PINN formulations, their coupling to a field equation and the network
architectures used for the maps. \Cref{sec:results} reports the numerical
results, and \cref{sec:concl} concludes.

\section{Flow Maps and the Semigroup Property}
\label{sec:nummeth}

In this section we define the mathematical framework
and the relevant properties of the flow map used in our Lagrangian PINN formulation in \cref{sec:pinns}.
The material is well known and we refer to \cite[Ch.~2]{Golse2013} for further details and rigorous proofs.

\subsection{Setting}
\label{sec:setting}
Let $\Omega\subset\R^d$ be the domain of the unknown, periodic in some of its
directions, and let $T>0$ be a final time. We look for
$\rho\colon[0,T]\times\Omega\to\R$ solving the Cauchy problem,
\begin{equation}
    \partial_t \rho + \nabla\cdot\big(\vec a[\rho] \, \rho  \big) = 0
    \quad\text{in }(0,T]\times\Omega,
    \qquad
    \rho(0,\cdot)=\rho_0 ,
    \label{eq:cauchy_problem}
\end{equation}
in which the velocity $\vec a[\rho]\colon[0,T]\times\Omega\to\R^d$ may depend nonlinearly and nonlocally on $\rho$.

\subsection{Flow maps}
\label{sec:flowmaps}

For $s,t\in[0,T]$, let $\submap_{s}^{t}\colon\Omega\to\Omega$ denote the flow map, defined for fixed $t$ by solving, as a function of $s$,
\begin{equation}
    \frac{\d}{\d s} \submap_{s}^{t}(\vec{x})
    =
    \vec a[\rho]\big(s,{\submap}_{s}^{t}(\vec{x})\big),
    \qquad
    {\submap}_{t}^{t}(\vec{x})=\vec{x} .
    \label{eq:characteristics}
\end{equation}
In words, $\submap_s^{t}(\vec x)$ is the position at time $s$ of the (backward, if $s<t$)
characteristic curve passing through $\vec x$ at time $t$.
The upper index $t$ is the time at which $\submap_t^{t}(\vec x) = \vec x$,
while the lower index $s$ represents
the time at which the characteristic curve is evaluated.
Hence, $\submap_0^{t}$ sends a point of the domain to its origin at $t=0$.
When $\vec a$ depends on $\rho$, the flow map is coupled to the solution of \eqref{eq:cauchy_problem} and must be solved simultaneously with it.

The map satisfies the group properties: for all
$r,s,t\in[0,T]$,
\begin{equation}
    {\submap}_{t}^{t} = \mathrm{Id}_\Omega,
    \qquad
    {\submap}_{r}^{s}\circ {\submap}_{s}^{t}
    =
    {\submap}_{r}^{t},
    \qquad
    \left({\submap}_{s}^{t}\right)^{-1}={\submap}_{t}^{s}.
    \label{eq:semi-group-properties}
\end{equation}
Moreover, viewed as a function of $(t,\vec x)$ for fixed $s$ rather than as the
solution of \eqref{eq:characteristics}, it also solves the transport equation
\begin{equation}
    \partial_t \submap_{s}^{t}(\vec x)
    + \big(\vec a[\rho](t,\vec x)\cdot\nabla_{\vec x}\big)\,\submap_{s}^{t}(\vec x)
    = 0,
    \qquad
    \submap_{s}^{s}(\vec x)=\vec x .
    \label{eq:map_transport}
\end{equation}

The solution of \eqref{eq:cauchy_problem} is then reconstructed from the map
alone by
\begin{equation}
    \rho(t,\vec{x})
    =
    \rho_0\!\big({\submap}_{0}^{t}(\vec{x})\big)
    \left|\det\!\left(\nabla {\submap}_{0}^{t}(\vec{x})\right)\right| ,
    \label{eq:comp_with_inicond}
\end{equation}
so knowledge of $\submap_0^{t}$ together with the initial condition is sufficient to recover $\rho$, see~\cite[Thm.~2.3.6]{Golse2013}. If
$\vec a(t,\cdot)$ is divergence-free, the Liouville formula gives
$|\det\nabla\submap_s^{t}|\equiv 1$, so the flow is
volume-preserving, see~\cite[Thm.~2.2.3]{Golse2013}. It is in addition symplectic
when the velocity is Hamiltonian, which is the case for the Vlasov models considered
here.
The reconstruction then reduces
to the pure pullback $\rho(t,\vec x)=\rho_0(\submap_0^{t}(\vec x))$, and every
functional $\int_\Omega G(\rho)\,\d\vec x$, in particular every $L^p$ norm,
is exactly conserved.

\begin{remark}
    \label{rem:absolute_value}
    The Liouville formula
    shows (see~\cite[Thm.~2.2.3]{Golse2013}) that the exact flow is orientation-preserving, so the absolute value in
    \eqref{eq:comp_with_inicond} is redundant. This is no longer the case for the
    network approximations of \cref{sec:flow_networks}, whose Jacobian determinant
    is positive by construction only for the symplectic family.
    We therefore keep the absolute value in \eqref{eq:comp_with_inicond}.
    In all the numerical tests of \cref{sec:results}, however, the velocity field
    is divergence-free, so that the exact flow has unit Jacobian. In the remainder
    of the paper, the solution is therefore reconstructed by the pure pullback
    $\rho(t,\vec x)=\rho_0(\submap_0^{t}(\vec x))$, i.e., with the determinants
    in \eqref{eq:density_reconstruction} and \eqref{eq:composed_density} set to
    one, also for networks that are not volume-preserving.
\end{remark}

\subsection{Flow decomposition into submaps}
\label{sec:submaps}

Let $0=t_0<t_1<\dots<t_{N_t}=T$ be a partition of the time interval, with
$\Delta t_i = t_{i+1}-t_i$. Applying the semigroup property repeatedly, the
backward map at a time $t\in[t_k,t_{k+1}]$ factorizes into submaps,
\begin{equation}
    {\submap}_{0}^{t}
    =
    \submap_{t_0}^{t_1}\circ \submap_{t_1}^{t_2}\circ \dots
    \circ \submap_{t_{k-1}}^{t_k}\circ \submap_{t_k}^{t} ,
    \label{eq:comp_structure}
\end{equation}
where each factor $\submap_{t_i}^{t_{i+1}}$ depends only on the dynamics over
$[t_i,t_{i+1}]$.

Let $t \in [t_k, t_{k+1}]$ for some $k \in \{0, \dots, N_t-1\}$.
We introduce the intermediate points
\begin{equation}
    \vec y_k = \submap_{t_k}^{t}(\vec x)
    \text{\qquad and \qquad}
    \vec y_i = \submap_{t_i}^{t_{i+1}}(\vec y_{i+1}),
    \text{\quad for all \quad}
    i=k-1,\dots,0 ,
    \label{eq:intermediate_points}
\end{equation}
so that $\vec y_0=\submap_0^{t}(\vec x)$.
Then, computing the Jacobian determinant of
\eqref{eq:comp_structure} gives
\begin{equation}
    \left|\det\!\left(\nabla {\submap}_{0}^{t}(\vec x)\right)\right|
    =
    \left|\det\!\left(\nabla \submap_{t_k}^{t}(\vec x)\right)\right| \,
    \prod_{i=0}^{k-1}
    \left|\det\!\left(\nabla {\submap}_{t_i}^{t_{i+1}}(\vec y_{i+1})\right)\right| ,
    \label{eq:jacobian_product}
\end{equation}
which turns \eqref{eq:comp_with_inicond} into the composed reconstruction
formula
\begin{equation}
    \rho(t,\vec x)
    =
    \rho_0(\vec y_0)\,
    \left|\det\!\left(\nabla \submap_{t_k}^{t}(\vec x)\right)\right|
    \prod_{i=0}^{k-1}
    \left|\det\!\left(\nabla {\submap}_{t_i}^{t_{i+1}}(\vec y_{i+1})\right)\right| .
    \label{eq:comp_reconstruction}
\end{equation}

To further justify the use of submaps, we note that a submap transports the
points only over the short interval $[t_i,t_{i+1}]$. Whatever the deformation
accumulated since $t=0$, it moves each point by a small distance and deforms
the domain only slightly, so it is a small perturbation of the identity and
therefore easy to learn.
Hence, the complexity of the total map, which grows in time as the flow stretches
and folds the domain, is distributed over the $N_t$ submaps instead of being
carried by a single object.
This is exploited in
\cref{sec:lagrangian_pinn_semigroup} below and
by grid-based characteristic
mapping methods, e.g.,~\cite{MercierYinNave2019,YinMercierYadavSchneiderNave2021,KrahYinBergmannNaveSchneider2024}.

\section{Approximating Flow Maps with Physics-Informed Neural Networks (PINNs)}
\label{sec:pinns}

We approximate \eqref{eq:cauchy_problem} using neural networks, whose general purpose and formulation are introduced in \cref{sec:neural_networks}.
We then introduce three neural-network-based approaches to approximate this problem.
The Eulerian PINN, presented in \cref{sec:eulerian_pinn}, directly approximates the solution $\rho$.
The Lagrangian PINN (\cref{sec:lagrangian_pinn}), in contrast, approximates the full flow map $\submap_{0}^{t}$ introduced in \cref{sec:flowmaps}.
Finally, the Lagrangian PINN with the semigroup property
(\cref{sec:lagrangian_pinn_semigroup}) approximates each submap
$\submap_{t_i}^{t_{i+1}}$ of the global map
$
    {\submap}_{0}^{t}=
    \submap_{t_0}^{t_1}\circ \submap_{t_1}^{t_2}\circ \dots
    \circ \submap_{t_{k-1}}^{t_k}\circ \submap_{t_k}^{t}
$
using a separate neural network.
For clarity, we assume throughout this discussion that the velocity field $\vec{a}$ in the advection equation is given, i.e., that the transport equation is linear.
The coupling of the advection equation \eqref{eq:cauchy_problem} with a field equation is discussed later in \cref{sec:field_coupling}.
We conclude by discussing different network architectures for approximating the flow maps in \cref{sec:flow_networks}.

\subsection{Approximation of PDE solutions with neural networks}
\label{sec:neural_networks}

Solutions of PDEs like \eqref{eq:cauchy_problem} or \eqref{eq:map_transport}
are naturally defined in infinite-dimensional function spaces.
In practice, we approximate such solutions by restricting the search
to a finite-dimensional space of functions.
In classical numerical methods, like the finite element method,
this space is spanned by a set of basis functions,
and the solution is represented as a linear combination of these basis functions.
The problem of approximation is then reduced to finding
the coefficients (called degrees of freedom) of this linear combination.
They are generally obtained by orthogonally projecting the PDE onto the finite-dimensional space,
leading to a linear system for the coefficients.
However, such representations are often limited by the curse of dimensionality,
as the number of basis functions required to accurately represent the solution
grows exponentially with the dimension of the problem.
This discussion is expanded upon in~\cite{FranckMichel-DansacNavoretVigon2026},
in the specific context of advection-dominated problems.

To avoid this limitation, we replace this linear space spanned by basis functions
with a nonlinear space of functions.
More specifically, we use neural networks, which are parametric functions
$y \mapsto \vec f_\theta(y)$,
where $\theta$ are the trainable parameters of the network,
to be determined.
Much like in linear spaces, the best approximation of the solution
is obtained by projecting the PDE onto the (nonlinear) space of functions
represented by the neural network.
This means that $\theta$ can be viewed as
the degrees of freedom of the neural network approximation.
This no longer leads to a linear system,
but rather to a nonlinear optimization problem for the parameters $\theta$.
In the following, we describe multiple ways to formulate this optimization problem,
starting from the well-known PINN approach from~\cite{RaissiPerdikarisKarniadakis2019}.

\subsection{Eulerian PINNs}
\label{sec:eulerian_pinn}

The classical PINN approach approximates the solution itself by a network
$\rho_\theta\colon[0,T]\times\Omega\to\R$ and fits the parameters by minimizing the
residual of \eqref{eq:cauchy_problem}.
This leads to the optimization problem
\begin{equation}
    \theta = \argmin_{\vartheta} \left(
    \int_{0}^{T} \int_{\Omega}
    \left|\partial_t \rho_\vartheta + \nabla\cdot\big(\vec a[\rho_\vartheta] \, \rho_\vartheta \big)\right|^2\,\dx\,\dt
    + \lambda \int_{\Omega}\left|\rho_\vartheta(0,\vec x)-\rho_0(\vec x)\right|^2\,\dx
    \right),
    \label{eq:eulerian_loss}
\end{equation}
where $\lambda>0$ weights the initial condition.
Both integrals are estimated
by Monte Carlo quadrature over collocation points sampled in
$(0,T)\times\Omega$ and in $\Omega$,
and the derivatives are obtained by automatic differentiation.

The approach is mesh-free, it is subject to no time step restriction
(as there is no time discretization),
and it applies in high dimensions.
Its main weaknesses are of a different nature.
First, the optimization problem couples all times at once and is stiff,
so that quasi-second-order optimizers such as
natural gradient methods~\cite{MueZei2023,SchFur2025,JniKiyShuUrbAhmMue2026}
are needed to reach acceptable accuracy.
Second, as mentioned in the introduction, $\rho(t,\cdot)$
develops finer and finer structures as $t$ grows,
which makes it harder and harder for the (smooth) network to represent the solution.

\subsection{Lagrangian PINNs}
\label{sec:lagrangian_pinn}

To alleviate the second drawback of PINNs,
instead of the solution, we now approximate the flow map.
We take a single time interval $[0,T]$.
Let $\submap_\theta\colon[0,T]\times\Omega\to\Omega$ be a network
with parameters $\theta$, intended to approximate the backward map, as
\begin{equation}
    \submap_\theta(t,\cdot) \approx \submap_0^{t},
    \qquad t\in[0,T].
    \label{eq:flow_ansatz}
\end{equation}

\paragraph{Training the approximate flow map.} The parameters of $\submap_\theta$ are obtained
by minimizing the residual of the Cauchy problem
\eqref{eq:map_transport} that the map satisfies.
This yields the optimization problem
\begin{equation}
    \theta = \argmin_{\vartheta} \left(
    \int_{0}^{T}\!\!\int_{\Omega}
    \left\|\partial_t \submap_\vartheta
    + \big(\vec a[\rho_{\vartheta}]\cdot \nabla_{\vec x}\big)\,\submap_\vartheta\right\|^2\,\dx\,\dt
    + \lambda \int_{\Omega}\left\|\submap_\vartheta(0,\vec x) - \vec x\right\|^2\,\dx
    \right).
    \label{eq:flow_loss}
\end{equation}
Both integrals are once again approximated by Monte Carlo quadrature.
In contrast with \eqref{eq:eulerian_loss},
the unknown is vector-valued and the residual has $d$ components,
one per Lagrangian invariant transported by the flow.

\paragraph{Reconstructing the density.} Following the pullback formula
\eqref{eq:comp_with_inicond}, the solution associated with the trained map is
\begin{equation}
    \rho_\theta(t,\vec x)
    =
    \rho_0\big(\submap_\theta(t,\vec x)\big)
    \left|\det\nabla_{\vec x}\submap_\theta(t,\vec x)\right| .
    \label{eq:density_reconstruction}
\end{equation}
Note that, when $\vec{a}$ depends on $\rho$,
the loss \eqref{eq:flow_loss} cannot be evaluated as such.
This is resolved by the procedure described in \cref{sec:field_coupling}.

\paragraph{Properties of the method.}
We exhibit three main properties of this method.
\begin{itemize}[nosep]
    \item The initial density $\rho_0$ is used exactly and is
          never projected onto a discrete space, so the fine structures of $\rho$ are
          produced by composition with the map, and not by the network itself.
    \item The solution \eqref{eq:density_reconstruction} is the image of $\rho_0$ by the inverse
          of $\submap_\theta(t,\cdot)$, so as soon as the latter is a bijection of
          $\Omega$ the mass of $\rho_\theta(t,\cdot)$ equals the mass of $\rho_0$
          exactly, whatever the error on $\theta$.
          No penalty term enforcing
          conservation is necessary in \eqref{eq:flow_loss}.
    \item If the network is volume-preserving, the determinant
          in \eqref{eq:density_reconstruction} equals one,
          so $\rho_\theta(t,\cdot)$ and
          $\rho_0$ share the same range and all $L^p$ norms are conserved exactly.
\end{itemize}
This shows the respective roles of the two forms of the equation. The
residual \eqref{eq:flow_loss} that is minimized is a property of the flow
map alone and is insensitive to the equation satisfied by $\rho$, while the
conservative form \eqref{eq:cauchy_problem} enters the method only through the
Jacobian factor of \eqref{eq:density_reconstruction}, which would be absent if
$\rho$ were transported in advective form instead. For the divergence-free
models considered in this paper, the exact flow has unit Jacobian and this
factor is set to one (see \cref{rem:absolute_value}). Keeping it would make
the conservation of mass exact even for a network that is not
volume-preserving, but it would also add the volume error of the trained map
to the solution.
Note that this formulation has already been studied in \cite{MontiMorozovskaPoggiShukla2026}
for the one-dimensional Vlasov-Poisson system.

\paragraph{Limitation.} A single network has to represent $\submap_0^{t}$ for
every $t$ up to $T$. This map moves further away from the identity as the
flow stretches and folds the domain, so it becomes increasingly hard to
represent as $T$ grows. This is what the
decomposition of the next section addresses.

\subsection{Lagrangian PINNs with a semigroup property}
\label{sec:lagrangian_pinn_semigroup}

We now use the decomposition of \cref{sec:submaps}. Let
$0=t_0<t_1<\dots<t_{N_t}=T$ be the partition introduced there, and associate
with each interval a network of its own,
\begin{equation}
    \submap_{\theta_k}\colon[t_k,t_{k+1}]\times\Omega\to\Omega,
    \qquad
    \submap_{\theta_k}(t,\cdot)\approx\submap_{t_k}^{t},
    \qquad k=0,\dots,N_t-1 .
    \label{eq:submap_ansatz}
\end{equation}
Each network approximates one factor of \eqref{eq:comp_structure} and is
trained exactly as in \cref{sec:lagrangian_pinn}, on its own interval and with
its own identity condition,
\begin{equation}
    \min_{\theta_k} \left(
    \int_{t_k}^{t_{k+1}}\int_{\Omega}
    \left\|\partial_t \submap_{\theta_k}
    + \big(\vec a[\rho_{\Theta}]\cdot \nabla_{\vec x}\big)\,\submap_{\theta_k}\right\|^2\,\dx\,\dt
    + \lambda \int_{\Omega}\left\|\submap_{\theta_k}(t_k,\vec x) - \vec x\right\|^2\,\dx
    \right),
    \label{eq:submap_loss}
\end{equation}
estimated by Monte Carlo quadrature with $N_{\mathrm{colloc}}$ collocation points drawn in
$(t_k,t_{k+1})\times\Omega$ for the first term and $N_{\mathrm{ic}}$ points drawn in $\Omega$ for the initial condition.

\paragraph{Composed map and density.} Write $\Theta=(\theta_0,\dots,\theta_k)$
for the parameters of the submaps trained so far. The discrete counterparts of
\eqref{eq:intermediate_points} and \eqref{eq:comp_reconstruction} are, for
$t\in[t_k,t_{k+1}]$,
\begin{align}
    \vec y_k = \submap_{\theta_k}(t,\vec x),
    \qquad
    \vec y_i = \submap_{\theta_i}(t_{i+1},\vec y_{i+1}),
    \quad i=k-1,\dots,0 ,
    \label{eq:discrete_composition} \\
    \rho_\Theta(t,\vec x)
    =
    \rho_0(\vec y_0)\,
    \left|\det\nabla_{\vec x}\submap_{\theta_k}(t,\vec x)\right|
    \prod_{i=0}^{k-1}
    \left|\det\nabla\submap_{\theta_i}(t_{i+1},\vec y_{i+1})\right| .
    \label{eq:composed_density}
\end{align}
Every earlier submap is evaluated at the end point $t_{i+1}$ of its own
interval, which is the only value of it entering \eqref{eq:comp_structure}.

\paragraph{Sequential training with frozen submaps.} The networks are trained
one interval after the other. While $\theta_k$ is being optimized, the
parameters $\theta_0,\dots,\theta_{k-1}$ are frozen. They are used to evaluate
$\rho_\Theta$ in the field update, but they receive no gradient and are not
modified. Three consequences are worth stating.
\begin{enumerate}[nosep]
    \item Each optimization problem involves the parameters of a single
          network, so its size does not grow with $k$ and the networks may be
          chosen smaller than the single network of
          \cref{sec:lagrangian_pinn}.
    \item The initial condition of interval $[t_k,t_{k+1}]$ is the identity and not a density fitted at $t_k$. No additional error is therefore committed
          when passing from one interval to the next, and the continuity in
          time is carried by the composition \eqref{eq:comp_structure} rather
          than by the loss.
    \item Evaluating $\rho_\Theta$ at one point requires $k+1$ network
          evaluations, so the cost of the field update grows linearly with
          $k$. This is the price of the composition and it is visible in the
          training times reported in \cref{sec:vlasov_poisson}.
\end{enumerate}
The resulting procedure is summarized in \cref{alg:lagrangian_pinn}, whose
field and flow update steps
alternate as explained in \cref{sec:field_coupling} below.

\begin{algorithm}[!ht]
    \caption{Lagrangian PINN with semigroup composition.}
    \label{alg:lagrangian_pinn}
    \begin{algorithmic}
        \STATE \textbf{Input:} initial density $\rho_0$, partition $0=t_0<\dots<t_{N_t}=T$, networks $\submap_{\theta_0},\dots,\submap_{\theta_{N_t-1}}$ initialized close to the identity, number of optimization steps $N_{\mathrm{epoch}}$ per interval, and, if $\vec a$ depends on $\rho$, a fixed spatial grid and a number $N_\tau$ of field times
        \FOR{$k=0$ \TO $N_t-1$}
        \STATE freeze $\theta_0,\dots,\theta_{k-1}$
        \IF{$\vec a$ depends on $\rho$}
        \STATE choose $N_\tau$ times $t_k\le\tau_1<\dots<\tau_{N_\tau}\le t_{k+1}$
        \ENDIF
        \FOR{$n=0$ \TO $N_{\mathrm{epoch}}-1$}
        \IF{$\vec a$ depends on $\rho$}
        \STATE \emph{field update:} evaluate $\rho_\Theta(\tau_j,\cdot)$ on the grid for $j=1,\dots,N_\tau$ by \eqref{eq:composed_density}, with the current $\theta_k$ (moments by quadrature in velocity for kinetic models)
        \STATE solve the field equation on the grid at each $\tau_j$ and store the resulting velocity $\vec a^n$
        \ENDIF
        \STATE draw $N_{\mathrm{colloc}}$ points in $(t_k,t_{k+1})\times\Omega$ and $N_{\mathrm{ic}}$ points in $\Omega$
        \STATE \emph{flow update:} interpolate $\vec a^n$ in space and time at the collocation points (the prescribed $\vec a$ otherwise), without propagating derivatives with respect to $\theta_k$ through it
        \STATE update $\theta_k$ by one optimizer step on \eqref{eq:submap_loss}
        \ENDFOR
        \ENDFOR
        \STATE\textbf{Output:} $\Theta=(\theta_0,\dots,\theta_{N_t-1})$, the solution at any $(t,\vec x)$ being given by \eqref{eq:composed_density}
    \end{algorithmic}
\end{algorithm}

\paragraph{Choice of the submap partitioning.} Taking $N_t=1$ recovers the formulation
of \cref{sec:lagrangian_pinn}. Increasing $N_t$ makes each submap closer to
the identity, so that smaller networks
are sufficient on each interval, at the price of the sequential training and
of the growing evaluation cost. The partition need not be uniform: short
intervals can be used where the dynamics is fast and long ones where it is
slow, and the influence of $N_t$ at a fixed network size is studied in
\cref{sec:simpl}.
Note that the composition in \eqref{eq:comp_structure} is used similarly in grid-based methods~\cite{MercierYinNave2019,YinMercierYadavSchneiderNave2021,KrahYinBergmannNaveSchneider2024}, in which each submap is represented on a grid and a new submap is introduced once an error indicator exceeds a prescribed threshold. In contrast, the formulation above represents each submap with a neural network and prescribes the partition in advance.
\subsection{Coupling PINNs to the field equation}
\label{sec:field_coupling}

When the velocity field is determined by a field equation, the loss associated with any of the formulations in \cref{sec:eulerian_pinn,sec:lagrangian_pinn,sec:lagrangian_pinn_semigroup} must be evaluated carefully. We therefore freeze the velocity during each optimization step and update it between steps, resulting in a Picard iteration. Given the current network parameters, one iteration consists of the following stages.
\begin{enumerate}[nosep]
    \item \textbf{Field update.} The current approximation is evaluated on a fixed spatial grid and at $N_\tau$ times over the time interval. For the Eulerian PINN of \cref{sec:eulerian_pinn}, this evaluation is performed directly from the network, whereas for the Lagrangian PINNs of \cref{sec:lagrangian_pinn,sec:lagrangian_pinn_semigroup}, it is obtained through the pullback relations \eqref{eq:density_reconstruction} and~\eqref{eq:composed_density}. For kinetic models, the required moments are computed by quadrature over the velocity variables. The field equation is then solved on the grid, and the resulting velocity field is stored.
    \item \textbf{Flow update.} The stored velocity is interpolated in space and time at the collocation points and treated as a prescribed external field. No derivatives with respect to the network parameters are propagated through this field. An optimizer step applied to the loss of the formulation under consideration then updates the network parameters.
\end{enumerate}
The velocity is therefore updated once per optimization step, and a fixed point of the iteration satisfies both the corresponding loss and the field equation defining the velocity. If the velocity is prescribed rather than determined by a field equation, the first stage is omitted and the velocity is computed only once.

\subsection{Architectures of the neural networks approximating the flow maps}
\label{sec:flow_networks}

The submap networks of \eqref{eq:submap_ansatz} are required to map $\Omega$
into itself while respecting the periodicity of the domain, to be close to the
identity at the beginning of their interval, and to provide the Jacobian
determinants of \eqref{eq:composed_density} at a reasonable cost. We use two
families, described below for a generic submap $\submap_\theta(t,\vec x)$. In
each of them the time is passed to the network as a conditioning input
alongside $\vec x$, and physical parameters can be appended to it.

\begin{description}
    \item[Perturbation-of-the-identity networks.] The simplest choice is
          \begin{equation}
              \submap_\theta(t,\vec x) = \vec x + h\,\mathcal{N}_\theta(t,\vec x),
              \label{eq:mlp_flow}
          \end{equation}
          with $\mathcal{N}_\theta$ a multilayer perceptron (MLP)~\cite{HornikStinchcombeWhite1989} of $\ell$ hidden layers of width $w$
          and activation $\sigma$, acting on $\vec z=(t,\vec x)$,
          $$
              \mathcal{N}_\theta = \mathcal{A}_{\ell+1} \circ \mathcal{L}_{\ell} \circ \cdots \circ \mathcal{L}_1,
              \quad \mathcal{L}_i(\vec z) = \sigma\left( W_i\vec z+\vec b_i \right),
              \quad \mathcal{A}_{\ell+1}(\vec z) = W_{\ell+1}\vec z+\vec b_{\ell+1},
          $$
          whose output layer $\mathcal{A}_{\ell+1}$ is affine, and $h>0$ a scaling.
          The weights and biases are
          $W_1\in\R^{w\times n_0}$, $W_i\in\R^{w\times w}$ for $2\le i\le\ell$,
          $W_{\ell+1}\in\R^{d\times w}$, $\vec b_i\in\R^{w}$ for $1\le i\le\ell$ and
          $\vec b_{\ell+1}\in\R^{d}$, where $n_0$ is the dimension of the input $\vec z$
          ($1+d$, plus one for each periodic direction because of the embedding
          described below, plus the number of physical parameters).
          The identity condition of
          \eqref{eq:submap_loss} is met by letting $\mathcal{N}_\theta$ vanish, which a small
          initialization nearly satisfies. Periodicity is imposed by construction,
          each periodic component $x_j$ of period $L_j$ being replaced by
          $(\cos(2\pi x_j/L_j),\sin(2\pi x_j/L_j))$ before the first layer, so that
          no penalty term is needed for it.

    \item[Symplectic networks.] When the velocity is Hamiltonian, the phase space splits as
          $\vec x=(\vec q,\vec p)$ and the flow is symplectic. This structure is
          preserved exactly by networks composed of gradient
          modules~\cite{jin2020sympnets}, which write
          $$
              \submap_{\theta}(t,\cdot) = \phi^{\ell}_{\mathrm{up}} \circ  \phi^{\ell}_{\mathrm{down}} \circ \cdots \circ \phi^{1}_{\mathrm{up}} \circ  \phi^{1}_{\mathrm{down}},
          $$
          where the modules $\phi^{i}_{\mathrm{down}}$ and $\phi^{i}_{\mathrm{up}}$, $i=1,\dots,\ell$, read
          \begin{equation}
              \phi_{\mathrm{down}}^{i}: \,  (\vec q,\vec p) \mapsto \big(\vec q,\ \vec p + h\,\nabla_{\vec q} V_i(\vec q,t)\big),
              \qquad
              V_i(\vec q,t) = \vec s_i(t)^{\!\top}\,\Sigma\big(W_i\vec q + \vec b_i(t) \big),
              \label{eq:sympnet_layer_down}
          \end{equation}
          \begin{equation}
              \phi_{\mathrm{up}}^{i}: \,  (\vec q,\vec p) \mapsto \big(\vec q+ h\,\nabla_{\vec p} T_i(\vec p,t),\ \vec p \big),
              \qquad
              T_i(\vec p,t) = \widehat{\vec s}_i(t)^{\!\top}\,\Sigma\big(\widehat W_i\vec p + \widehat{\vec b}_i(t) \big).
              \label{eq:sympnet_layer_up}
          \end{equation}
          Here $\vec q,\vec p\in\R^{n}$ with $d=2n$, the weights are
          $W_i,\widehat W_i\in\R^{w\times n}$, the vectors
          $\vec s_i,\vec b_i,\vec c_i$ and their counterparts in $T_i$ belong to $\R^{w}$,
          $w$ being the width of the modules,
          and $\Sigma$ is an antiderivative of the activation $\sigma$, applied componentwise,
          so that $\nabla_{\vec q} V_i = W_i^{\top}\big(\vec s_i(t)\odot\sigma(W_i\vec q + \vec b_i(t))\big)$,
          and the time enters through $\vec s_i(t) = (t-t_k)\,\vec s_i$ and
          $\vec b_i(t) = \vec b_i + (t-t_k)\,\vec c_i$ (and similarly for $T_i$).
          Since $\vec s_i(t_k)=\widehat{\vec s}_i(t_k)=0$, every module, hence the
          network, is exactly the identity at $t=t_k$.
          The Jacobian of \eqref{eq:sympnet_layer_down} is
          $\left(\begin{smallmatrix} \mathrm{I} & 0\\ h\nabla^2_{\vec q}V_i & \mathrm{I}\end{smallmatrix}\right)$,
          which is symplectic since $\nabla^2_{\vec q}V_i$ is symmetric, and
          has determinant one. Therefore, one has
          $\det\nabla\submap_\theta\equiv 1$, so the products of determinants in
          \eqref{eq:composed_density} disappear and, as noted in \cref{sec:flowmaps}
          for a divergence-free velocity, every $L^p$ norm of $\rho_\theta$ is
          conserved exactly.

          When some directions of $\Omega$ are periodic, the network needs to be adapted accordingly. We explain here a periodic variant of the symplectic network that was introduced in~\cite{fournet2026nonlinearcontrolvariatedeltaf}. In this case, the potential $V_i$ is replaced by
          $$
              \widetilde{V}_{i}(\vec q,t) = \vec s_i(t)^{\!\top} \Sigma \left( W^{c}_i \cos\left( \frac{2\pi \vec q}{L} \right) + \vec b^{c}_i(t) \right) + \vec s_i(t)^{\!\top} \Sigma \left( W^{s}_i \sin\left( \frac{2\pi \vec q}{L} \right) + \vec b^{s}_i(t) \right),
          $$ with $W^{c}_i,W^{s}_i\in\R^{w\times n}$ and
          $\vec b^{c}_i(t),\vec b^{s}_i(t)\in\R^{w}$ depending on $t$ as $\vec b_i(t)$ above, so that the
          module maps the torus into itself, and the network is periodic by construction.

\end{description}

Both network families provide the same two operations, evaluation and
Jacobian determinant, so the residual \eqref{eq:map_transport} and the
loss functions \eqref{eq:flow_loss} and \eqref{eq:submap_loss} are written once and for all.
Changing the family changes the approximation space of the submaps and nothing
else in the training procedure.

\section{Results}
\label{sec:results}
This section presents the numerical results obtained with the Eulerian PINNs and Lagrangian PINNs detailed in \cref{sec:pinns}.
We start with simple advection-diffusion test cases in \cref{subsec:advect-diff-eq}.
Thereafter, we consider physical systems such as the Euler equations (\cref{sec:euler}) and kinetic equations (\cref{sec:vlasov_poisson,sec:gyrokinetic-test-case}).

\subsection{Advection-diffusion equations}
\label{subsec:advect-diff-eq}
The following test cases investigate the performance of the methodology with a given external velocity field. In \cref{sec:simpl} this velocity field and the corresponding analytic map are known and the error on the map can be directly assessed.
In \cref{sec:simpl-diffusion} a small diffusion is added to the same problem. The density is then no longer constant along the characteristics, and the flow map is combined with an Eulerian field that carries the diffusion.

\subsubsection{Simple advection}
\label{sec:simpl}

\begin{figure}[htp!]
    \centering
    \includegraphics[width=0.9\linewidth]{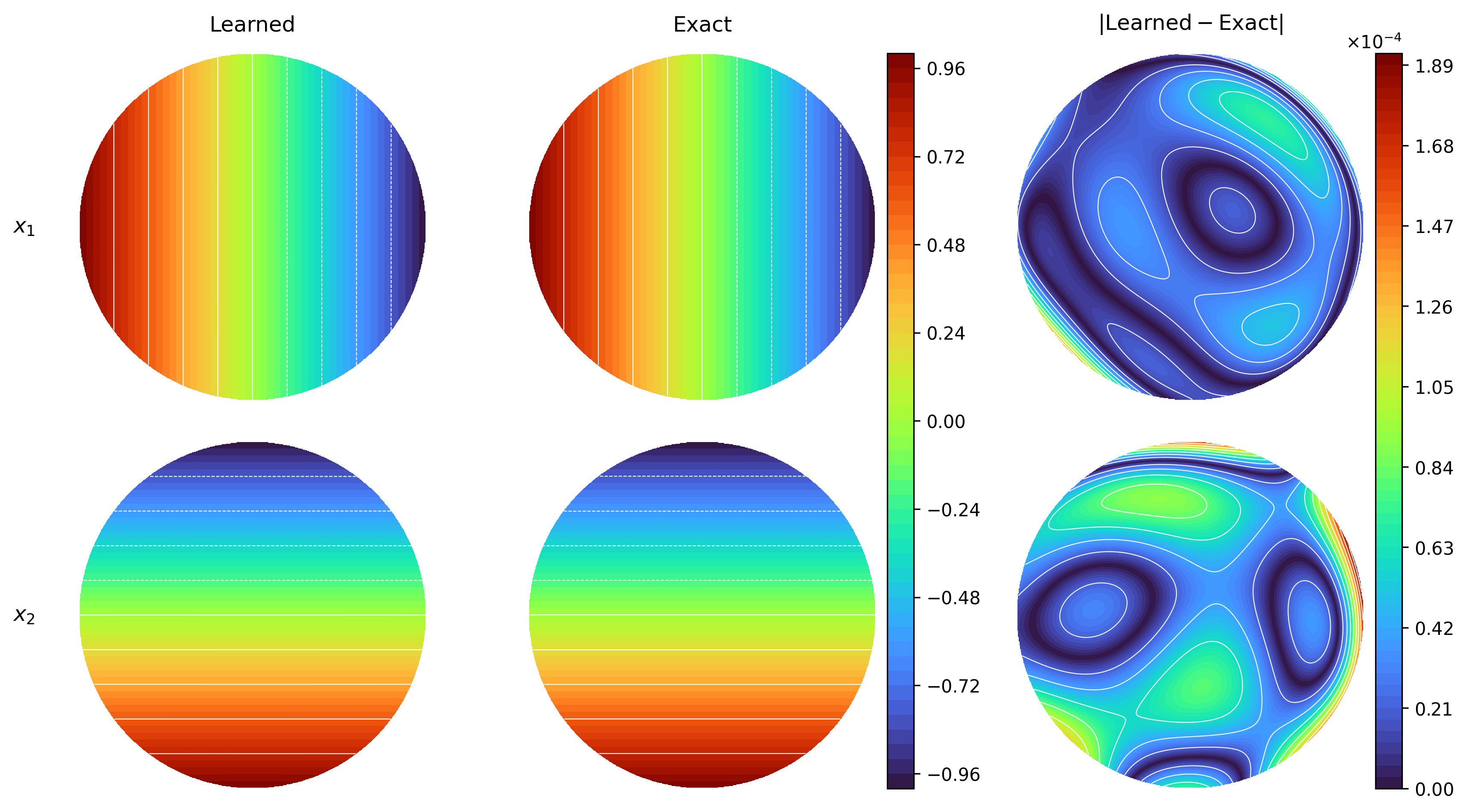}
    \caption{Learned and exact backward flow map $\submap_0^{0.5}$ at $t=0.5$ (end of the first of $N_t=2$ subintervals), one row per component ($(\submap_0^{0.5})_1$, $(\submap_0^{0.5})_2$), plotted over the unit disk in the $(x_1,x_2)$ plane: learned (left), exact (middle) and pointwise difference $|\tilde{\submap}_0^{0.5}-\submap_0^{0.5}|$ (right). The learned and exact columns share one color scale; the difference column, of order $10^{-4}$, has its own.}
    \label{fig:map-comparison}
\end{figure}

\begin{figure}[htp!]
    \centering
    \includegraphics[width=\linewidth]{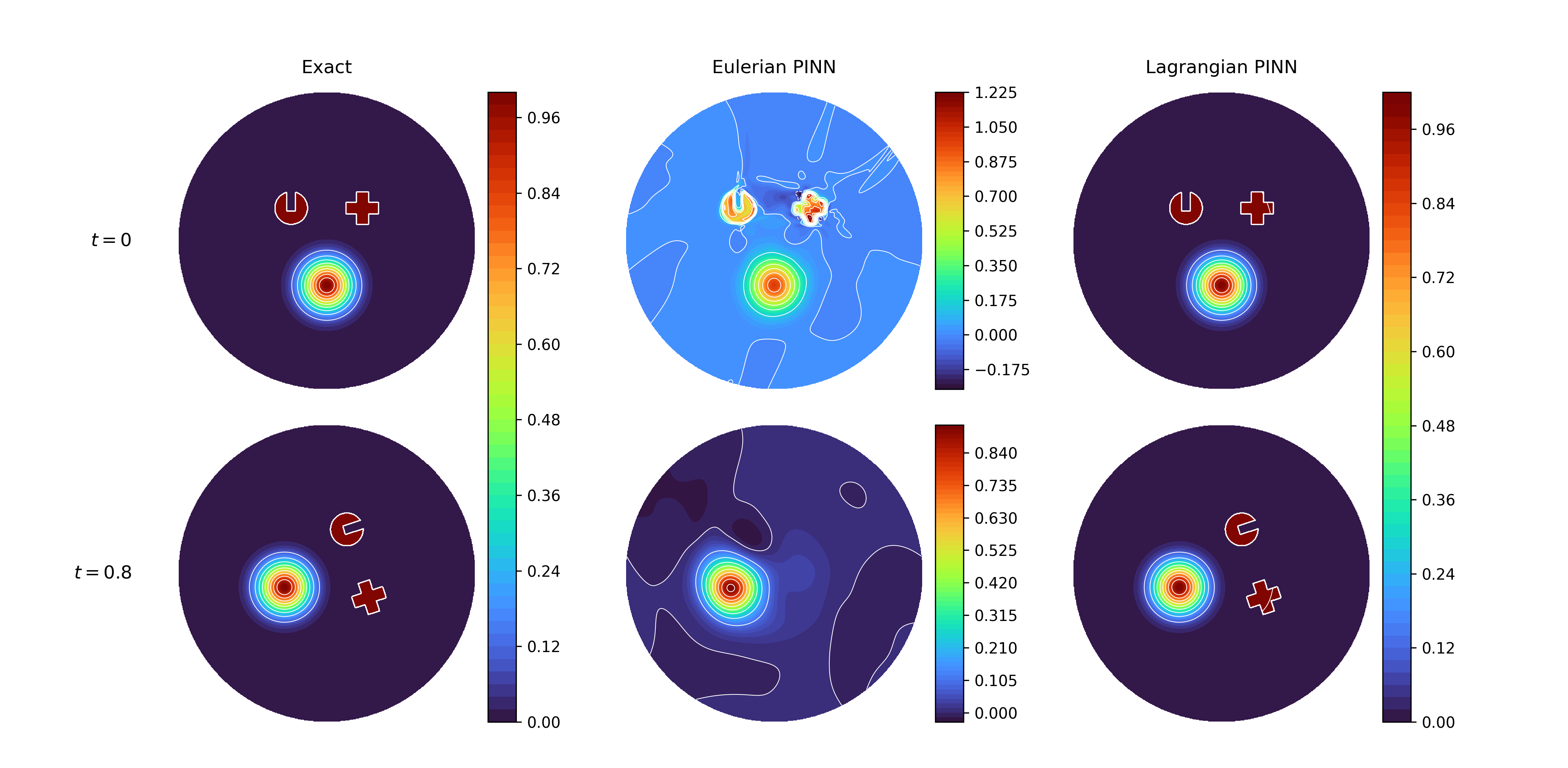}
    \caption{Exact solution (left), Eulerian PINN (middle) and Lagrangian PINN (right) for the simple advection test case at $t=0$ and $t=0.8$.
        Exact and Lagrangian PINN share the same color scale, the Eulerian PINN uses its own per-row scale, since it preserves neither the density range nor the periodicity of the exact solution.}
    \label{fig:comparison-smiley}
\end{figure}

We consider the 2D linear transport equation
\begin{equation}
    \partial_t \rho(t,\vec{x}) + \vec{a}(\vec{x})\cdot\nabla\rho(t,\vec{x}) = 0,
    \qquad \vec{x}=(x_1,x_2),
\end{equation}
with the rigid-rotation velocity field $\vec{a}(x_1,x_2) = (-2\pi x_2,\; 2\pi x_1)$  in the unit disk, up to the final time $T=1$, which corresponds to one full rotation.
The initial density is a composite ``smiley-face'' shape, the pointwise
maximum of three components: a Gaussian mouth and two eyes,
\begin{equation}
    \label{eq:initial-density-smiley}
    \rho_0(\vec x) = \max\big(\rho_{\text{mouth}}(\vec x),\, \rho_{\text{cross}}(\vec x),\, \rho_{\text{slot}}(\vec x)\big),
\end{equation}
where $\rho_{\text{mouth}}$ is a Gaussian bump of width $0.11$ centered at
$(0,-0.30)$, $\rho_{\text{cross}}$ is the indicator function of a plus-shaped
region (two perpendicular bars of half-length $0.11$ and half-width $0.04$)
centered at $(0.24,0.22)$, and $\rho_{\text{slot}}$ is the indicator function
of a disk of radius $0.11$ centered at $(-0.24,0.22)$ with a rectangular
notch removed near its top (a slotted disk). This composite shape, with its mix of smooth and
sharp features, is used throughout this section.
The exact flow map used for evaluation is
\begin{equation}
    \label{eq:exact-flow-map-advection}
    \submap_0^{t}(x_1,x_2)
    =
    \begin{pmatrix}
        \cos(2\pi t)  & \sin(2\pi t) \\
        -\sin(2\pi t) & \cos(2\pi t)
    \end{pmatrix}
    \begin{pmatrix}x_1\\x_2\end{pmatrix}.
\end{equation}
The exact solution is the composition of the initial density \eqref{eq:initial-density-smiley} with the flow map \eqref{eq:exact-flow-map-advection}: $\rho(t,\vec{x})=\rho_0(\submap_0^t(\vec{x}))$.
The exact flow map and the resulting solution are shown in \cref{fig:map-comparison,fig:comparison-smiley}.
To compare the Eulerian and Lagrangian approaches directly on this
initial condition, we train an Eulerian PINN
with a single time-continuous MLP with $4$ hidden
layers of width $20$ and $\tanh$ activation, minimizing \eqref{eq:eulerian_loss}
directly over the full time interval $[0,T]$ and compare it against the Lagrangian
PINN with $N_t=2$ submap networks of depth $\ell=5$ and hidden width $w=14$.
Both models use
$N_{\text{colloc}}=8000$ and $N_{\text{ic}}=4000$. The Lagrangian PINN
trains each of its $N_t=2$ flow segments for $100$ epochs, and the Eulerian PINN is trained
for $400$ epochs.

\cref{fig:comparison-smiley} shows the exact solution together with both
predictions at $t=0$ and $t=0.8$. The Lagrangian PINN matches the exact
solution closely at both times: since it only evaluates the flow map
learned on $[0,T]$ and reconstructs the density through the pullback
relation \eqref{eq:density_reconstruction}, the sharp contours of the
initial shape are transported to their rotated position at $t=0.8$
without being smeared.
The Eulerian PINN, despite the larger training budget, still struggles to reproduce the sharp contours of the initial condition,
and at $t=0.8$ the two eyes are lost altogether, only the smooth mouth being transported.

To study architecture and temporal decomposition, we run a grid of experiments with:
\begin{itemize}[nosep]
    \item network depths (hidden-layer counts): $\ell\in\{2,3,4,6,8\}$,
    \item number of time subintervals/submaps: $N_t\in\{1,2,4\}$,
    \item collocation parameters: $N_{\text{colloc}}=8000$, $N_{\text{ic}}=4000$,
    \item epochs: $100$, hidden width: $14$, activation: $\tanh$.
\end{itemize}
As a baseline, we also sweep the Eulerian PINN of \cref{sec:eulerian_pinn}
over the same depths (hidden width $20$, $\tanh$ activation), trained over
the full time interval $[0,T]$ rather than a single subinterval, since there is
no flow decomposition to exploit, and for $400$ epochs
as above, so that it receives the same total training budget as the most
expensive Lagrangian configuration in the sweep ($N_t=4$).
As error metric, we use the mean pointwise error with respect to the exact
solution, computed for both the backward map $\submap_0^{t}$ (in Euclidean norm)
and the density $\rho$ (in absolute value). The pointwise error is averaged over
a uniform evaluation grid restricted to the unit disk, and then over the end
times $t_k$ of the subintervals. For the Eulerian PINN, which has no
subintervals, it is averaged over $t\in\{0.25,0.5,0.75,1\}$.
Both errors are shown in \cref{fig:arch-error-vs-depth-advection}.
We see that the flow-map approach is more than one order of magnitude more precise than the Eulerian one, independently of the architecture.
Furthermore, we note that the one-map approach ($N_t=1$), similar to \cite{MontiMorozovskaPoggiShukla2026}, is less accurate than the approach using submaps ($N_t>1$).

\begin{figure}[htp!]
    \centering
    \begin{subfigure}[t]{0.49\linewidth}
        \centering
        \setlength{\figurewidth}{\linewidth}
        \setlength{\figureheight}{0.8\linewidth}
        %
    \definecolor{mplcolor0}{rgb}{0.121570,0.466670,0.705880}
\definecolor{mplcolor1}{rgb}{0.172550,0.627450,0.172550}
\begin{tikzpicture}
\begin{axis}[
  width=\figurewidth,
  height=\figureheight,
  ymode=log, log basis y=10,
  xlabel={Network depth (number of hidden layers)},
  ylabel={Mean error},
  xmin=1.7, xmax=8.3,
  ymin=5e-06, ymax=0.1,
  xmajorgrids=true,
  ymajorgrids=true,
  legend to name=archlegend,
  legend columns=2,
  legend style={font=\footnotesize, legend cell align={left}, column sep=3pt}
]
\addplot [color=mplcolor0, mark=*, mark size=2.5pt, mark options={solid}, line width=1.00pt]
  coordinates {(2,0.00129892) (3,0.000208129) (4,0.000271809) (6,0.000159864) (8,0.000935826)};
\addplot [color=mplcolor0, dashed, mark=o, mark size=2.5pt, mark options={solid, fill=white}, line width=1.00pt]
  coordinates {(2,0.00129377) (3,0.000118043) (4,0.000139851) (6,4.65634e-05) (8,0.000207218)};
\addplot [color=mplcolor0, mark=square*, mark size=2.5pt, mark options={solid}, line width=1.00pt]
  coordinates {(2,0.000331907) (3,0.000140464) (4,7.78328e-05) (6,4.84438e-05) (8,0.000206836)};
\addplot [color=mplcolor0, dashed, mark=square, mark size=2.5pt, mark options={solid, fill=white}, line width=1.00pt]
  coordinates {(2,0.000479466) (3,1.64215e-05) (4,1.96593e-05) (6,8.56795e-06) (8,3.55276e-05)};
\addplot [color=mplcolor0, mark=triangle*, mark size=2.5pt, mark options={solid}, line width=1.00pt]
  coordinates {(2,0.000138614) (3,8.63881e-05) (4,8.6724e-05) (6,6.61665e-05) (8,9.68466e-05)};
\addplot [color=mplcolor0, dashed, mark=triangle, mark size=2.5pt, mark options={solid, fill=white}, line width=1.00pt]
  coordinates {(2,2.09005e-05) (3,2.04239e-05) (4,2.26607e-05) (6,1.77926e-05) (8,2.30102e-05)};
\addplot [color=mplcolor1, dashed, mark=diamond*, mark size=2.5pt, mark options={solid}, line width=1.00pt]
  coordinates {(2,0.0621648) (3,0.0367366) (4,0.0295108) (6,0.0255319) (8,0.0279254)};
\addlegendentry{Lagrangian, flow map, $N_t=1$}
\addlegendentry{Lagrangian, density, $N_t=1$}
\addlegendentry{Lagrangian, flow map, $N_t=2$}
\addlegendentry{Lagrangian, density, $N_t=2$}
\addlegendentry{Lagrangian, flow map, $N_t=4$}
\addlegendentry{Lagrangian, density, $N_t=4$}
\addlegendentry{Eulerian, density}
\end{axis}
\end{tikzpicture}%

        \caption{Simple advection (\cref{sec:simpl}).}
        \label{fig:arch-error-vs-depth-advection}
    \end{subfigure}
    \hfill
    \begin{subfigure}[t]{0.49\linewidth}
        \centering
        \setlength{\figurewidth}{\linewidth}
        \setlength{\figureheight}{0.8\linewidth}
        %
    \definecolor{mplcolor0}{rgb}{0.121570,0.466670,0.705880}
\definecolor{mplcolor1}{rgb}{0.172550,0.627450,0.172550}
\begin{tikzpicture}
\begin{axis}[
  width=\figurewidth,
  height=\figureheight,
  ymode=log, log basis y=10,
  xlabel={Network depth (number of hidden layers)},
  ylabel={Mean error},
  xmin=1.7, xmax=8.3,
  ymin=5e-06, ymax=0.1,
  xmajorgrids=true,
  ymajorgrids=true,
]
\addplot [color=mplcolor0, mark=*, mark size=2.5pt, mark options={solid}, line width=1.00pt]
  coordinates {(2,0.00634612) (3,0.00268289) (4,0.00243823) (6,0.00424842) (8,0.00183019)};
\addplot [color=mplcolor0, dashed, mark=o, mark size=2.5pt, mark options={solid, fill=white}, line width=1.00pt]
  coordinates {(2,0.0165633) (3,0.00706029) (4,0.00401787) (6,0.00382138) (8,0.00533295)};
\addplot [color=mplcolor0, mark=square*, mark size=2.5pt, mark options={solid}, line width=1.00pt]
  coordinates {(2,0.00513464) (3,0.00117019) (4,0.00106341) (6,0.00110299) (8,0.000316885)};
\addplot [color=mplcolor0, dashed, mark=square, mark size=2.5pt, mark options={solid, fill=white}, line width=1.00pt]
  coordinates {(2,0.0132424) (3,0.00467103) (4,0.00398055) (6,0.00188585) (8,0.00191106)};
\addplot [color=mplcolor0, mark=triangle*, mark size=2.5pt, mark options={solid}, line width=1.00pt]
  coordinates {(2,0.0028444) (3,0.00071386) (4,0.000938823) (6,0.00135168) (8,0.000716451)};
\addplot [color=mplcolor0, dashed, mark=triangle, mark size=2.5pt, mark options={solid, fill=white}, line width=1.00pt]
  coordinates {(2,0.00972677) (3,0.00421828) (4,0.00602362) (6,0.00269331) (8,0.00159581)};
\addplot [color=mplcolor1, dashed, mark=diamond*, mark size=2.5pt, mark options={solid}, line width=1.00pt]
  coordinates {(2,0.0465756) (3,0.0298439) (4,0.0198403) (6,0.0209679) (8,0.0211014)};
\end{axis}
\end{tikzpicture}%

        \caption{Advection with diffusion (\cref{sec:simpl-diffusion}).}
        \label{fig:arch-error-vs-depth-diffusion}
    \end{subfigure}
    \\[4pt]
    \input{figures/advection/architecture_study_colors}\ref*{archlegend}
    \caption{Mean flow-map error (solid lines, filled markers) and mean density
        error (dashed lines, open markers) versus network depth $\ell$, for the Lagrangian PINN
        (blue, marker shape encodes $N_t$) and the Eulerian PINN baseline (green
        diamonds, density error only).}
    \label{fig:arch-error-vs-depth}
\end{figure}

\subsubsection{Simple advection with diffusion}
\label{sec:simpl-diffusion}

We now add a small diffusion to the transport equation of \cref{sec:simpl},
\begin{equation}
    \label{eq:advection-diffusion}
    \partial_t \rho(t,\vec{x}) + \vec{a}(\vec{x})\cdot\nabla\rho(t,\vec{x}) = \nu\,\Delta\rho(t,\vec{x}),
    \qquad \nu = 10^{-3},
\end{equation}
with the same rigid rotation $\vec a$, the same unit disk and the same final time $T=1$.
Over $[0,T]$ the diffusion length $\sqrt{2\nu T}\approx 0.045$ remains
below the size of the eyes of the smiley, so that the problem is advection-dominated.
To have a closed-form reference, the initial density is the smiley of
\cref{sec:simpl} (with the three components summed instead of taking their maximum, which differs by less than
$10^{-4}$ since the components are almost disjoint) convolved with a Gaussian of
standard deviation $\sigma_0=0.005$, so its contours remain sharp at the
scale of the plots.
Since the rotation commutes with the Laplacian, the exact solution is the
initial shape convolved with a Gaussian of growing width and rotated,
\begin{equation}
    \label{eq:exact-advection-diffusion}
    \rho(t,\vec x) = w\big(t,\submap_0^t(\vec x)\big),
    \qquad
    w(t,\cdot) = G_{\sigma(t)} * \rho_{\text{smiley}},
    \qquad
    \sigma(t)^2 = \sigma_0^2 + 2\nu t,
\end{equation}
where $\submap_0^t$ is the rotation \eqref{eq:exact-flow-map-advection},
$G_\sigma$ the centered Gaussian of standard deviation $\sigma$, and
$\rho_{\text{smiley}}$ the sum of the three components of \eqref{eq:initial-density-smiley}.
The density stays negligible near the boundary of the disk, so no boundary condition is needed.

\paragraph{Flow map and co-moving field.}
With diffusion, the density is no longer constant along the characteristics,
and the reconstruction \eqref{eq:density_reconstruction} from $\rho_0$ alone no longer applies.
Instead, the flow map is kept for the advection only, and the remaining
evolution is carried by an Eulerian field written in the co-moving
coordinates $\vec y = \submap_0^t(\vec x)$,
\begin{equation}
    \label{eq:split-ansatz}
    \rho_{\theta,\eta}(t,\vec x) = w_\eta\big(t,\submap_\theta(t,\vec x)\big),
    \qquad
    \submap_\theta(t,\vec x) = \vec x + t\,\mathcal N_\theta(t,\vec x),
\end{equation}
where $\mathcal N_\theta$ and $w_\eta$ are two MLPs.
The identity condition of the map holds by construction and $w_\eta(0,\cdot)$
is fitted to $\rho_0$.
The map takes care of the transport. The field $w_\eta$ therefore only has to
describe the slow diffusion, which does not move with the rotation and
changes smoothly in time.
In model order reduction, methods that move to a co-moving frame are
reviewed in \cite{HesthavenPeherstorferUnger2026}. A prominent variant is
the shifted POD \cite{ReissSchulzeSesterhennMehrmann2018,Krah_P_2024_PP_robust_spodpmdfmt},
in which the data is first transformed into a co-moving frame where it is
nearly stationary and therefore simpler to approximate with linear or
nonlinear methods.
In our test case the flow map is a rotation, which does not change
distances. The field $w$ then solves the plain heat equation
$\partial_t w = \nu\Delta w$. For a general velocity field, $w$ solves a
diffusion equation whose diffusion depends on the direction.
The training does not use any of this. Both networks are trained jointly,
with a single optimizer, on the transport residual of the map
(the first term of \eqref{eq:flow_loss}), on the residual of
\eqref{eq:advection-diffusion} evaluated on the composition
$\rho_{\theta,\eta}$, and on the initial condition $w_\eta(0,\cdot)=\rho_0$,
\begin{equation}
    \label{eq:split-loss}
    \begin{split}
        \min_{\theta,\eta} \Bigg(
         & \int_{0}^{T}\!\!\int_{\Omega}
        \left\|\partial_t \submap_\theta
        + \big(\vec a\cdot \nabla_{\vec x}\big)\,\submap_\theta\right\|^2
        + \left|\partial_t \rho_{\theta,\eta} + \vec a\cdot\nabla\rho_{\theta,\eta} - \nu\Delta\rho_{\theta,\eta}\right|^2\,\dx\,\dt \\
         & + \lambda \int_{\Omega}\left|w_\eta(0,\vec x) - \rho_0(\vec x)\right|^2\,\dx
        \Bigg).
    \end{split}
\end{equation}
For $N_t>1$ subintervals, the same ansatz is used on each interval
$[t_k,t_{k+1}]$ with its own pair of networks: $\submap_{\theta_k}$ maps
back to $t_k$ as in \eqref{eq:submap_ansatz}, the co-moving field
$w_{\eta_k}$ is written in the coordinates of time $t_k$, and its initial
condition is the density reached at the end of the previous interval,
$w_{\eta_k}(t_k,\cdot) = \rho_{\theta_{k-1},\eta_{k-1}}(t_k,\cdot)$.
The intervals are trained one after the other, as in
\cref{sec:lagrangian_pinn_semigroup}.

We run the same parameter study as in \cref{sec:simpl}, with network depths
$\ell\in\{2,3,4,6,8\}$ and $N_t\in\{1,2,4\}$ subintervals, $100$ epochs per
subinterval for the Lagrangian PINN and $400$ epochs for the Eulerian PINN of
\cref{sec:eulerian_pinn}, which is trained on the same residual with the same
initial-condition weight.
The depth $\ell$ applies to all networks. The maps have hidden width $16$, the
co-moving field and the Eulerian PINN have hidden width $20$, all with
$\tanh$ activations, and we use $N_{\text{colloc}}=N_{\text{ic}}=10\,000$ and
$\lambda=8$. The errors are measured as in \cref{sec:simpl} and shown in
\cref{fig:arch-error-vs-depth-diffusion}.
For all architectures the Lagrangian PINN is more accurate than the Eulerian
one. Its mean density error ranges from $1.6\cdot10^{-3}$ to
$1.7\cdot10^{-2}$, while the Eulerian PINN does not get below $2.0\cdot10^{-2}$.
The best configurations are thus about one order of magnitude apart.
As for pure advection, subintervals help. With $N_t>1$ the flow-map error drops
to $3\cdot10^{-4}$--$1.4\cdot10^{-3}$ for $\ell\geq3$, and the best density errors,
below $2\cdot10^{-3}$, are all obtained with $N_t=2$ or $N_t=4$ and deep networks.
The density errors are larger than without diffusion. There, the density is
reconstructed exactly from $\rho_0$, whereas here the co-moving field also has
to learn the sharp initial condition.

\begin{figure}[htp!]
    \centering
    \includegraphics[width=\linewidth]{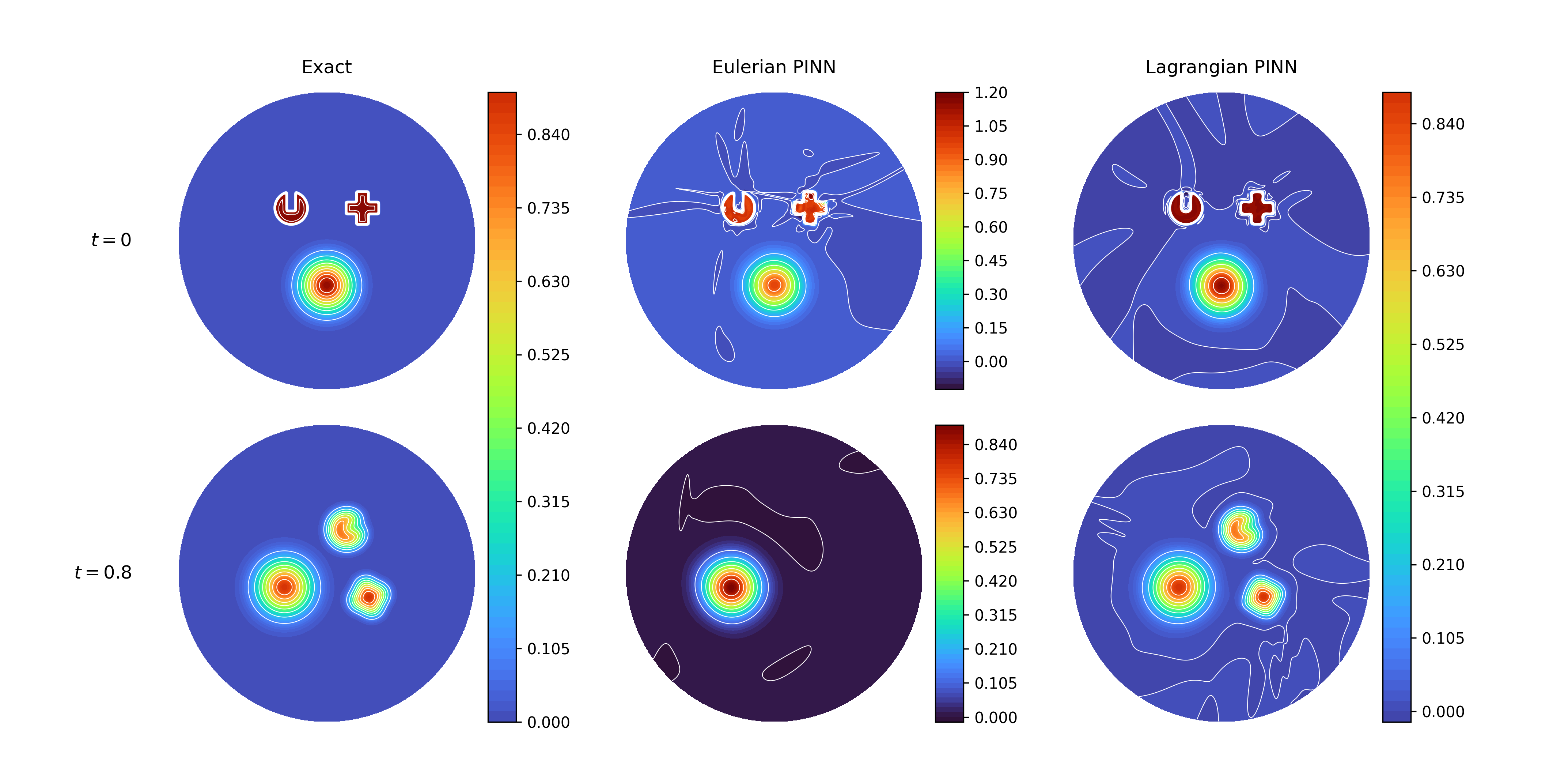}
    \caption{Exact solution (left), Eulerian PINN (middle) and Lagrangian PINN (right) for the advection-diffusion test case at $t=0$ and $t=0.8$,
        for the best networks of the parameter study (Lagrangian PINN: $N_t=4$, $\ell=8$, Eulerian PINN: $\ell=4$).
        Exact and Lagrangian PINN share the same color scale, the Eulerian PINN uses its own per-row scale.}
    \label{fig:comparison-smiley-diffusion}
\end{figure}

\cref{fig:comparison-smiley-diffusion} compares the best networks of the
study with the exact solution at $t=0$ and $t=0.8$. These are the Lagrangian
PINN with $N_t=4$ and $\ell=8$ and the Eulerian PINN with $\ell=4$.
The Lagrangian PINN reproduces the rotated and diffused smiley, including the
notch of the slotted disk and the arms of the cross, which have started to
blur by $t=0.8$. Its relative $L^2$ error is $6.0\cdot10^{-2}$ at $t=0$ and
$2.8\cdot10^{-2}$ at $t=0.8$. It does not grow along the rotation, because
the transport is handled by the map. The error mostly comes from fitting the
sharp initial condition, and it decreases as diffusion smooths the solution.
The remaining error consists of low-amplitude oscillations of the background,
visible in the iso-contours, and small defects along the sharp edges of the eyes.
The Eulerian PINN, in contrast, fits the initial condition only roughly and
loses the eyes altogether as soon as they start to move. At $t=0.8$ only the
smooth mouth is transported, and its relative $L^2$ error is $0.68$.

\subsection{Incompressible Euler equations}
\label{sec:euler}
We study the vorticity-streamfunction formulation of the incompressible Euler equations on a  periodic box:
\begin{align}
     & \partial_t \omega + \vec{u} \cdot \nabla \omega = 0                               \\
     & \nabla^2 \psi = \omega, \qquad \vec{u} = (-\partial_y \psi, \partial_x \psi)  \,.
\end{align}
For both test cases below, we monitor the conservation of enstrophy $ \|\omega(t,\cdot)\|_{L^2(\Omega)}^2$ and (kinetic) energy $ \|\vec u(t,\cdot)\|_{L^2(\Omega)}^2$,
\begin{equation}
    \|\omega(t,\cdot)\|_{L^2(\Omega)}^2 = \int_\Omega \omega(t,\vec x)^2 \,\dx ,
    \quad
    \|\vec u(t,\cdot)\|_{L^2(\Omega)}^2 = \int_\Omega |\vec u(t,\vec x)|^2 \,\dx ,
    \label{eq:enstrophy_energy_def}
\end{equation}
following the convention of eqs.~(3.32c)--(3.32d) in \cite{YinMercierYadavSchneiderNave2021}; both quantities are exactly conserved by the continuous Euler flow, so the corresponding conservation error at time $t$, $\operatorname{Enstrophy}(t) - \operatorname{Enstrophy}(0)$ (resp.\ $\operatorname{Energy}(t) - \operatorname{Energy}(0)$), quantifies the numerical/approximation error of a given method.

Both test cases below use the same networks and training parameters, run until $T=8$.
The Lagrangian PINN uses the semigroup formulation of \cref{sec:lagrangian_pinn_semigroup} with
\begin{itemize}[nosep]
    \item Number of time intervals: $N_t = 8$, of equal width $1$;
    \item Submap network: perturbation of the identity \eqref{eq:mlp_flow}, MLP with $3$ hidden layers of $32$ neurons and $\tanh$ activation, periodic in $x$ and $y$ by construction ($N_\theta = 2370$ per submap);
    \item Number of collocation points: $N_{\text{colloc}} = 10000$, $N_{\text{ic}} = 4000$;
    \item Number of training steps: $500$ natural-gradient steps per submap (regularization $10^{-6}$, logarithmic grid line search).
\end{itemize}
The single-map variant ($N_t=1$) uses the same network over $[0,T]$, trained for $4000$ steps so that it receives the same total budget.
The velocity is coupled through the Picard iteration of \cref{sec:field_coupling}: at every optimization step the vorticity is evaluated on a $128\times128$ grid at $N_\tau = 80$ times per subinterval, the Biot-Savart law is solved by FFT on this grid, and the resulting velocity is interpolated bilinearly in space and linearly in time at the collocation points.
The vorticity is reconstructed without the Jacobian determinant of the flow, $\omega(t,\cdot) = \omega_0\circ\submap_0^{t}$, both in this field update and in all results shown below. For the divergence-free Euler velocity, the determinant is one, but the MLP submaps do not preserve volume exactly, so including the determinant of the learned map would only add an error to the reconstructed vorticity.
The Eulerian PINN of \cref{sec:eulerian_pinn} is a single space-time MLP approximating $\omega(t,\vec x)$ directly, with $3$ hidden layers of $64$ neurons, $\tanh$ activation and the same periodic embedding ($N_\theta = 8769$). It is trained for $1200$ natural-gradient steps with $N_{\text{colloc}} = 10000$ and $N_{\text{ic}} = 4000$, using the same Picard iteration for the velocity. For the vortex merger, its regularization is increased to $10^{-4}$ and its collocation times are drawn as $t = T u^2$ with $u$ uniform in $[0,1]$, concentrating them near $t=0$; without these changes the Eulerian PINN collapses to $\omega = 0$ for $t>0$.
All PINNs are trained on a single NVIDIA V100 GPU.

\subsubsection{Four-modes test case}
\label{sec:euler_four_modes}

In this section we consider the ``four-modes'' test case of
\cite{YinMercierYadavSchneiderNave2021}, with initial vorticity
\begin{equation}
    \omega_0(x,y) = \cos(x) + \cos(y) + 0.6 \cos(2x) + 0.2 \cos(3x) ,
    \label{eq:euler_four_modes_ic}
\end{equation} inside the domain $[0,2\pi]^2$ with periodic boundary conditions.
The flow can be characterized as two vortices of opposite sign, which roll up into increasingly thin vortex sheets as they interact.

\begin{figure}[htp!]
    \centering
    \setlength{\figurewidth}{1\linewidth}
    \setlength{\figureheight}{0.2206\linewidth}
    %
    {\pgfplotsset{compat=1.18}
\pgfplotsset{colormap={turbo}{rgb255(0.0000)=(48,18,59) rgb255(0.0159)=(53,30,88) rgb255(0.0317)=(57,42,115) rgb255(0.0476)=(61,53,139) rgb255(0.0635)=(64,64,162) rgb255(0.0794)=(66,75,181) rgb255(0.0952)=(68,86,199) rgb255(0.1111)=(70,97,214) rgb255(0.1270)=(70,107,227) rgb255(0.1429)=(71,118,238) rgb255(0.1587)=(70,128,246) rgb255(0.1746)=(69,138,252) rgb255(0.1905)=(66,148,255) rgb255(0.2063)=(61,158,254) rgb255(0.2222)=(55,168,250) rgb255(0.2381)=(47,178,244) rgb255(0.2540)=(39,190,233) rgb255(0.2698)=(32,199,223) rgb255(0.2857)=(27,208,213) rgb255(0.3016)=(24,215,202) rgb255(0.3175)=(24,222,192) rgb255(0.3333)=(26,228,182) rgb255(0.3492)=(32,234,172) rgb255(0.3651)=(42,239,161) rgb255(0.3810)=(53,243,148) rgb255(0.3968)=(67,247,135) rgb255(0.4127)=(82,250,122) rgb255(0.4286)=(97,252,108) rgb255(0.4444)=(113,254,95) rgb255(0.4603)=(128,255,83) rgb255(0.4762)=(143,255,73) rgb255(0.4921)=(156,254,64) rgb255(0.5079)=(169,251,57) rgb255(0.5238)=(180,248,54) rgb255(0.5397)=(190,244,52) rgb255(0.5556)=(200,239,52) rgb255(0.5714)=(210,233,53) rgb255(0.5873)=(219,226,54) rgb255(0.6032)=(227,219,56) rgb255(0.6190)=(235,211,57) rgb255(0.6349)=(241,203,58) rgb255(0.6508)=(246,195,58) rgb255(0.6667)=(250,186,57) rgb255(0.6825)=(252,177,54) rgb255(0.6984)=(254,167,50) rgb255(0.7143)=(254,155,45) rgb255(0.7302)=(254,144,41) rgb255(0.7460)=(252,132,35) rgb255(0.7619)=(249,117,29) rgb255(0.7778)=(245,105,24) rgb255(0.7937)=(241,93,19) rgb255(0.8095)=(236,83,15) rgb255(0.8254)=(231,73,12) rgb255(0.8413)=(225,65,9) rgb255(0.8571)=(218,57,7) rgb255(0.8730)=(210,49,5) rgb255(0.8889)=(202,42,4) rgb255(0.9048)=(193,35,2) rgb255(0.9206)=(183,29,2) rgb255(0.9365)=(172,23,1) rgb255(0.9524)=(161,18,1) rgb255(0.9683)=(149,13,1) rgb255(0.9841)=(136,8,2) rgb255(1.0000)=(122,4,3)}}
\begin{tikzpicture}
\begin{axis}[
  xlabel={$x$},
  ylabel={$y$},
  at={(0.0172\figurewidth, 0.0779\figureheight)},
  anchor=south west,
  width=0.1834\figurewidth,
  height=0.8312\figureheight,
  title={Spectral},
  xmin=0, xmax=0.992188,
  ymin=0, ymax=0.992188,
  xtick=\empty,
  ytick=\empty,
  point meta min=-2.8,
  point meta max=2.8,
  colormap name=turbo,
  title style={font=\tiny},
  label style={font=\tiny},
  scale only axis
]
\addplot graphics [
  xmin=0, xmax=0.992188,
  ymin=0, ymax=0.992188
] {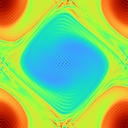};
\end{axis}

\begin{axis}[
  xlabel={$x$},
  at={(0.2023\figurewidth, 0.0779\figureheight)},
  anchor=south west,
  width=0.1834\figurewidth,
  height=0.8312\figureheight,
  title={Lag-PINN ($N_t=1$, mlp)},
  xmin=0, xmax=0.992188,
  ymin=0, ymax=0.992188,
  xtick=\empty,
  ytick=\empty,
  point meta min=-2.8,
  point meta max=2.8,
  colormap name=turbo,
  title style={font=\tiny},
  label style={font=\tiny},
  scale only axis
]
\addplot graphics [
  xmin=0, xmax=0.992188,
  ymin=0, ymax=0.992188
] {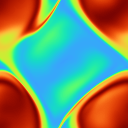};
\end{axis}

\begin{axis}[
  xlabel={$x$},
  at={(0.3874\figurewidth, 0.0779\figureheight)},
  anchor=south west,
  width=0.1834\figurewidth,
  height=0.8312\figureheight,
  title={Lag-PINN ($N_t=8$, mlp)},
  xmin=0, xmax=0.992188,
  ymin=0, ymax=0.992188,
  xtick=\empty,
  ytick=\empty,
  point meta min=-2.8,
  point meta max=2.8,
  colormap name=turbo,
  title style={font=\tiny},
  label style={font=\tiny},
  scale only axis
]
\addplot graphics [
  xmin=0, xmax=0.992188,
  ymin=0, ymax=0.992188
] {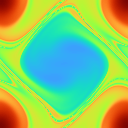};
\end{axis}

\begin{axis}[
  xlabel={$x$},
  at={(0.5725\figurewidth, 0.0779\figureheight)},
  anchor=south west,
  width=0.1834\figurewidth,
  height=0.8312\figureheight,
  title={PINN (classical)},
  xmin=0, xmax=0.992188,
  ymin=0, ymax=0.992188,
  xtick=\empty,
  ytick=\empty,
  point meta min=-2.8,
  point meta max=2.8,
  colormap name=turbo,
  colorbar,
  colorbar style={width=0.0101\figurewidth},
  title style={font=\tiny},
  label style={font=\tiny},
  scale only axis
]
\addplot graphics [
  xmin=0, xmax=0.992188,
  ymin=0, ymax=0.992188
] {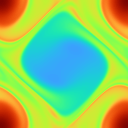};
\end{axis}
\end{tikzpicture}
}%

    \caption{Four-modes test case: vorticity
        $\omega(t=8,x,y)$ for the pseudo-spectral reference, the Lagrangian PINN
        with $N_t=1$ and $N_t=8$ flow-map segments, and the classical PINN.}
    \label{fig:euler-4modes-vorticity}
\end{figure}

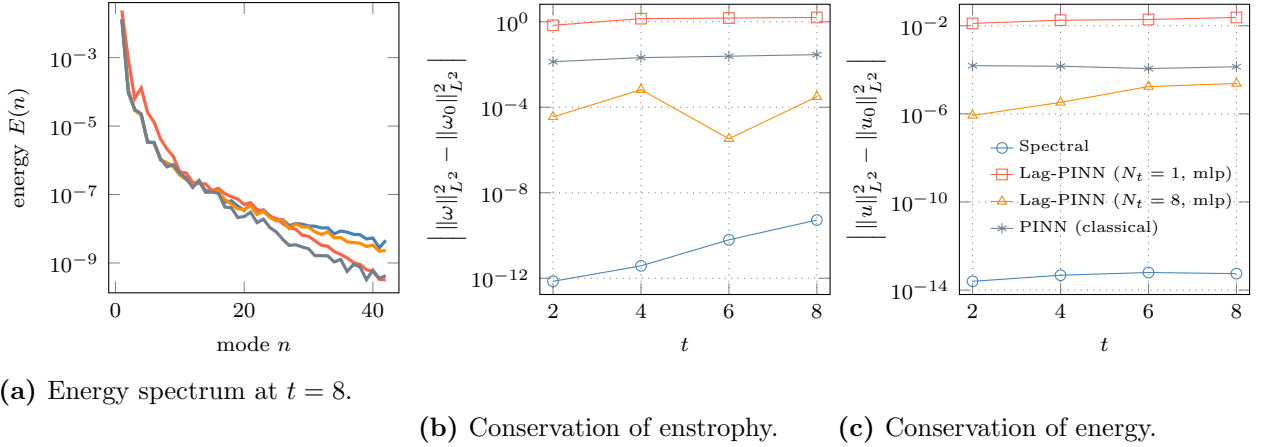
\begin{figure}[htp!]
    \centering
    \begin{subfigure}[t]{0.32\linewidth}
        \centering
        \setlength{\figurewidth}{1\linewidth}
        \setlength{\figureheight}{1\linewidth}
        %
    \definecolor{mplcolor0}{rgb}{0.274510,0.509800,0.705880}
\definecolor{mplcolor1}{rgb}{1.000000,0.388240,0.278430}
\definecolor{mplcolor2}{rgb}{1.000000,0.549020,0.000000}
\definecolor{mplcolor3}{rgb}{0.439220,0.501960,0.564710}
\definecolor{mplcolor4}{rgb}{0.501960,0.501960,0.501960}
\begin{tikzpicture}
\begin{axis}[
  width=\figurewidth,
  height=\figureheight,
  ymode=log, log basis y=10,
  xlabel={mode $n$},
  ylabel={energy $E(n)$},
  xmin=-1.05, xmax=44.05,
  legend style={
    font=\tiny,
    legend cell align={left},
    at={(0.99,0.99)}, anchor=north east, draw=none,
    fill=none, fill opacity=0.85, draw opacity=1, text opacity=1
  },
  ymin=1.34104e-10, ymax=0.0405845
]
\addplot [color=mplcolor0, mark=none, line width=1.20pt]
  coordinates {(1,0.0130901) (2,0.000110824) (3,2.9669e-05) (4,2.13033e-05) (5,3.34839e-06) (6,3.3207e-06) (7,1.02133e-06) (8,8.27465e-07) (9,6.90184e-07) (10,3.95647e-07) (11,3.03224e-07) (12,1.79115e-07) (13,2.55014e-07) (14,1.28603e-07) (15,1.31082e-07) (16,1.27574e-07) (17,5.94546e-08) (18,6.54164e-08) (19,4.63988e-08) (20,3.43537e-08) (21,5.43647e-08) (22,3.03211e-08) (23,3.55089e-08) (24,2.38328e-08) (25,1.79583e-08) (26,1.85159e-08) (27,1.23241e-08) (28,1.42021e-08) (29,1.36258e-08) (30,1.20947e-08) (31,1.18597e-08) (32,1.05568e-08) (33,1.05999e-08) (34,7.11054e-09) (35,8.54924e-09) (36,7.91175e-09) (37,6.85637e-09) (38,6.69155e-09) (39,5.02927e-09) (40,5.36214e-09) (41,2.78057e-09) (42,4.58832e-09)};
\addplot [color=mplcolor1, mark=none, line width=1.20pt]
  coordinates {(1,0.023952) (2,0.000898691) (3,6.37897e-05) (4,0.000129927) (5,2.46587e-05) (6,1.31454e-05) (7,4.24785e-06) (8,2.16713e-06) (9,1.28564e-06) (10,4.53751e-07) (11,4.24606e-07) (12,2.26117e-07) (13,2.08639e-07) (14,2.04115e-07) (15,1.30072e-07) (16,1.50946e-07) (17,1.06602e-07) (18,8.69269e-08) (19,7.43602e-08) (20,5.10653e-08) (21,5.7585e-08) (22,3.61336e-08) (23,3.36821e-08) (24,2.19958e-08) (25,2.01121e-08) (26,1.42125e-08) (27,1.08086e-08) (28,9.74528e-09) (29,6.52103e-09) (30,6.01921e-09) (31,4.46439e-09) (32,3.29761e-09) (33,3.0605e-09) (34,2.1606e-09) (35,1.7762e-09) (36,1.41296e-09) (37,1.02269e-09) (38,9.09279e-10) (39,6.37873e-10) (40,5.63437e-10) (41,3.32874e-10) (42,3.18921e-10)};
\addplot [color=mplcolor2, mark=none, line width=1.20pt]
  coordinates {(1,0.0131041) (2,0.000110255) (3,2.86269e-05) (4,2.11251e-05) (5,3.37236e-06) (6,3.23143e-06) (7,1.05478e-06) (8,8.35519e-07) (9,6.50365e-07) (10,3.69281e-07) (11,2.68473e-07) (12,1.84702e-07) (13,2.39094e-07) (14,1.20334e-07) (15,1.1874e-07) (16,1.26182e-07) (17,7.14962e-08) (18,5.29618e-08) (19,3.8273e-08) (20,3.65094e-08) (21,4.70387e-08) (22,2.60211e-08) (23,3.26753e-08) (24,2.12013e-08) (25,1.71064e-08) (26,1.81338e-08) (27,1.01183e-08) (28,1.38726e-08) (29,1.1354e-08) (30,1.12993e-08) (31,8.18145e-09) (32,6.7984e-09) (33,8.14678e-09) (34,5.39836e-09) (35,5.54505e-09) (36,4.17069e-09) (37,3.94548e-09) (38,3.78746e-09) (39,3.03726e-09) (40,3.30677e-09) (41,2.20869e-09) (42,2.31202e-09)};
\addplot [color=mplcolor3, mark=none, line width=1.20pt]
  coordinates {(1,0.0130665) (2,9.25704e-05) (3,3.04755e-05) (4,2.27983e-05) (5,3.37661e-06) (6,3.30356e-06) (7,1.11563e-06) (8,6.4245e-07) (9,7.38372e-07) (10,4.14172e-07) (11,2.99043e-07) (12,1.60759e-07) (13,2.52384e-07) (14,1.1599e-07) (15,1.1874e-07) (16,1.08152e-07) (17,4.27851e-08) (18,4.13337e-08) (19,2.12982e-08) (20,2.26785e-08) (21,3.038e-08) (22,1.50551e-08) (23,1.89709e-08) (24,1.08115e-08) (25,8.02721e-09) (26,4.87976e-09) (27,3.2638e-09) (28,3.40145e-09) (29,2.93391e-09) (30,2.61057e-09) (31,1.6344e-09) (32,1.63672e-09) (33,1.61224e-09) (34,1.07356e-09) (35,1.27474e-09) (36,5.67769e-10) (37,9.07644e-10) (38,7.67336e-10) (39,3.2579e-10) (40,6.14314e-10) (41,3.54045e-10) (42,4.26641e-10)};
\addplot [color=mplcolor4, dotted, mark=none, line width=0.80pt, forget plot]
  coordinates {(42,0) (42,1)};
\end{axis}
\end{tikzpicture}%

        \caption{Energy spectrum at $t=8$.}
        \label{fig:euler-4modes-energy-spectra}
    \end{subfigure}
    \begin{subfigure}[t]{0.32\linewidth}
        \centering
        \setlength{\figurewidth}{1\linewidth}
        \setlength{\figureheight}{1\linewidth}
        %
    \definecolor{mplcolor0}{rgb}{0.274510,0.509800,0.705880}
\definecolor{mplcolor1}{rgb}{1.000000,0.388240,0.278430}
\definecolor{mplcolor2}{rgb}{1.000000,0.549020,0.000000}
\definecolor{mplcolor3}{rgb}{0.439220,0.501960,0.564710}
\begin{tikzpicture}
\begin{axis}[
  at={(0.0\figurewidth, 0.1367\figureheight)},
  anchor=south west,
  width=\figurewidth,
  height=\figureheight,
  ymode=log, log basis y=10,
  ylabel style={yshift=-0.2cm},
  xlabel={$t$},
  ylabel={$\left|\,\|\omega\|_{L^2}^2 - \|\omega_0\|_{L^2}^2\,\right|$},
  xmin=1.7, xmax=8.3,
  ymin=1.73826e-13, ymax=6.56353,
  xmajorgrids=true,
  ymajorgrids=true
]
\addplot [color=mplcolor0, mark=o]
  coordinates {(2,7.19869e-13) (4,3.80429e-12) (6,6.28564e-11) (8,5.29081e-10)};
\addplot [color=mplcolor1, mark=square]
  coordinates {(2,0.681367) (4,1.3959) (6,1.49186) (8,1.58489)};
\addplot [color=mplcolor2, mark=triangle]
  coordinates {(2,3.52959e-05) (4,0.000665025) (6,3.35129e-06) (8,0.000310673)};
\addplot [color=mplcolor3, mark=asterisk]
  coordinates {(2,0.0134698) (4,0.0209549) (6,0.0244203) (8,0.0289285)};
\end{axis}
\end{tikzpicture}%

        \caption{Conservation of enstrophy.}
        \label{fig:euler-4modes-conservation-enstrophy}
    \end{subfigure}
    \begin{subfigure}[t]{0.32\linewidth}
        \centering
        \setlength{\figurewidth}{1\linewidth}
        \setlength{\figureheight}{1\linewidth}
        %
    \definecolor{mplcolor0}{rgb}{0.274510,0.509800,0.705880}
\definecolor{mplcolor1}{rgb}{1.000000,0.388240,0.278430}
\definecolor{mplcolor2}{rgb}{1.000000,0.549020,0.000000}
\definecolor{mplcolor3}{rgb}{0.439220,0.501960,0.564710}
\begin{tikzpicture}
\begin{axis}[
  at={(0.44\figurewidth, 0.1367\figureheight)},
  anchor=south west,
  width=\figurewidth,
  height=\figureheight,
  ymode=log, log basis y=10,
  ylabel style={yshift=-0.2cm},
  xlabel={$t$},
  ylabel={$\left|\,\|u\|_{L^2}^2 - \|u_0\|_{L^2}^2\,\right|$},
  xmin=1.7, xmax=8.3,
  ymin=6.28083e-15, ymax=0.0944935,
  xmajorgrids=true,
  ymajorgrids=true,
  legend style={
    font=\tiny,
    legend cell align={left},
    at={(0.99,0.57)}, anchor=north east, draw=none,
    fill=none, fill opacity=0.85, draw opacity=1, text opacity=1
  }
]
\addplot [color=mplcolor0, mark=o]
  coordinates {(2,2.49453e-14) (4,4.71359e-14) (6,6.27068e-14) (8,5.56846e-14)};
\addplot [color=mplcolor1, mark=square]
  coordinates {(2,0.0127662) (4,0.0179168) (6,0.0192665) (8,0.023792)};
\addplot [color=mplcolor2, mark=triangle]
  coordinates {(2,8.50131e-07) (4,3.3137e-06) (6,1.7355e-05) (8,2.39521e-05)};
\addplot [color=mplcolor3, mark=asterisk]
  coordinates {(2,0.000153394) (4,0.000145376) (6,0.000113994) (8,0.000137707)};
\addlegendentry{Spectral}
\addlegendentry{Lag-PINN ($N_t=1$, mlp)}
\addlegendentry{Lag-PINN ($N_t=8$, mlp)}
\addlegendentry{PINN (classical)}
\end{axis}
\end{tikzpicture}%

        \caption{Conservation of energy.}
        \label{fig:euler-4modes-conservation-energy}
    \end{subfigure}
    \caption{Four-modes test case: pseudo-spectral reference vs Lagrangian PINN with
        $N_t=1$ and $N_t=8$ flow-map segments, and classical PINN.}
    \label{fig:euler-4modes-spectra-analyticity}
\end{figure}

We compare: (i) a Fourier pseudo-spectral reference
solve~\cite{CanutoHussainiQuarteroniZang2007}, (ii) a
Lagrangian PINN using a single flow map over $[0,T]$ ($N_t=1$) similar to \cite{MontiMorozovskaPoggiShukla2026}, (iii) a
Lagrangian PINN composed of $N_t=8$ backward-flow submaps, and (iv) a classical Eulerian PINN.
The reference solver evaluates the nonlinear term $\vec u\cdot\nabla\omega$ in physical space, removes aliasing errors with the $2/3$ rule~\cite{Orszag1971}, and advances the Fourier coefficients with the classical fourth-order Runge--Kutta scheme, using $N=128$ points per direction and $\Delta t = 5\cdot10^{-3}$.
\Cref{fig:euler-4modes-vorticity} shows the vorticity field of all methods
at $t=8$, and \cref{fig:euler-4modes-energy-spectra} the corresponding Fourier energy spectra. Furthermore, the conservation of enstrophy and energy is plotted over time in \cref{fig:euler-4modes-conservation-energy,fig:euler-4modes-conservation-enstrophy}.
The errors at $t=8$ are summarized in \cref{tab:euler-4modes}. With a single map ($N_t=1$), the enstrophy error reaches $1.6$, whereas the composition of $N_t=8$ submaps reduces it to $3.1\cdot10^{-4}$, about two orders of magnitude below the Eulerian PINN ($2.9\cdot10^{-2}$). The energy errors follow the same ranking ($2.4\cdot10^{-2}$, $2.4\cdot10^{-5}$ and $1.4\cdot10^{-4}$, respectively). The characteristic mapping method of~\cite{YinMercierYadavSchneiderNave2021} conserves both quantities more accurately, with $105$ submaps instead of $8$.

\begin{table}[htp!]
    \centering
    \caption{Four-modes test case at $t=8$: runtime,
        enstrophy/energy conservation error \eqref{eq:enstrophy_energy_def}, and
        number of submaps. The values for the characteristic mapping method (CMM) are
        literature results from \cite{YinMercierYadavSchneiderNave2021},
        Tables~4.1--4.2, obtained on an older Intel Core i5-2320 CPU. All other
        methods run on an NVIDIA V100 GPU, so the runtimes are not directly
        comparable.}
    \label{tab:euler-4modes}
    \begin{tabular}{lrrrr}
        \toprule
        Method                                      & Runtime [s] & Enstrophy error        & Energy error          & \# maps \\
        \midrule
        Spectral                                    & 4.9         & $-5.291\times10^{-10}$ & $5.568\times10^{-14}$ & 0       \\
        CMM \cite{YinMercierYadavSchneiderNave2021} & 359.0       & $1.370\times10^{-5}$   & $-7.510\times10^{-6}$ & 105     \\
        Lagrangian PINN ($N_t=1$)                   & 448.7       & $1.585$                & $2.379\times10^{-2}$  & 1       \\
        Lagrangian PINN ($N_t=8$)                   & 757.7       & $-3.107\times10^{-4}$  & $2.395\times10^{-5}$  & 8       \\
        Classical PINN                              & 548.8       & $-2.893\times10^{-2}$  & $-1.377\times10^{-4}$ & 0       \\
        \bottomrule
    \end{tabular}
\end{table}

\subsubsection{Vortex merger and subgrid resolution property}
\label{sec:euler_vortex_merger}

We consider the ``vortex merger'' test case from \cite{YinMercierYadavSchneiderNave2021}, with initial vorticity given by two co-rotating Gaussian blobs,
\begin{equation}
    \omega_0(x,y) = \frac{1}{\sigma\sqrt{2\pi}}\left[
        \exp\!\left(-\frac{(x-x_1)^2+(y-0.5)^2}{2\sigma^2}\right)
        + \exp\!\left(-\frac{(x-x_2)^2+(y-0.5)^2}{2\sigma^2}\right)
        \right],
    \label{eq:euler_vortex_merger_ic}
\end{equation}
with blob centers $x_{1,2} = 0.5 \mp 0.15$ (separation $0.3$) and width $\sigma = 0.07$, inside the domain $[0,1]^2$ with periodic boundary conditions. Unlike the four-modes case, the two same-sign, co-rotating blobs merge into a single, increasingly tightly wound spiral vortex.

\begin{figure}[htp!]
    \centering
    \setlength{\figurewidth}{1\linewidth}
    \setlength{\figureheight}{0.2206\linewidth}
    %
    {\pgfplotsset{compat=1.18}
\pgfplotsset{colormap={turbo}{rgb255(0.0000)=(48,18,59) rgb255(0.0159)=(53,30,88) rgb255(0.0317)=(57,42,115) rgb255(0.0476)=(61,53,139) rgb255(0.0635)=(64,64,162) rgb255(0.0794)=(66,75,181) rgb255(0.0952)=(68,86,199) rgb255(0.1111)=(70,97,214) rgb255(0.1270)=(70,107,227) rgb255(0.1429)=(71,118,238) rgb255(0.1587)=(70,128,246) rgb255(0.1746)=(69,138,252) rgb255(0.1905)=(66,148,255) rgb255(0.2063)=(61,158,254) rgb255(0.2222)=(55,168,250) rgb255(0.2381)=(47,178,244) rgb255(0.2540)=(39,190,233) rgb255(0.2698)=(32,199,223) rgb255(0.2857)=(27,208,213) rgb255(0.3016)=(24,215,202) rgb255(0.3175)=(24,222,192) rgb255(0.3333)=(26,228,182) rgb255(0.3492)=(32,234,172) rgb255(0.3651)=(42,239,161) rgb255(0.3810)=(53,243,148) rgb255(0.3968)=(67,247,135) rgb255(0.4127)=(82,250,122) rgb255(0.4286)=(97,252,108) rgb255(0.4444)=(113,254,95) rgb255(0.4603)=(128,255,83) rgb255(0.4762)=(143,255,73) rgb255(0.4921)=(156,254,64) rgb255(0.5079)=(169,251,57) rgb255(0.5238)=(180,248,54) rgb255(0.5397)=(190,244,52) rgb255(0.5556)=(200,239,52) rgb255(0.5714)=(210,233,53) rgb255(0.5873)=(219,226,54) rgb255(0.6032)=(227,219,56) rgb255(0.6190)=(235,211,57) rgb255(0.6349)=(241,203,58) rgb255(0.6508)=(246,195,58) rgb255(0.6667)=(250,186,57) rgb255(0.6825)=(252,177,54) rgb255(0.6984)=(254,167,50) rgb255(0.7143)=(254,155,45) rgb255(0.7302)=(254,144,41) rgb255(0.7460)=(252,132,35) rgb255(0.7619)=(249,117,29) rgb255(0.7778)=(245,105,24) rgb255(0.7937)=(241,93,19) rgb255(0.8095)=(236,83,15) rgb255(0.8254)=(231,73,12) rgb255(0.8413)=(225,65,9) rgb255(0.8571)=(218,57,7) rgb255(0.8730)=(210,49,5) rgb255(0.8889)=(202,42,4) rgb255(0.9048)=(193,35,2) rgb255(0.9206)=(183,29,2) rgb255(0.9365)=(172,23,1) rgb255(0.9524)=(161,18,1) rgb255(0.9683)=(149,13,1) rgb255(0.9841)=(136,8,2) rgb255(1.0000)=(122,4,3)}}
\begin{tikzpicture}
\begin{axis}[
  xlabel={$x$},
  ylabel={$y$},
  at={(0.0172\figurewidth, 0.0779\figureheight)},
  anchor=south west,
  width=0.1834\figurewidth,
  height=0.8312\figureheight,
  title={Spectral ($N=512$)},
  xmin=0, xmax=0.998047,
  ymin=0, ymax=0.998047,
  xtick=\empty,
  ytick=\empty,
  point meta min=0,
  point meta max=5.69966,
  colormap name=turbo,
  title style={font=\tiny},
  label style={font=\tiny},
  scale only axis
]
\addplot graphics [
  xmin=0, xmax=0.998047,
  ymin=0, ymax=0.998047
] {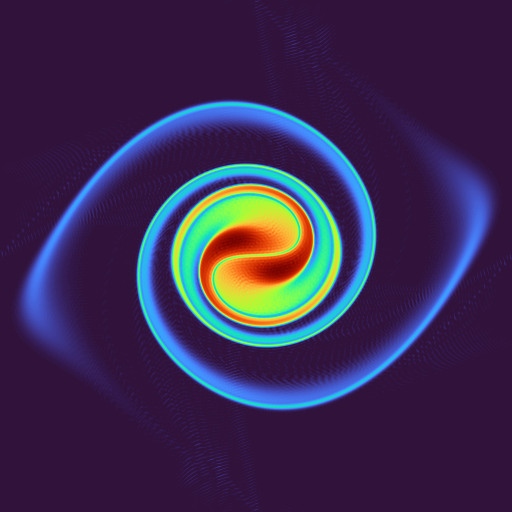};
\end{axis}

\begin{axis}[
  xlabel={$x$},
  at={(0.2023\figurewidth, 0.0779\figureheight)},
  anchor=south west,
  width=0.1834\figurewidth,
  height=0.8312\figureheight,
  title={Spectral ($N=128$)},
  xmin=0, xmax=0.998047,
  ymin=0, ymax=0.998047,
  xtick=\empty,
  ytick=\empty,
  point meta min=0,
  point meta max=5.69966,
  colormap name=turbo,
  title style={font=\tiny},
  label style={font=\tiny},
  scale only axis
]
\addplot graphics [
  xmin=0, xmax=0.998047,
  ymin=0, ymax=0.998047
] {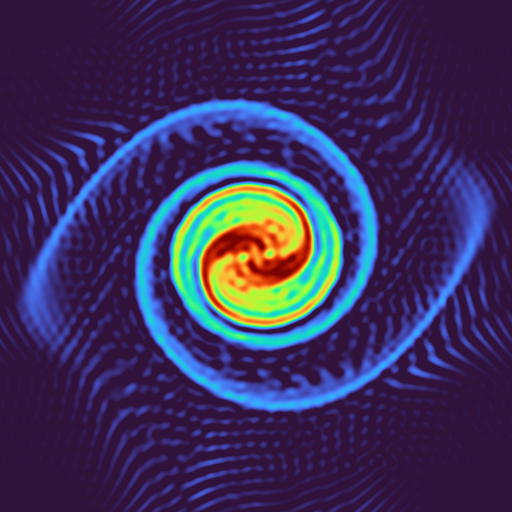};
\end{axis}

\begin{axis}[
  xlabel={$x$},
  at={(0.3874\figurewidth, 0.0779\figureheight)},
  anchor=south west,
  width=0.1834\figurewidth,
  height=0.8312\figureheight,
  title={Lag-PINN ($N_t=1$, mlp)},
  xmin=0, xmax=0.998047,
  ymin=0, ymax=0.998047,
  xtick=\empty,
  ytick=\empty,
  point meta min=0,
  point meta max=5.69966,
  colormap name=turbo,
  title style={font=\tiny},
  label style={font=\tiny},
  scale only axis
]
\addplot graphics [
  xmin=0, xmax=0.998047,
  ymin=0, ymax=0.998047
] {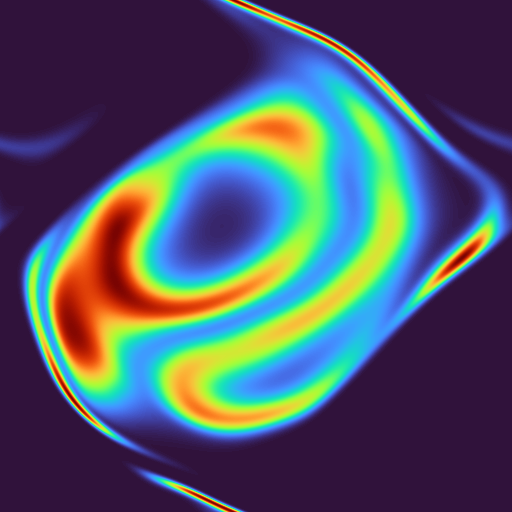};
\end{axis}

\begin{axis}[
  xlabel={$x$},
  at={(0.5725\figurewidth, 0.0779\figureheight)},
  anchor=south west,
  width=0.1834\figurewidth,
  height=0.8312\figureheight,
  title={Lag-PINN ($N_t=8$, mlp)},
  xmin=0, xmax=0.998047,
  ymin=0, ymax=0.998047,
  xtick=\empty,
  ytick=\empty,
  point meta min=0,
  point meta max=5.69966,
  colormap name=turbo,
  title style={font=\tiny},
  label style={font=\tiny},
  scale only axis
]
\addplot graphics [
  xmin=0, xmax=0.998047,
  ymin=0, ymax=0.998047
] {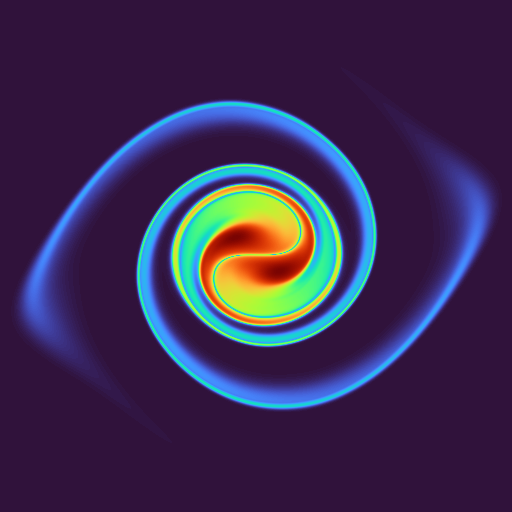};
\end{axis}

\begin{axis}[
  xlabel={$x$},
  at={(0.7576\figurewidth, 0.0779\figureheight)},
  anchor=south west,
  width=0.1834\figurewidth,
  height=0.8312\figureheight,
  title={PINN (classical)},
  xmin=0, xmax=0.998047,
  ymin=0, ymax=0.998047,
  xtick=\empty,
  ytick=\empty,
  point meta min=0,
  point meta max=5.69966,
  colormap name=turbo,
  colorbar,
  colorbar style={width=0.0101\figurewidth},
  title style={font=\tiny},
  label style={font=\tiny},
  scale only axis
]
\addplot graphics [
  xmin=0, xmax=0.998047,
  ymin=0, ymax=0.998047
] {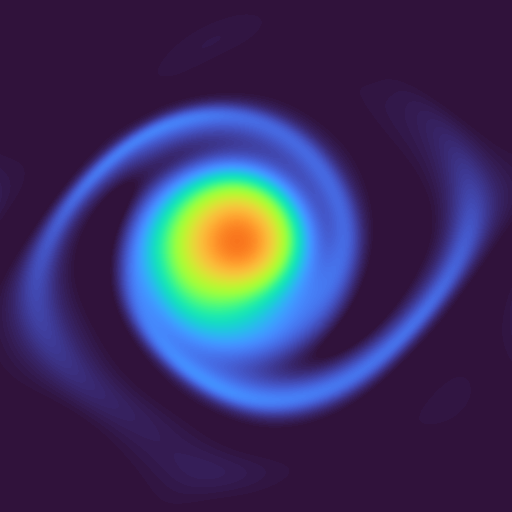};
\end{axis}
\end{tikzpicture}
}%

    \caption{Vortex-merger test case: vorticity
        $\omega(t=8,x,y)$ for the pseudo-spectral reference at $N=512$ and at $N=128$, the Lagrangian PINN
        with $N_t=1$ and $N_t=8$ flow-map segments, and the classical PINN.}
    \label{fig:euler-vortex-vorticity}
\end{figure}

\begin{figure}[htp!]
    \centering
    \setlength{\figurewidth}{0.5\linewidth}
    \setlength{\figureheight}{0.4\linewidth}
    %
    \input{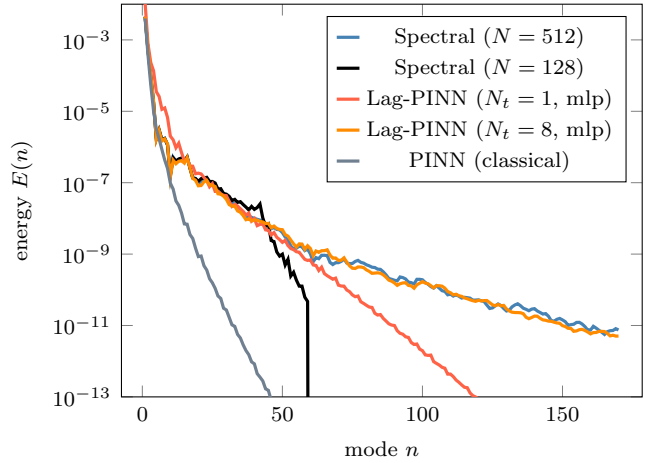}%

    \caption{Vortex-merger test case: Comparison of the Fourier spectra $E(n)$ for the pseudo-spectral references ($N=512$, $N=128$) vs Lagrangian PINN with
        $N_t=1$ and $N_t=8$ flow-map segments, and classical PINN.}
    \label{fig:euler-vortex-spectra-analyticity}
\end{figure}

We use the same final time and the same methods as in \cref{sec:euler_four_modes}: a pseudo-spectral reference, a Lagrangian PINN using a single flow map over $[0,T]$ ($N_t=1$), a Lagrangian PINN composed of $N_t=8$ backward-flow submaps, and a classical Eulerian PINN. The final snapshot at $t=T=8$ is shown in \cref{fig:euler-vortex-vorticity}, its Fourier energy spectrum in \cref{fig:euler-vortex-spectra-analyticity}.

Unlike the smoother four-modes flow, the vortex-merger vorticity keeps sharpening into thinner filaments over time, so a fixed-resolution spectral solve eventually falls behind the true solution: at $t=8$, an $N=128$ pseudo-spectral run (matching the resolution of the periodic FFT Biot-Savart grid the PINN methods are trained against) already differs from a properly resolved $N=512$ run by about $20\%$ in relative $L^2$ norm, even though its enstrophy and energy remain conserved to near machine precision. We report both resolutions, labeled ``Spectral ($N=512$)'' and ``Spectral ($N=128$)'', in \cref{fig:euler-vortex-vorticity,tab:euler-vortex}, and use the $N=512$ run as the ground truth. Every PINN method is still trained against the coarser $N=128$ Biot-Savart grid (as in \cref{sec:euler_four_modes}); only its post-training field is re-evaluated on the finer $N=512$ grid for a fair comparison.

As already observed in \cite{YinMercierYadavSchneiderNave2021} the use of submaps helps to resolve the fine scales of the solution without having to resolve them entirely during the simulation, i.e., during PINN training. This is illustrated in \cref{fig:euler-vortex-spectra-analyticity}, where a PINN trained with a $128^2$ Poisson grid resolves the Fourier spectra of the $512^2$ pseudo-spectral simulation, up to a regime where the $128^2$ pseudo-spectral simulation already breaks down.

\begin{table}[htp!]
    \centering
    \caption{Vortex-merger test case at $t=8$: runtime,
        enstrophy/energy conservation error \eqref{eq:enstrophy_energy_def}, and
        number of submaps. All methods run on an NVIDIA V100 GPU.}
    \label{tab:euler-vortex}
    \begin{tabular}{lrrrr}
        \toprule
        Method                    & Runtime [s] & Enstrophy error        & Energy error           & \# maps \\
        \midrule
        Spectral ($N=512$)        & 14.9        & $-7.271\times10^{-8}$  & $-9.080\times10^{-14}$ & 0       \\
        Spectral ($N=128$)        & 9.2         & $-3.886\times10^{-13}$ & $-1.561\times10^{-17}$ & 0       \\
        Lagrangian PINN ($N_t=1$) & 441.4       & $2.097$                & $1.338\times10^{-2}$   & 1       \\
        Lagrangian PINN ($N_t=8$) & 709.6       & $-6.578\times10^{-3}$  & $7.482\times10^{-6}$   & 8       \\
        Classical PINN            & 538.1       & $-2.771\times10^{-1}$  & $-1.186\times10^{-3}$  & 0       \\
        \bottomrule
    \end{tabular}
\end{table}

\subsection{Vlasov-Poisson}
\label{sec:vlasov_poisson}
We now consider the nonlinear Vlasov-Poisson system in 1D1V~\cite{vlasov1938vibrational,krall1986principles}
\begin{equation*}
    \partial_t f(t,x,v) + v\,\partial_x f(t,x,v) + E(t,x)\,\partial_v f(t,x,v) = 0,
\end{equation*}
in the domain $[0,T] \times L \mathbb{T} \times \mathbb{R}$, with $L\mathbb{T} = \mathbb{R} / L \mathbb{Z}$.
The electric field $E = -\partial_x \phi$ is generated
self-consistently by the distribution function through the Poisson equation
\begin{equation*}
    -\partial_{xx} \phi(t,x) =\int_{\mathbb{R}} f(t,x,v)\,\dv- 1.
\end{equation*}
The coupling with the Poisson equation is handled by the Picard iteration described in~\cref{sec:field_coupling}. At each optimization step, the density $\rho$ is
evaluated by a Gauss quadrature in $v$, and the Poisson equation is solved by a spectral scheme~\cite{BirdsallLangdon2018} on a uniform grid in $x$. The reference solutions are given by a backward semi-Lagrangian
scheme~\cite{SonnendrueckerRocheBertrandGhizzo1999}.

\newcount\vpdepth
\newcommand{\phasespacebbox}{\tikzset{every picture/.append style={
                execute at begin picture={\advance\vpdepth by 1\relax},
                execute at end picture={\ifnum\vpdepth=1
                            \path (0,0) rectangle (\figurewidth,\figureheight);\fi}}}%
    \pgfplotsset{every colorbar/.append style={yticklabel style={
                        /pgf/number format/fixed, /pgf/number format/precision=2}}}}

\paragraph{Landau damping.} We first consider a classical Landau damping~\cite{ChenKnorr1976,SonnendrueckerRocheBertrandGhizzo1999,FilbetSonnendrucker2003}.
The initial condition is a perturbed Maxwellian,
\begin{equation*}
    f_0(x,v) = \big(1 + \varepsilon \cos(k x)\big)\, \frac{e^{-v^2/2}}{\sqrt{2\pi}},
    \qquad \varepsilon = 0.1, \quad k = 0.5,
\end{equation*}
the computational domain is $[0,L] \times [-6,6]$ with  $L = 2\pi/k$. A linear analysis predicts that the electric
field decays exponentially, with a damping rate
$\gamma \approx -0.1533$~\cite{Landau1946,Canosa1973}. The goal of this test is to recover this rate of decay with Lagrangian PINNs.

The final time is $T = 24$ and the time interval $[0,T]$ is split into $N_t = 6$ subintervals of length $4$.
Each MLP submap has $2$ hidden layers of $30$ neurons (so the number of parameters per submap is $N_\theta = 1142$), while each
SympNet submap has $13$ layers of width $8$ ($N_\theta = 1144$). We use  $400$ optimization steps for each submap. We use $5000$ collocation points and
$2000$ initial-condition points at each optimization step. In the Picard iteration, $\rho$ is computed
with $128$ quadrature points in $v$ and the Poisson equation is solved on
$32$ points in $x$, at $N_\tau = 80$ times per subinterval. The Eulerian PINN
has $3$ hidden layers of $64$ neurons ($N_\theta = 8705$) and is trained for
$2000$ steps on $10000$ collocation points and $4000$ initial-condition points,
the field being computed at $200$ times over $[0, T]$. The reference uses a
$512 \times 512$ grid and a time step $\Delta t = 0.05$.

\begin{figure}[!ht]
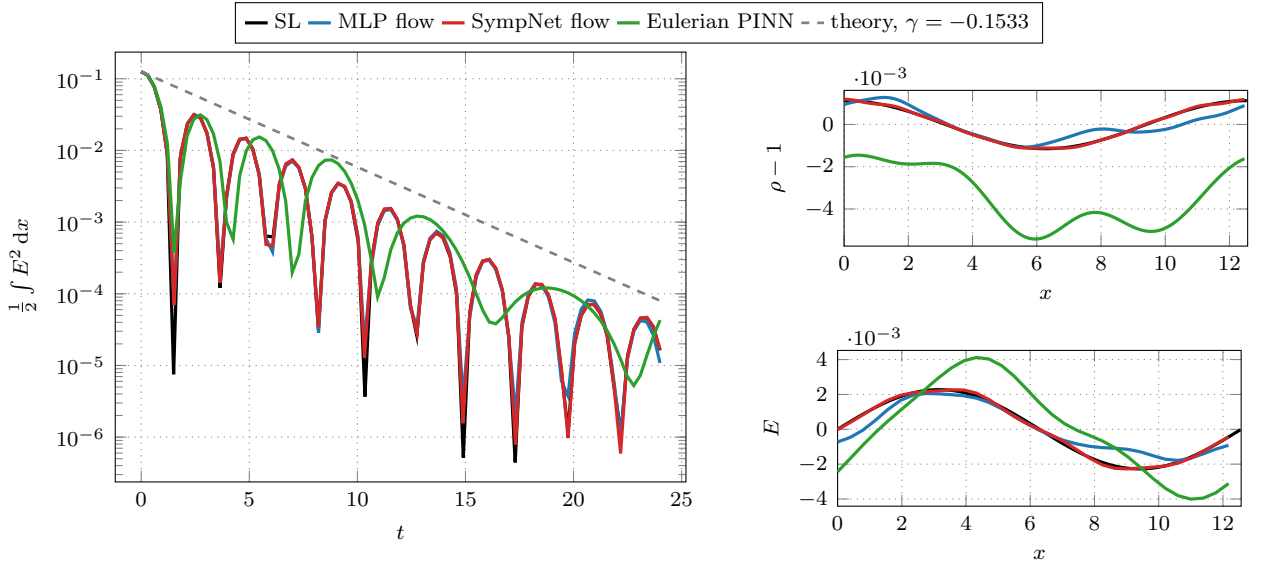

    \centering
    \input{figures/vlasov_poisson/vp_colors}\ref*{vplegendlandaudampingenergy}\\[2pt]
    \begin{minipage}[c]{0.55\linewidth}\centering
        \setlength{\figurewidth}{0.98\linewidth}
        \setlength{\figureheight}{0.78\linewidth}
        %
    \definecolor{mplcolor0}{rgb}{0.000000,0.000000,0.000000}
\definecolor{mplcolor1}{rgb}{0.121570,0.466670,0.705880}
\definecolor{mplcolor2}{rgb}{0.839220,0.152940,0.156860}
\definecolor{mplcolor3}{rgb}{0.172550,0.627450,0.172550}
\definecolor{mplcolor4}{rgb}{0.501960,0.501960,0.501960}
\begin{tikzpicture}
\begin{axis}[
  width=\figurewidth,
  height=\figureheight,
  ymode=log, log basis y=10,
  xlabel={$t$},
  ylabel={$\frac{1}{2}\int E^2\,\mathrm{d}x$},
  xmin=-1.2, xmax=25.2,
  ymin=2.36986e-07, ymax=0.235399,
  xmajorgrids=true,
  ymajorgrids=true,
  legend to name=vplegendlandaudampingenergy,
  legend columns=5
]
\addplot [color=mplcolor0, mark=none, line width=1.20pt]
  coordinates {(0,0.125664) (0.3,0.112096) (0.6,0.0777318) (0.9,0.0380645) (1.2,0.00929485) (1.5,7.53788e-06) (1.8,0.00771553) (2.15,0.024023) (2.45,0.0316184) (2.75,0.0282417) (3.05,0.0170657) (3.35,0.00555425) (3.65,0.000120769) (3.95,0.00235321) (4.25,0.00886052) (4.55,0.0142483) (4.85,0.014832) (5.15,0.0106282) (5.45,0.00468519) (5.75,0.000639127) (6.1,0.000621189) (6.4,0.00352754) (6.7,0.00646878) (7,0.00734322) (7.3,0.00574379) (7.6,0.00290169) (7.9,0.000617223) (8.2,3.16706e-05) (8.5,0.00103638) (8.8,0.00257069) (9.1,0.00345704) (9.4,0.00315138) (9.7,0.00195836) (10.05,0.000515617) (10.35,3.68184e-06) (10.65,0.000292381) (10.95,0.00098125) (11.25,0.00151092) (11.55,0.0015344) (11.85,0.00108493) (12.15,0.000479502) (12.45,7.03631e-05) (12.75,2.59976e-05) (13.05,0.000267282) (13.35,0.000564767) (13.65,0.000705206) (13.95,0.000613733) (14.3,0.000323896) (14.6,9.38555e-05) (14.9,5.18898e-07) (15.2,5.50938e-05) (15.5,0.000183469) (15.8,0.000286196) (16.1,0.000299943) (16.4,0.000225388) (16.7,0.000113333) (17,2.60169e-05) (17.3,4.43933e-07) (17.6,3.32805e-05) (17.9,9.07499e-05) (18.25,0.000136255) (18.55,0.000132042) (18.85,9.30943e-05) (19.15,4.30862e-05) (19.45,7.76964e-06) (19.75,1.1843e-06) (20.05,2.033e-05) (20.35,4.94379e-05) (20.65,7.00435e-05) (20.95,7.08807e-05) (21.25,5.26886e-05) (21.55,2.62641e-05) (21.85,5.6289e-06) (22.2,1.12048e-06) (22.5,1.37133e-05) (22.8,3.23553e-05) (23.1,4.56578e-05) (23.4,4.61238e-05) (23.7,3.39176e-05) (24,1.62375e-05)};
\addplot [color=mplcolor1, mark=none, line width=1.20pt]
  coordinates {(0,0.125653) (0.303797,0.111841) (0.607595,0.076851) (0.911392,0.0368156) (1.21519,0.00837843) (1.51899,7.03547e-05) (1.82278,0.00868693) (2.12658,0.0230212) (2.43038,0.0315015) (2.73418,0.0288418) (3.03797,0.0178315) (3.34177,0.00603506) (3.64557,0.000202098) (3.94937,0.00213047) (4.25316,0.00859198) (4.55696,0.01405) (4.86076,0.0146861) (5.16456,0.0105056) (5.46835,0.00457719) (5.77215,0.000586221) (6.07595,0.000384675) (6.37975,0.00306426) (6.68354,0.00604012) (6.98734,0.00709428) (7.29114,0.00570932) (7.59494,0.00296376) (7.89873,0.000638846) (8.20253,2.84623e-05) (8.50633,0.00102068) (8.81013,0.00255325) (9.11392,0.0034332) (9.41772,0.00312318) (9.72152,0.00192179) (10.0253,0.000643541) (10.3291,2.01483e-05) (10.6329,0.00022653) (10.9367,0.000884034) (11.2405,0.00143464) (11.5443,0.00150239) (11.8481,0.00106686) (12.1519,0.000464203) (12.4557,6.20984e-05) (12.7595,3.54975e-05) (13.0633,0.000293631) (13.3671,0.000600966) (13.6709,0.000744792) (13.9747,0.000643165) (14.2785,0.000376391) (14.5823,0.000118189) (14.8861,3.40324e-06) (15.1899,4.50861e-05) (15.4937,0.000168707) (15.7975,0.000283446) (16.1013,0.00029347) (16.4051,0.000219547) (16.7089,0.000108204) (17.0127,2.2526e-05) (17.3165,2.01649e-06) (17.6203,4.11724e-05) (17.9241,0.000100692) (18.2278,0.000137857) (18.5316,0.000131904) (18.8354,9.05431e-05) (19.1392,3.96176e-05) (19.443,6.05644e-06) (19.7468,3.75115e-06) (20.0506,2.80905e-05) (20.3544,6.07126e-05) (20.6582,8.19389e-05) (20.962,7.94452e-05) (21.2658,5.65756e-05) (21.5696,2.72088e-05) (21.8734,5.92664e-06) (22.1772,1.0102e-06) (22.481,1.1864e-05) (22.7848,2.98081e-05) (23.0886,4.22902e-05) (23.3924,4.04871e-05) (23.6962,2.6274e-05) (24,1.08366e-05)};
\addplot [color=mplcolor2, mark=none, line width=1.20pt]
  coordinates {(0,0.125664) (0.303797,0.111789) (0.607595,0.0767622) (0.911392,0.0367296) (1.21519,0.00836131) (1.51899,7.08284e-05) (1.82278,0.00868707) (2.12658,0.0229781) (2.43038,0.0313937) (2.73418,0.0286492) (3.03797,0.0175796) (3.34177,0.00581661) (3.64557,0.000143639) (3.94937,0.00233545) (4.25316,0.00892556) (4.55696,0.014321) (4.86076,0.0147628) (5.16456,0.0103674) (5.46835,0.00436609) (5.77215,0.000487801) (6.07595,0.000466539) (6.37975,0.00328933) (6.68354,0.00634357) (6.98734,0.00736709) (7.29114,0.005832) (7.59494,0.00296074) (7.89873,0.00062725) (8.20253,3.42131e-05) (8.50633,0.00107001) (8.81013,0.0026293) (9.11392,0.00349502) (9.41772,0.0031292) (9.72152,0.00188411) (10.0253,0.000607668) (10.3291,1.27235e-05) (10.6329,0.000255105) (10.9367,0.000947923) (11.2405,0.00150806) (11.5443,0.00155288) (11.8481,0.0011037) (12.1519,0.000486697) (12.4557,6.98232e-05) (12.7595,2.80445e-05) (13.0633,0.000276688) (13.3671,0.000576421) (13.6709,0.00070858) (13.9747,0.000603338) (14.2785,0.00034936) (14.5823,0.000107684) (14.8861,1.5607e-06) (15.1899,5.05859e-05) (15.4937,0.000180442) (15.7975,0.000286967) (16.1013,0.000302954) (16.4051,0.000227824) (16.7089,0.000113014) (17.0127,2.46401e-05) (17.3165,7.90463e-07) (17.6203,3.62427e-05) (17.9241,9.50014e-05) (18.2278,0.000135891) (18.5316,0.00013558) (18.8354,9.77623e-05) (19.1392,4.62827e-05) (19.443,8.79824e-06) (19.7468,9.69664e-07) (20.0506,1.99839e-05) (20.3544,4.98465e-05) (20.6582,7.09244e-05) (20.962,7.15574e-05) (21.2658,5.26594e-05) (21.5696,2.55786e-05) (21.8734,4.98794e-06) (22.1772,5.94095e-07) (22.481,1.23031e-05) (22.7848,3.13021e-05) (23.0886,4.56374e-05) (23.3924,4.69659e-05) (23.6962,3.48484e-05) (24,1.66296e-05)};
\addplot [color=mplcolor3, mark=none, line width=1.20pt]
  coordinates {(0,0.125577) (0.303797,0.111344) (0.607595,0.0783063) (0.911392,0.0407133) (1.21519,0.0122944) (1.51899,0.000374677) (1.82278,0.00398673) (2.12658,0.0159522) (2.43038,0.0272475) (2.73418,0.031399) (3.03797,0.0269431) (3.34177,0.0170557) (3.64557,0.00697882) (3.94937,0.000988498) (4.25316,0.000579457) (4.55696,0.00446907) (4.86076,0.00987486) (5.16456,0.0140616) (5.46835,0.0154009) (5.77215,0.0136932) (6.07595,0.00989069) (6.37975,0.00550406) (6.68354,0.00197286) (6.98734,0.000199199) (7.29114,0.000353476) (7.59494,0.00196004) (7.89873,0.00418047) (8.20253,0.0061588) (8.50633,0.00729536) (8.81013,0.0073704) (9.11392,0.00651333) (9.41772,0.00507184) (9.72152,0.00345443) (10.0253,0.00200456) (10.3291,0.000933334) (10.6329,0.000309767) (10.9367,9.05469e-05) (11.2405,0.00016697) (11.5443,0.000410256) (11.8481,0.000704481) (12.1519,0.000964229) (12.4557,0.00113936) (12.7595,0.0012113) (13.0633,0.00118526) (13.3671,0.00108128) (13.6709,0.000926196) (13.9747,0.000747341) (14.2785,0.000568231) (14.5823,0.000406293) (14.8861,0.000272209) (15.1899,0.000170479) (15.4937,0.00010071) (15.7975,5.91836e-05) (16.1013,4.03669e-05) (16.4051,3.81541e-05) (16.7089,4.67467e-05) (17.0127,6.11786e-05) (17.3165,7.75486e-05) (17.6203,9.30427e-05) (17.9241,0.000105828) (18.2278,0.000114882) (18.5316,0.000119805) (18.8354,0.00012064) (19.1392,0.000117725) (19.443,0.000111574) (19.7468,0.000102794) (20.0506,9.20227e-05) (20.3544,7.98995e-05) (20.6582,6.7046e-05) (20.962,5.40618e-05) (21.2658,4.15259e-05) (21.5696,3.00029e-05) (21.8734,2.00502e-05) (22.1772,1.22261e-05) (22.481,7.09659e-06) (22.7848,5.24123e-06) (23.0886,7.25817e-06) (23.3924,1.37678e-05) (23.6962,2.54157e-05) (24,4.28745e-05)};
\addplot [color=mplcolor4, dashed, mark=none, line width=1.00pt]
  coordinates {(0,0.125664) (0.3,0.114621) (0.6,0.104548) (0.9,0.095361) (1.2,0.086981) (1.5,0.0793374) (1.8,0.0723655) (2.15,0.0650021) (2.45,0.0592899) (2.75,0.0540797) (3.05,0.0493274) (3.35,0.0449927) (3.65,0.0410389) (3.95,0.0374325) (4.25,0.0341431) (4.55,0.0311427) (4.85,0.028406) (5.15,0.0259097) (5.45,0.0236329) (5.75,0.0215561) (6.1,0.0193627) (6.4,0.0176612) (6.7,0.0161092) (7,0.0146936) (7.3,0.0134023) (7.6,0.0122246) (7.9,0.0111503) (8.2,0.0101705) (8.5,0.00927673) (8.8,0.00846152) (9.1,0.00771795) (9.4,0.00703973) (9.7,0.0064211) (10.05,0.00576773) (10.35,0.00526089) (10.65,0.00479858) (10.95,0.00437689) (11.25,0.00399227) (11.55,0.00364144) (11.85,0.00332144) (12.15,0.00302957) (12.45,0.00276334) (12.75,0.00252051) (13.05,0.00229901) (13.35,0.00209698) (13.65,0.00191271) (13.95,0.00174462) (14.3,0.0015671) (14.6,0.00142939) (14.9,0.00130378) (15.2,0.00118921) (15.5,0.00108471) (15.8,0.000989386) (16.1,0.000902442) (16.4,0.000823139) (16.7,0.000750804) (17,0.000684826) (17.3,0.000624646) (17.6,0.000569754) (17.9,0.000519686) (18.25,0.000466806) (18.55,0.000425785) (18.85,0.000388368) (19.15,0.00035424) (19.45,0.000323111) (19.75,0.000294717) (20.05,0.000268818) (20.35,0.000245195) (20.65,0.000223648) (20.95,0.000203995) (21.25,0.000186068) (21.55,0.000169717) (21.85,0.000154803) (22.2,0.000139052) (22.5,0.000126832) (22.8,0.000115687) (23.1,0.00010552) (23.4,9.62477e-05) (23.7,8.77897e-05) (24,8.00751e-05)};
\addlegendentry{SL}
\addlegendentry{MLP flow}
\addlegendentry{SympNet flow}
\addlegendentry{Eulerian PINN}
\addlegendentry{theory, $\gamma=-0.1533$}
\end{axis}
\end{tikzpicture}%

    \end{minipage}\hfill
    \begin{minipage}[c]{0.43\linewidth}\centering
        \setlength{\figurewidth}{0.95\linewidth}
        \setlength{\figureheight}{0.5\linewidth}
        %
    \input{figures/vlasov_poisson/landau_damping_final_rho.tex}%
\\[6pt]
        %
    \input{figures/vlasov_poisson/landau_damping_final_E.tex}%

    \end{minipage}
    \caption{Landau damping from \cref{sec:vlasov_poisson}.
        Left: electric energy $\frac{1}{2}\int E^2\,\dx$ over time, with the theoretical damping rate $\gamma \approx -0.1533$.
        Right: density perturbation $\rho-1$ and electric field $E$
        at the final time $T = 24$.}
    \label{fig:vp-landau}
\end{figure}

\begin{figure}[!ht]
    \centering
    \pgfplotsset{every axis title/.append style={text height=1.6ex, text depth=0.25ex}}
    \input{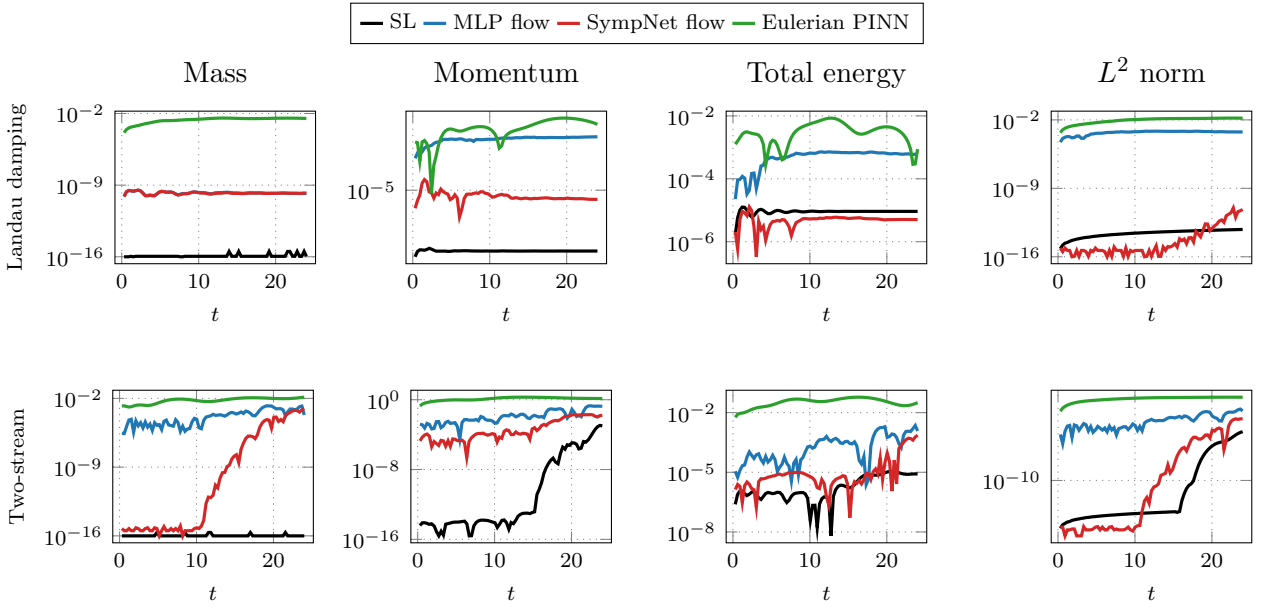}\ref*{vplegendlandaudamping}\\[2pt]
    \begin{minipage}[t]{0.25\linewidth}\centering
        \setlength{\figurewidth}{1\linewidth}
        \setlength{\figureheight}{0.85\linewidth}
        %
    \definecolor{mplcolor0}{rgb}{0.000000,0.000000,0.000000}
\definecolor{mplcolor1}{rgb}{0.121570,0.466670,0.705880}
\definecolor{mplcolor2}{rgb}{0.839220,0.152940,0.156860}
\definecolor{mplcolor3}{rgb}{0.172550,0.627450,0.172550}
\begin{tikzpicture}
\begin{axis}[
  width=\figurewidth,
  height=\figureheight,
  ymode=log, log basis y=10,
  xlabel={$t$},
  ylabel={Landau damping},
  title={Mass},
  xmin=-0.885, xmax=25.185,
  ymin=2.10075e-17, ymax=0.0169879,
  xmajorgrids=true,
  ymajorgrids=true,
  legend to name=vplegendlandaudamping,
  legend columns=4
]
\addplot [color=mplcolor0, mark=none, line width=1.20pt]
  coordinates {(0.3,1e-16) (0.6,1e-16) (0.9,1e-16) (1.2,1.11022e-16) (1.5,1e-16) (1.8,1.11022e-16) (2.15,1.11022e-16) (2.45,1.11022e-16) (2.75,1.11022e-16) (3.05,1.11022e-16) (3.35,1.11022e-16) (3.65,1.11022e-16) (3.95,1.11022e-16) (4.25,1.11022e-16) (4.55,1.11022e-16) (4.85,1.11022e-16) (5.15,1.11022e-16) (5.45,1.11022e-16) (5.75,1.11022e-16) (6.1,1.11022e-16) (6.4,1.11022e-16) (6.7,1.11022e-16) (7,1.11022e-16) (7.3,1.11022e-16) (7.6,1e-16) (7.9,1e-16) (8.2,1.11022e-16) (8.5,1.11022e-16) (8.8,1.11022e-16) (9.1,1.11022e-16) (9.4,1.11022e-16) (9.7,1.11022e-16) (10.05,1.11022e-16) (10.35,1.11022e-16) (10.65,1.11022e-16) (10.95,1.11022e-16) (11.25,1.11022e-16) (11.55,1.11022e-16) (11.85,1.11022e-16) (12.15,1.11022e-16) (12.45,1.11022e-16) (12.75,1.11022e-16) (13.05,1.11022e-16) (13.35,1.11022e-16) (13.65,1.11022e-16) (13.95,3.33067e-16) (14.3,1.11022e-16) (14.6,1.11022e-16) (14.9,1.11022e-16) (15.2,3.33067e-16) (15.5,1.11022e-16) (15.8,1.11022e-16) (16.1,1.11022e-16) (16.4,1.11022e-16) (16.7,1.11022e-16) (17,1.11022e-16) (17.3,1.11022e-16) (17.6,1.11022e-16) (17.9,1.11022e-16) (18.25,1.11022e-16) (18.55,1.11022e-16) (18.85,1.11022e-16) (19.15,3.33067e-16) (19.45,1.11022e-16) (19.75,1.11022e-16) (20.05,1.11022e-16) (20.35,1.11022e-16) (20.65,1.11022e-16) (20.95,1.11022e-16) (21.25,1.11022e-16) (21.55,3.33067e-16) (21.85,3.33067e-16) (22.2,1.11022e-16) (22.5,1.11022e-16) (22.8,3.33067e-16) (23.1,1.11022e-16) (23.4,1.11022e-16) (23.7,3.33067e-16) (24,1.11022e-16)};
\addplot [color=mplcolor1, mark=none, line width=1.20pt]
  coordinates {(0.303797,7.80145e-11) (0.607595,2.25503e-10) (0.911392,3.0385e-10) (1.21519,2.76817e-10) (1.51899,2.49754e-10) (1.82278,2.85465e-10) (2.12658,3.21439e-10) (2.43038,2.77717e-10) (2.73418,1.72337e-10) (3.03797,9.67648e-11) (3.34177,9.4996e-11) (3.64557,1.0606e-10) (3.94937,8.98734e-11) (4.25316,8.85515e-11) (4.55696,1.30182e-10) (4.86076,2.02925e-10) (5.16456,2.43654e-10) (5.46835,2.29435e-10) (5.77215,2.04897e-10) (6.07595,2.08067e-10) (6.37975,2.15403e-10) (6.68354,1.8832e-10) (6.98734,1.4046e-10) (7.29114,1.15233e-10) (7.59494,1.28916e-10) (7.89873,1.51668e-10) (8.20253,1.53574e-10) (8.50633,1.52514e-10) (8.81013,1.71355e-10) (9.11392,2.02837e-10) (9.41772,2.18999e-10) (9.72152,2.08351e-10) (10.0253,1.88701e-10) (10.3291,1.81853e-10) (10.6329,1.83359e-10) (10.9367,1.74796e-10) (11.2405,1.56542e-10) (11.5443,1.45283e-10) (11.8481,1.5439e-10) (12.1519,1.76103e-10) (12.4557,1.82207e-10) (12.7595,1.81921e-10) (13.0633,1.86434e-10) (13.3671,1.97517e-10) (13.6709,2.03161e-10) (13.9747,1.94872e-10) (14.2785,1.80383e-10) (14.5823,1.72355e-10) (14.8861,1.71398e-10) (15.1899,1.68609e-10) (15.4937,1.62903e-10) (15.7975,1.61148e-10) (16.1013,1.64848e-10) (16.4051,1.73082e-10) (16.7089,1.76421e-10) (17.0127,1.75178e-10) (17.3165,1.7393e-10) (17.6203,1.75158e-10) (17.9241,1.75696e-10) (18.2278,1.71414e-10) (18.5316,1.64648e-10) (18.8354,1.6112e-10) (19.1392,1.61276e-10) (19.443,1.60663e-10) (19.7468,1.58505e-10) (20.0506,1.59877e-10) (20.3544,1.61893e-10) (20.6582,1.66312e-10) (20.962,1.68408e-10) (21.2658,1.65954e-10) (21.5696,1.61682e-10) (21.8734,1.59219e-10) (22.1772,1.58328e-10) (22.481,1.56312e-10) (22.7848,1.53255e-10) (23.0886,1.52136e-10) (23.3924,1.53577e-10) (23.6962,1.55651e-10) (24,1.57295e-10)};
\addplot [color=mplcolor2, mark=none, line width=1.20pt]
  coordinates {(0.303797,7.62069e-11) (0.607595,2.25282e-10) (0.911392,3.03101e-10) (1.21519,2.79756e-10) (1.51899,2.55125e-10) (1.82278,2.91334e-10) (2.12658,3.27337e-10) (2.43038,2.79129e-10) (2.73418,1.69549e-10) (3.03797,9.48731e-11) (3.34177,9.08161e-11) (3.64557,1.07121e-10) (3.94937,9.9868e-11) (4.25316,9.60837e-11) (4.55696,1.39577e-10) (4.86076,2.11736e-10) (5.16456,2.50778e-10) (5.46835,2.35108e-10) (5.77215,2.08125e-10) (6.07595,2.0724e-10) (6.37975,2.1219e-10) (6.68354,1.85578e-10) (6.98734,1.35932e-10) (7.29114,1.07447e-10) (7.59494,1.16884e-10) (7.89873,1.36812e-10) (8.20253,1.41703e-10) (8.50633,1.42944e-10) (8.81013,1.63525e-10) (9.11392,1.96747e-10) (9.41772,2.12717e-10) (9.72152,1.9961e-10) (10.0253,1.78559e-10) (10.3291,1.70877e-10) (10.6329,1.7009e-10) (10.9367,1.58781e-10) (11.2405,1.39022e-10) (11.5443,1.30291e-10) (11.8481,1.40672e-10) (12.1519,1.56371e-10) (12.4557,1.62836e-10) (12.7595,1.63285e-10) (13.0633,1.69456e-10) (13.3671,1.81033e-10) (13.6709,1.85378e-10) (13.9747,1.76252e-10) (14.2785,1.62898e-10) (14.5823,1.56751e-10) (14.8861,1.56957e-10) (15.1899,1.55193e-10) (15.4937,1.49904e-10) (15.7975,1.48442e-10) (16.1013,1.55187e-10) (16.4051,1.64736e-10) (16.7089,1.69098e-10) (17.0127,1.67812e-10) (17.3165,1.66726e-10) (17.6203,1.68372e-10) (17.9241,1.68748e-10) (18.2278,1.63985e-10) (18.5316,1.56942e-10) (18.8354,1.53712e-10) (19.1392,1.55374e-10) (19.443,1.58105e-10) (19.7468,1.59038e-10) (20.0506,1.59373e-10) (20.3544,1.62079e-10) (20.6582,1.66442e-10) (20.962,1.68763e-10) (21.2658,1.66438e-10) (21.5696,1.62953e-10) (21.8734,1.60172e-10) (22.1772,1.60528e-10) (22.481,1.60833e-10) (22.7848,1.61053e-10) (23.0886,1.58022e-10) (23.3924,1.60628e-10) (23.6962,1.6741e-10) (24,1.62599e-10)};
\addplot [color=mplcolor3, mark=none, line width=1.20pt]
  coordinates {(0.303797,0.000134767) (0.607595,0.000271092) (0.911392,0.000386117) (1.21519,0.000473622) (1.51899,0.000537716) (1.82278,0.000596037) (2.12658,0.000669734) (2.43038,0.000767367) (2.73418,0.00088039) (3.03797,0.000994948) (3.34177,0.00110705) (3.64557,0.00122608) (3.94937,0.00136447) (4.25316,0.00152394) (4.55696,0.00168959) (4.86076,0.00183611) (5.16456,0.00194125) (5.46835,0.00199762) (5.77215,0.00201586) (6.07595,0.00201858) (6.37975,0.00202963) (6.68354,0.00206475) (6.98734,0.00212728) (7.29114,0.00220952) (7.59494,0.00229796) (7.89873,0.00237959) (8.20253,0.00244677) (8.50633,0.00249922) (8.81013,0.00254278) (9.11392,0.00258604) (9.41772,0.00263682) (9.72152,0.00269999) (10.0253,0.00277699) (10.3291,0.00286664) (10.6329,0.00296598) (10.9367,0.00307101) (11.2405,0.003177) (11.5443,0.00327879) (11.8481,0.00337117) (12.1519,0.0034494) (12.4557,0.0035097) (12.7595,0.0035497) (13.0633,0.00356873) (13.3671,0.0035677) (13.6709,0.003549) (13.9747,0.00351602) (14.2785,0.00347273) (14.5823,0.00342331) (14.8861,0.00337174) (15.1899,0.00332159) (15.4937,0.00327587) (15.7975,0.00323694) (16.1013,0.00320649) (16.4051,0.00318553) (16.7089,0.0031745) (17.0127,0.00317327) (17.3165,0.00318124) (17.6203,0.00319742) (17.9241,0.00322056) (18.2278,0.00324918) (18.5316,0.0032817) (18.8354,0.00331656) (19.1392,0.00335223) (19.443,0.00338729) (19.7468,0.00342047) (20.0506,0.00345066) (20.3544,0.00347693) (20.6582,0.00349847) (20.962,0.00351464) (21.2658,0.00352487) (21.5696,0.0035287) (21.8734,0.00352571) (22.1772,0.00351549) (22.481,0.00349765) (22.7848,0.00347176) (23.0886,0.00343734) (23.3924,0.00339386) (23.6962,0.00334074) (24,0.0032773)};
\addlegendentry{SL}
\addlegendentry{MLP flow}
\addlegendentry{SympNet flow}
\addlegendentry{Eulerian PINN}
\end{axis}
\end{tikzpicture}%

    \end{minipage}%
    \begin{minipage}[t]{0.25\linewidth}\centering
        \setlength{\figurewidth}{1\linewidth}
        \setlength{\figureheight}{0.85\linewidth}
        %
    \definecolor{mplcolor0}{rgb}{0.000000,0.000000,0.000000}
\definecolor{mplcolor1}{rgb}{0.121570,0.466670,0.705880}
\definecolor{mplcolor2}{rgb}{0.839220,0.152940,0.156860}
\definecolor{mplcolor3}{rgb}{0.172550,0.627450,0.172550}
\begin{tikzpicture}
\begin{axis}[
  width=\figurewidth,
  height=\figureheight,
  ymode=log, log basis y=10,
  xlabel={$t$},
  title={Momentum},
  xmin=-0.885, xmax=25.185,
  ymin=1.1372e-09, ymax=0.173804,
  xmajorgrids=true,
  ymajorgrids=true
]
\addplot [color=mplcolor0, mark=none, line width=1.20pt, forget plot]
  coordinates {(0.3,2.67823e-09) (0.6,5.10937e-09) (0.9,6.31105e-09) (1.2,6.24912e-09) (1.5,6.00942e-09) (1.8,6.51949e-09) (2.15,7.63834e-09) (2.45,6.82402e-09) (2.75,5.96538e-09) (3.05,5.36923e-09) (3.35,5.36482e-09) (3.65,5.37858e-09) (3.95,5.39808e-09) (4.25,5.3783e-09) (4.55,5.24096e-09) (4.85,5.21441e-09) (5.15,5.44289e-09) (5.45,5.59552e-09) (5.75,5.53516e-09) (6.1,5.53751e-09) (6.4,5.61325e-09) (6.7,5.64324e-09) (7,5.63217e-09) (7.3,5.58647e-09) (7.6,5.51519e-09) (7.9,5.47352e-09) (8.2,5.46935e-09) (8.5,5.45832e-09) (8.8,5.41574e-09) (9.1,5.41528e-09) (9.4,5.4767e-09) (9.7,5.46866e-09) (10.05,5.42261e-09) (10.35,5.42121e-09) (10.65,5.41658e-09) (10.95,5.41368e-09) (11.25,5.46294e-09) (11.55,5.53338e-09) (11.85,5.53018e-09) (12.15,5.48365e-09) (12.45,5.45733e-09) (12.75,5.45463e-09) (13.05,5.44706e-09) (13.35,5.42756e-09) (13.65,5.41604e-09) (13.95,5.4093e-09) (14.3,5.41959e-09) (14.6,5.4358e-09) (14.9,5.43895e-09) (15.2,5.44266e-09) (15.5,5.46183e-09) (15.8,5.48224e-09) (16.1,5.47512e-09) (16.4,5.45035e-09) (16.7,5.43405e-09) (17,5.43045e-09) (17.3,5.43135e-09) (17.6,5.42816e-09) (17.9,5.42373e-09) (18.25,5.42609e-09) (18.55,5.43846e-09) (18.85,5.45128e-09) (19.15,5.45239e-09) (19.45,5.44746e-09) (19.75,5.4463e-09) (20.05,5.44801e-09) (20.35,5.44447e-09) (20.65,5.43518e-09) (20.95,5.42865e-09) (21.25,5.42826e-09) (21.55,5.4324e-09) (21.85,5.43619e-09) (22.2,5.4368e-09) (22.5,5.43841e-09) (22.8,5.44389e-09) (23.1,5.4493e-09) (23.4,5.44747e-09) (23.7,5.4404e-09) (24,5.43565e-09)};
\addplot [color=mplcolor1, mark=none, line width=1.20pt, forget plot]
  coordinates {(0.303797,0.000515649) (0.607595,0.000977739) (0.911392,0.000843857) (1.21519,0.0011377) (1.51899,0.00182072) (1.82278,0.00221606) (2.12658,0.00228171) (2.43038,0.00254211) (2.73418,0.0031981) (3.03797,0.00372935) (3.34177,0.00385609) (3.64557,0.00430325) (3.94937,0.00529131) (4.25316,0.00535287) (4.55696,0.00525176) (4.86076,0.00551551) (5.16456,0.00547833) (5.46835,0.00530078) (5.77215,0.00539099) (6.07595,0.00567553) (6.37975,0.00561297) (6.68354,0.00503198) (6.98734,0.00466075) (7.29114,0.0049822) (7.59494,0.00489202) (7.89873,0.00414763) (8.20253,0.00521695) (8.50633,0.00522276) (8.81013,0.00508881) (9.11392,0.00512018) (9.41772,0.00527486) (9.72152,0.00537477) (10.0253,0.00540357) (10.3291,0.00539208) (10.6329,0.00539105) (10.9367,0.00556346) (11.2405,0.00594656) (11.5443,0.00622148) (11.8481,0.00629789) (12.1519,0.00612532) (12.4557,0.00627282) (12.7595,0.00612229) (13.0633,0.00629196) (13.3671,0.00656975) (13.6709,0.00653705) (13.9747,0.00645762) (14.2785,0.00665495) (14.5823,0.00690833) (14.8861,0.00692867) (15.1899,0.00683434) (15.4937,0.0068435) (15.7975,0.00696964) (16.1013,0.00656202) (16.4051,0.00664318) (16.7089,0.00663361) (17.0127,0.00656463) (17.3165,0.00659312) (17.6203,0.00669331) (17.9241,0.00674895) (18.2278,0.00671786) (18.5316,0.00666584) (18.8354,0.00665994) (19.1392,0.00667366) (19.443,0.00664464) (19.7468,0.00661636) (20.0506,0.00677436) (20.3544,0.0067887) (20.6582,0.00691134) (20.962,0.0069194) (21.2658,0.00688247) (21.5696,0.00690054) (21.8734,0.00698319) (22.1772,0.00707459) (22.481,0.00713255) (22.7848,0.00716609) (23.0886,0.00720528) (23.3924,0.00725302) (23.6962,0.00728734) (24,0.00730888)};
\addplot [color=mplcolor2, mark=none, line width=1.20pt, forget plot]
  coordinates {(0.303797,1.07333e-06) (0.607595,2.76283e-06) (0.911392,5.47096e-06) (1.21519,2.38988e-05) (1.51899,3.69212e-05) (1.82278,2.84172e-05) (2.12658,8.10232e-06) (2.43038,7.39916e-06) (2.73418,3.74498e-06) (3.03797,7.06105e-06) (3.34177,5.81843e-06) (3.64557,1.1361e-05) (3.94937,8.34784e-06) (4.25316,6.94364e-06) (4.55696,6.18388e-06) (4.86076,3.12035e-06) (5.16456,3.66925e-06) (5.46835,5.63824e-06) (5.77215,3.80611e-06) (6.07595,3.66649e-07) (6.37975,9.73459e-07) (6.68354,3.59121e-06) (6.98734,5.9761e-06) (7.29114,6.19941e-06) (7.59494,6.06235e-06) (7.89873,6.73246e-06) (8.20253,7.16546e-06) (8.50633,6.68709e-06) (8.81013,8.25652e-06) (9.11392,8.14541e-06) (9.41772,6.83672e-06) (9.72152,6.51897e-06) (10.0253,5.35854e-06) (10.3291,4.2059e-06) (10.6329,4.4216e-06) (10.9367,5.03458e-06) (11.2405,5.58222e-06) (11.5443,6.00649e-06) (11.8481,5.40386e-06) (12.1519,5.48839e-06) (12.4557,4.95048e-06) (12.7595,3.79829e-06) (13.0633,3.67884e-06) (13.3671,3.98511e-06) (13.6709,3.78091e-06) (13.9747,3.48843e-06) (14.2785,3.55046e-06) (14.5823,3.76238e-06) (14.8861,3.66155e-06) (15.1899,3.44687e-06) (15.4937,3.48317e-06) (15.7975,3.64901e-06) (16.1013,3.56342e-06) (16.4051,3.66188e-06) (16.7089,3.73962e-06) (17.0127,3.82275e-06) (17.3165,3.88629e-06) (17.6203,3.81211e-06) (17.9241,3.65787e-06) (18.2278,3.52421e-06) (18.5316,3.48717e-06) (18.8354,3.56556e-06) (19.1392,3.55882e-06) (19.443,3.51209e-06) (19.7468,3.53848e-06) (20.0506,3.77188e-06) (20.3544,3.63884e-06) (20.6582,3.54862e-06) (20.962,3.51727e-06) (21.2658,3.46773e-06) (21.5696,3.49101e-06) (21.8734,3.52423e-06) (22.1772,3.47237e-06) (22.481,3.36691e-06) (22.7848,3.28379e-06) (23.0886,3.24282e-06) (23.3924,3.22903e-06) (23.6962,3.23745e-06) (24,3.22688e-06)};
\addplot [color=mplcolor3, mark=none, line width=1.20pt, forget plot]
  coordinates {(0.303797,0.00388666) (0.607595,0.0035255) (0.911392,0.000364675) (1.21519,0.00387553) (1.51899,0.00466204) (1.82278,0.00312548) (2.12658,0.00110948) (2.43038,7.21216e-06) (2.73418,0.000220188) (3.03797,0.000706518) (3.34177,0.00261852) (3.64557,0.00627274) (3.94937,0.0109531) (4.25316,0.0153657) (4.55696,0.0183318) (4.86076,0.0193179) (5.16456,0.0185864) (5.46835,0.0169695) (5.77215,0.0154405) (6.07595,0.0147177) (6.37975,0.0150719) (6.68354,0.0163626) (6.98734,0.0182156) (7.29114,0.0202207) (7.59494,0.0220667) (7.89873,0.0235853) (8.20253,0.024718) (8.50633,0.0254412) (8.81013,0.0256932) (9.11392,0.0253419) (9.41772,0.0242137) (9.72152,0.0221668) (10.0253,0.0191701) (10.3291,0.0153448) (10.6329,0.0109496) (10.9367,0.00632143) (11.2405,0.00179882) (11.5443,0.0023437) (11.8481,0.00593447) (12.1519,0.00891341) (12.4557,0.0113166) (12.7595,0.0132478) (13.0633,0.0148474) (13.3671,0.0162674) (13.6709,0.0176552) (13.9747,0.019144) (14.2785,0.020849) (14.5823,0.0228644) (14.8861,0.025259) (15.1899,0.0280727) (15.4937,0.0313127) (15.7975,0.0349517) (16.1013,0.0389286) (16.4051,0.0431529) (16.7089,0.0475115) (17.0127,0.0518773) (17.3165,0.0561184) (17.6203,0.0601075) (17.9241,0.0637292) (18.2278,0.0668856) (18.5316,0.0695003) (18.8354,0.0715194) (19.1392,0.0729121) (19.443,0.0736686) (19.7468,0.0737987) (20.0506,0.0733282) (20.3544,0.072296) (20.6582,0.0707502) (20.962,0.0687454) (21.2658,0.0663385) (21.5696,0.0635864) (21.8734,0.0605428) (22.1772,0.0572567) (22.481,0.0537699) (22.7848,0.0501162) (23.0886,0.0463211) (23.3924,0.0424005) (23.6962,0.0383619) (24,0.0342044)};
\end{axis}
\end{tikzpicture}%

    \end{minipage}%
    \begin{minipage}[t]{0.25\linewidth}\centering
        \setlength{\figurewidth}{1\linewidth}
        \setlength{\figureheight}{0.85\linewidth}
        %
    \definecolor{mplcolor0}{rgb}{0.000000,0.000000,0.000000}
\definecolor{mplcolor1}{rgb}{0.121570,0.466670,0.705880}
\definecolor{mplcolor2}{rgb}{0.839220,0.152940,0.156860}
\definecolor{mplcolor3}{rgb}{0.172550,0.627450,0.172550}
\begin{tikzpicture}
\begin{axis}[
  width=\figurewidth,
  height=\figureheight,
  ymode=log, log basis y=10,
  xlabel={$t$},
  title={Total energy},
  xmin=-0.885, xmax=25.185,
  ymin=1.99227e-07, ymax=0.0141476,
  xmajorgrids=true,
  ymajorgrids=true
]
\addplot [color=mplcolor0, mark=none, line width=1.20pt, forget plot]
  coordinates {(0.3,1.64518e-06) (0.6,5.70099e-06) (0.9,1.00547e-05) (1.2,1.25972e-05) (1.5,1.23982e-05) (1.8,1.01268e-05) (2.15,7.14137e-06) (2.45,6.0591e-06) (2.75,6.80003e-06) (3.05,8.60093e-06) (3.35,1.01958e-05) (3.65,1.06604e-05) (3.95,9.92443e-06) (4.25,8.67629e-06) (4.55,7.80316e-06) (4.85,7.79955e-06) (5.15,8.52826e-06) (5.45,9.41793e-06) (5.75,9.89538e-06) (6.1,9.66933e-06) (6.4,9.07967e-06) (6.7,8.58426e-06) (7,8.48382e-06) (7.3,8.77648e-06) (7.6,9.21459e-06) (7.9,9.50456e-06) (8.2,9.49859e-06) (8.5,9.26018e-06) (8.8,8.98515e-06) (9.1,8.85474e-06) (9.4,8.92581e-06) (9.7,9.11952e-06) (10.05,9.31563e-06) (10.35,9.34946e-06) (10.65,9.26191e-06) (10.95,9.1323e-06) (11.25,9.04971e-06) (11.55,9.05669e-06) (11.85,9.13276e-06) (12.15,9.21969e-06) (12.45,9.2644e-06) (12.75,9.24986e-06) (13.05,9.19777e-06) (13.35,9.14758e-06) (13.65,9.1293e-06) (13.95,9.1478e-06) (14.3,9.1916e-06) (14.6,9.21998e-06) (14.9,9.22476e-06) (15.2,9.20859e-06) (15.5,9.18544e-06) (15.8,9.17022e-06) (16.1,9.17014e-06) (16.4,9.18239e-06) (16.7,9.19788e-06) (17,9.20764e-06) (17.3,9.20769e-06) (17.6,9.20012e-06) (17.9,9.19065e-06) (18.25,9.18446e-06) (18.55,9.18581e-06) (18.85,9.19161e-06) (19.15,9.19811e-06) (19.45,9.20188e-06) (19.75,9.20148e-06) (20.05,9.19784e-06) (20.35,9.1934e-06) (20.65,9.19062e-06) (20.95,9.19076e-06) (21.25,9.19341e-06) (21.55,9.19693e-06) (21.85,9.1994e-06) (22.2,9.19947e-06) (22.5,9.19735e-06) (22.8,9.19458e-06) (23.1,9.19273e-06) (23.4,9.19276e-06) (23.7,9.19453e-06) (24,9.19697e-06)};
\addplot [color=mplcolor1, mark=none, line width=1.20pt, forget plot]
  coordinates {(0.303797,2.26616e-05) (0.607595,9.08704e-05) (0.911392,9.56784e-05) (1.21519,0.000119862) (1.51899,0.000105024) (1.82278,3.27665e-05) (2.12658,4.30698e-05) (2.43038,0.000147557) (2.73418,0.000152028) (3.03797,4.8726e-05) (3.34177,0.000105214) (3.64557,0.000255257) (3.94937,0.000269478) (4.25316,0.000413423) (4.55696,0.000427492) (4.86076,0.000475567) (5.16456,0.000481872) (5.46835,0.000442098) (5.77215,0.000425347) (6.07595,0.000452267) (6.37975,0.000480782) (6.68354,0.000484166) (6.98734,0.000484587) (7.29114,0.000513181) (7.59494,0.000557194) (7.89873,0.0005839) (8.20253,0.000650973) (8.50633,0.000632534) (8.81013,0.000616159) (9.11392,0.00059964) (9.41772,0.000605677) (9.72152,0.000637746) (10.0253,0.00065906) (10.3291,0.000650184) (10.6329,0.000642254) (10.9367,0.000668435) (11.2405,0.000710281) (11.5443,0.000719605) (11.8481,0.000696374) (12.1519,0.000704115) (12.4557,0.000716999) (12.7595,0.000708878) (13.0633,0.000699687) (13.3671,0.00069906) (13.6709,0.00069403) (13.9747,0.000684107) (14.2785,0.0006803) (14.5823,0.000681419) (14.8861,0.000677359) (15.1899,0.000672413) (15.4937,0.000677626) (15.7975,0.000680445) (16.1013,0.000686432) (16.4051,0.000683078) (16.7089,0.000663497) (17.0127,0.000658602) (17.3165,0.000673906) (17.6203,0.000696692) (17.9241,0.000709082) (18.2278,0.000701872) (18.5316,0.000681789) (18.8354,0.000665524) (19.1392,0.000663194) (19.443,0.000666369) (19.7468,0.000660631) (20.0506,0.000640032) (20.3544,0.000620675) (20.6582,0.000622614) (20.962,0.000636112) (21.2658,0.000635) (21.5696,0.000618221) (21.8734,0.000604451) (22.1772,0.000606624) (22.481,0.000617398) (22.7848,0.000619071) (23.0886,0.000606432) (23.3924,0.000596614) (23.6962,0.000606904) (24,0.000607833)};
\addplot [color=mplcolor2, mark=none, line width=1.20pt, forget plot]
  coordinates {(0.303797,2.03891e-06) (0.607595,5.93e-07) (0.911392,5.17075e-06) (1.21519,9.28865e-06) (1.51899,8.2739e-06) (1.82278,6.70321e-06) (2.12658,1.27029e-05) (2.43038,9.72486e-06) (2.73418,2.95801e-06) (3.03797,3.31026e-07) (3.34177,4.40985e-06) (3.64557,2.4416e-06) (3.94937,3.10013e-06) (4.25316,1.1194e-06) (4.55696,2.21681e-06) (4.86076,3.56417e-06) (5.16456,4.36789e-06) (5.46835,4.46077e-06) (5.77215,4.75851e-06) (6.07595,4.49314e-06) (6.37975,5.0273e-06) (6.68354,4.92796e-06) (6.98734,4.046e-06) (7.29114,2.98364e-06) (7.59494,1.55728e-06) (7.89873,2.70175e-06) (8.20253,4.02458e-06) (8.50633,4.42207e-06) (8.81013,4.87124e-06) (9.11392,5.39716e-06) (9.41772,5.42368e-06) (9.72152,5.33571e-06) (10.0253,5.24933e-06) (10.3291,5.33981e-06) (10.6329,5.35868e-06) (10.9367,5.3237e-06) (11.2405,5.59865e-06) (11.5443,5.51772e-06) (11.8481,5.56256e-06) (12.1519,5.84489e-06) (12.4557,5.75506e-06) (12.7595,5.73322e-06) (13.0633,5.94411e-06) (13.3671,6.0133e-06) (13.6709,5.89682e-06) (13.9747,5.79606e-06) (14.2785,5.80206e-06) (14.5823,5.70616e-06) (14.8861,5.52421e-06) (15.1899,5.48326e-06) (15.4937,5.58482e-06) (15.7975,5.65763e-06) (16.1013,5.69624e-06) (16.4051,5.7972e-06) (16.7089,5.61974e-06) (17.0127,5.46082e-06) (17.3165,5.39747e-06) (17.6203,5.34213e-06) (17.9241,5.34479e-06) (18.2278,5.34108e-06) (18.5316,5.33243e-06) (18.8354,5.28088e-06) (19.1392,5.15868e-06) (19.443,5.02869e-06) (19.7468,4.93415e-06) (20.0506,4.99189e-06) (20.3544,5.02789e-06) (20.6582,5.02124e-06) (20.962,5.02425e-06) (21.2658,5.05133e-06) (21.5696,5.0698e-06) (21.8734,5.07658e-06) (22.1772,5.08247e-06) (22.481,5.08514e-06) (22.7848,5.07486e-06) (23.0886,5.06543e-06) (23.3924,5.06412e-06) (23.6962,5.04301e-06) (24,5.01256e-06)};
\addplot [color=mplcolor3, mark=none, line width=1.20pt, forget plot]
  coordinates {(0.303797,0.00126347) (0.607595,0.0015822) (0.911392,0.002059) (1.21519,0.0026346) (1.51899,0.00302006) (1.82278,0.00309401) (2.12658,0.0029535) (2.43038,0.00278092) (2.73418,0.00268311) (3.03797,0.00261855) (3.34177,0.0024398) (3.64557,0.00199934) (3.94937,0.00124869) (4.25316,0.00028433) (4.55696,0.000677905) (4.86076,0.00139089) (5.16456,0.00168701) (5.46835,0.00155471) (5.77215,0.0011409) (6.07595,0.000681643) (6.37975,0.000401825) (6.68354,0.000434143) (6.98734,0.0007894) (7.29114,0.00137971) (7.59494,0.00207298) (7.89873,0.00274965) (8.20253,0.00333896) (8.50633,0.0038259) (8.81013,0.00423362) (9.11392,0.00459615) (9.41772,0.0049376) (9.72152,0.005267) (10.0253,0.00558619) (10.3291,0.00590098) (10.6329,0.00622586) (10.9367,0.00657892) (11.2405,0.00697048) (11.5443,0.00739205) (11.8481,0.00781149) (12.1519,0.00817712) (12.4557,0.00842916) (12.7595,0.00851464) (13.0633,0.00840046) (13.3671,0.00808076) (13.6709,0.00757724) (13.9747,0.0069334) (14.2785,0.00620545) (14.5823,0.00545266) (14.8861,0.00472938) (15.1899,0.00407977) (15.4937,0.0035353) (15.7975,0.00311466) (16.1013,0.00282496) (16.4051,0.00266384) (16.7089,0.00262159) (17.0127,0.00268306) (17.3165,0.00282932) (17.6203,0.00303903) (17.9241,0.0032896) (18.2278,0.00355835) (18.5316,0.00382348) (18.8354,0.00406503) (19.1392,0.00426559) (19.443,0.00441086) (19.7468,0.00448997) (20.0506,0.00449555) (20.3544,0.00442365) (20.6582,0.00427347) (20.962,0.004047) (21.2658,0.00374859) (21.5696,0.00338443) (21.8734,0.00296213) (22.1772,0.0024902) (22.481,0.0019777) (22.7848,0.0014338) (23.0886,0.000867561) (23.3924,0.000287671) (23.6962,0.000297696) (24,0.000880959)};
\end{axis}
\end{tikzpicture}%

    \end{minipage}%
    \begin{minipage}[t]{0.25\linewidth}\centering
        \setlength{\figurewidth}{1\linewidth}
        \setlength{\figureheight}{0.85\linewidth}
        %
    \definecolor{mplcolor0}{rgb}{0.000000,0.000000,0.000000}
\definecolor{mplcolor1}{rgb}{0.121570,0.466670,0.705880}
\definecolor{mplcolor2}{rgb}{0.839220,0.152940,0.156860}
\definecolor{mplcolor3}{rgb}{0.172550,0.627450,0.172550}
\begin{tikzpicture}
\begin{axis}[
  width=\figurewidth,
  height=\figureheight,
  ymode=log, log basis y=10,
  xlabel={$t$},
  title={$L^2$ norm},
  xmin=-0.885, xmax=25.185,
  ymin=1.95317e-17, ymax=0.0784223,
  xmajorgrids=true,
  ymajorgrids=true
]
\addplot [color=mplcolor0, mark=none, line width=1.20pt, forget plot]
  coordinates {(0.3,6.66134e-16) (0.6,1.55431e-15) (0.9,2.44249e-15) (1.2,2.88658e-15) (1.5,3.77476e-15) (1.8,4.66294e-15) (2.15,5.55112e-15) (2.45,6.43929e-15) (2.75,7.10543e-15) (3.05,7.99361e-15) (3.35,8.88178e-15) (3.65,9.54792e-15) (3.95,9.99201e-15) (4.25,1.08802e-14) (4.55,1.17684e-14) (4.85,1.26565e-14) (5.15,1.33227e-14) (5.45,1.42109e-14) (5.75,1.4877e-14) (6.1,1.57652e-14) (6.4,1.64313e-14) (6.7,1.73195e-14) (7,1.84297e-14) (7.3,1.88738e-14) (7.6,1.95399e-14) (7.9,2.04281e-14) (8.2,2.10942e-14) (8.5,2.17604e-14) (8.8,2.28706e-14) (9.1,2.35367e-14) (9.4,2.4647e-14) (9.7,2.53131e-14) (10.05,2.62013e-14) (10.35,2.68674e-14) (10.65,2.75335e-14) (10.95,2.84217e-14) (11.25,2.93099e-14) (11.55,2.9976e-14) (11.85,3.06422e-14) (12.15,3.13083e-14) (12.45,3.21965e-14) (12.75,3.28626e-14) (13.05,3.37508e-14) (13.35,3.4639e-14) (13.65,3.53051e-14) (13.95,3.59712e-14) (14.3,3.68594e-14) (14.6,3.75255e-14) (14.9,3.81917e-14) (15.2,3.90799e-14) (15.5,3.9968e-14) (15.8,4.04121e-14) (16.1,4.15223e-14) (16.4,4.21885e-14) (16.7,4.30767e-14) (17,4.37428e-14) (17.3,4.44089e-14) (17.6,4.52971e-14) (17.9,4.64073e-14) (18.25,4.70735e-14) (18.55,4.79616e-14) (18.85,4.88498e-14) (19.15,4.92939e-14) (19.45,5.01821e-14) (19.75,5.12923e-14) (20.05,5.17364e-14) (20.35,5.24025e-14) (20.65,5.32907e-14) (20.95,5.39568e-14) (21.25,5.50671e-14) (21.55,5.57332e-14) (21.85,5.63993e-14) (22.2,5.75096e-14) (22.5,5.81757e-14) (22.8,5.90639e-14) (23.1,5.973e-14) (23.4,6.06182e-14) (23.7,6.12843e-14) (24,6.21725e-14)};
\addplot [color=mplcolor1, mark=none, line width=1.20pt, forget plot]
  coordinates {(0.303797,5.37335e-05) (0.607595,0.000126964) (0.911392,0.000169225) (1.21519,0.000197766) (1.51899,0.000189713) (1.82278,0.000147737) (2.12658,0.000158003) (2.43038,0.000225564) (2.73418,0.000223548) (3.03797,0.000130624) (3.34177,0.000131084) (3.64557,0.000226444) (3.94937,0.000256126) (4.25316,0.000353519) (4.55696,0.000368082) (4.86076,0.00041264) (5.16456,0.0004362) (5.46835,0.000429402) (5.77215,0.000438457) (6.07595,0.000479285) (6.37975,0.000517144) (6.68354,0.00052518) (6.98734,0.000523799) (7.29114,0.000547119) (7.59494,0.000579137) (7.89873,0.000572203) (8.20253,0.000612664) (8.50633,0.000593556) (8.81013,0.000575636) (9.11392,0.000567364) (9.41772,0.000576963) (9.72152,0.000599495) (10.0253,0.000614834) (10.3291,0.000617147) (10.6329,0.000624064) (10.9367,0.000651569) (11.2405,0.000688798) (11.5443,0.00070231) (11.8481,0.00068447) (12.1519,0.000689127) (12.4557,0.000701343) (12.7595,0.000703274) (13.0633,0.000696501) (13.3671,0.000690172) (13.6709,0.000683416) (13.9747,0.000676794) (14.2785,0.000675074) (14.5823,0.000677799) (14.8861,0.000678333) (15.1899,0.000675474) (15.4937,0.000675323) (15.7975,0.000675974) (16.1013,0.000684074) (16.4051,0.000688246) (16.7089,0.000668945) (17.0127,0.000655741) (17.3165,0.000663387) (17.6203,0.000684784) (17.9241,0.000700624) (18.2278,0.00069687) (18.5316,0.000676375) (18.8354,0.000655291) (19.1392,0.000647284) (19.443,0.000649765) (19.7468,0.000647341) (20.0506,0.000622841) (20.3544,0.000605967) (20.6582,0.000608886) (20.962,0.000620453) (21.2658,0.000621896) (21.5696,0.00061027) (21.8734,0.000596181) (22.1772,0.000589755) (22.481,0.00059089) (22.7848,0.000591739) (23.0886,0.0005875) (23.3924,0.000583427) (23.6962,0.000587017) (24,0.00058562)};
\addplot [color=mplcolor2, mark=none, line width=1.20pt, forget plot]
  coordinates {(0.303797,7.77156e-16) (0.607595,3.33067e-16) (0.911392,7.77156e-16) (1.21519,4.44089e-16) (1.51899,4.44089e-16) (1.82278,1e-16) (2.12658,1.11022e-16) (2.43038,4.44089e-16) (2.73418,1.11022e-16) (3.03797,4.44089e-16) (3.34177,4.44089e-16) (3.64557,3.33067e-16) (3.94937,1.11022e-16) (4.25316,4.44089e-16) (4.55696,3.33067e-16) (4.86076,7.77156e-16) (5.16456,4.44089e-16) (5.46835,4.44089e-16) (5.77215,4.44089e-16) (6.07595,3.33067e-16) (6.37975,1e-16) (6.68354,4.44089e-16) (6.98734,1.11022e-16) (7.29114,9.99201e-16) (7.59494,4.44089e-16) (7.89873,4.44089e-16) (8.20253,4.44089e-16) (8.50633,1.11022e-16) (8.81013,7.77156e-16) (9.11392,1.11022e-16) (9.41772,4.44089e-16) (9.72152,3.33067e-16) (10.0253,4.44089e-16) (10.3291,4.44089e-16) (10.6329,4.44089e-16) (10.9367,4.44089e-16) (11.2405,1e-16) (11.5443,1.11022e-16) (11.8481,4.44089e-16) (12.1519,1.11022e-16) (12.4557,4.44089e-16) (12.7595,1.33227e-15) (13.0633,1.11022e-16) (13.3671,1.11022e-16) (13.6709,9.99201e-16) (13.9747,4.44089e-16) (14.2785,1.11022e-15) (14.5823,9.99201e-16) (14.8861,1.33227e-15) (15.1899,2.22045e-15) (15.4937,9.99201e-16) (15.7975,6.88338e-15) (16.1013,4.44089e-15) (16.4051,6.43929e-15) (16.7089,2.22045e-15) (17.0127,8.77076e-15) (17.3165,9.21485e-15) (17.6203,3.21965e-15) (17.9241,4.03011e-14) (18.2278,2.09832e-14) (18.5316,1.08802e-14) (18.8354,4.84057e-14) (19.1392,2.30038e-13) (19.443,2.21601e-13) (19.7468,4.59632e-14) (20.0506,8.57092e-14) (20.3544,1.45439e-13) (20.6582,2.22933e-13) (20.962,4.54636e-13) (21.2658,5.85199e-13) (21.5696,3.33622e-13) (21.8734,1.57141e-12) (22.1772,2.25886e-12) (22.481,3.14027e-12) (22.7848,9.8056e-12) (23.0886,1.09979e-12) (23.3924,4.1912e-12) (23.6962,5.07183e-12) (24,6.89004e-12)};
\addplot [color=mplcolor3, mark=none, line width=1.20pt, forget plot]
  coordinates {(0.303797,0.000495555) (0.607595,0.000981951) (0.911392,0.00140763) (1.21519,0.00177401) (1.51899,0.00209788) (1.82278,0.00239144) (2.12658,0.00266203) (2.43038,0.0029198) (2.73418,0.00317811) (3.03797,0.00344369) (3.34177,0.00371031) (3.64557,0.00396756) (3.94937,0.00421717) (4.25316,0.00447755) (4.55696,0.00477058) (4.86076,0.00510313) (5.16456,0.00545998) (5.46835,0.00581287) (5.77215,0.00613826) (6.07595,0.00643195) (6.37975,0.00671265) (6.68354,0.00701294) (6.98734,0.00736282) (7.29114,0.00777449) (7.59494,0.00823572) (7.89873,0.00871451) (8.20253,0.00917185) (8.50633,0.00957616) (8.81013,0.00991296) (9.11392,0.0101868) (9.41772,0.0104163) (9.72152,0.0106254) (10.0253,0.0108352) (10.3291,0.0110588) (10.6329,0.0112999) (10.9367,0.0115543) (11.2405,0.0118138) (11.5443,0.0120685) (11.8481,0.0123105) (12.1519,0.0125341) (12.4557,0.0127364) (12.7595,0.0129167) (13.0633,0.0130757) (13.3671,0.0132147) (13.6709,0.013335) (13.9747,0.013438) (14.2785,0.0135249) (14.5823,0.0135967) (14.8861,0.0136549) (15.1899,0.0137011) (15.4937,0.0137372) (15.7975,0.0137656) (16.1013,0.0137889) (16.4051,0.0138099) (16.7089,0.0138313) (17.0127,0.0138558) (17.3165,0.0138858) (17.6203,0.0139233) (17.9241,0.0139699) (18.2278,0.0140266) (18.5316,0.014094) (18.8354,0.0141721) (19.1392,0.0142603) (19.443,0.0143576) (19.7468,0.0144625) (20.0506,0.014573) (20.3544,0.0146869) (20.6582,0.0148012) (20.962,0.014913) (21.2658,0.0150189) (21.5696,0.015115) (21.8734,0.0151973) (22.1772,0.0152616) (22.481,0.0153032) (22.7848,0.0153172) (23.0886,0.0152986) (23.3924,0.0152417) (23.6962,0.015141) (24,0.0149905)};
\end{axis}
\end{tikzpicture}%

    \end{minipage}\\[4pt]
    \begin{minipage}[t]{0.25\linewidth}\centering
        \setlength{\figurewidth}{1\linewidth}
        \setlength{\figureheight}{0.85\linewidth}
        %
    \definecolor{mplcolor0}{rgb}{0.000000,0.000000,0.000000}
\definecolor{mplcolor1}{rgb}{0.121570,0.466670,0.705880}
\definecolor{mplcolor2}{rgb}{0.839220,0.152940,0.156860}
\definecolor{mplcolor3}{rgb}{0.172550,0.627450,0.172550}
\begin{tikzpicture}
\begin{axis}[
  width=\figurewidth,
  height=\figureheight,
  ymode=log, log basis y=10,
  xlabel={$t$},
  ylabel={Two-stream},
  xmin=-0.885, xmax=25.185,
  ymin=1.96731e-17, ymax=0.0673999,
  xmajorgrids=true,
  ymajorgrids=true,
  legend to name=vplegendtwostream,
  legend columns=4
]
\addplot [color=mplcolor0, mark=none, line width=1.20pt]
  coordinates {(0.3,1e-16) (0.6,1e-16) (0.9,1e-16) (1.2,1e-16) (1.5,1e-16) (1.8,1e-16) (2.15,1e-16) (2.45,1e-16) (2.75,1e-16) (3.05,1e-16) (3.35,1e-16) (3.65,1e-16) (3.95,1e-16) (4.25,1e-16) (4.55,1e-16) (4.85,2.22045e-16) (5.15,1e-16) (5.45,1e-16) (5.75,1e-16) (6.1,1e-16) (6.4,1e-16) (6.7,1e-16) (7,1e-16) (7.3,1e-16) (7.6,1e-16) (7.9,1e-16) (8.2,1e-16) (8.5,2.22045e-16) (8.8,1e-16) (9.1,1e-16) (9.4,1e-16) (9.7,1e-16) (10.05,1e-16) (10.35,1e-16) (10.65,1e-16) (10.95,1e-16) (11.25,1e-16) (11.55,2.22045e-16) (11.85,2.22045e-16) (12.15,1e-16) (12.45,1e-16) (12.75,1e-16) (13.05,1e-16) (13.35,1e-16) (13.65,1e-16) (13.95,1e-16) (14.3,1e-16) (14.6,1e-16) (14.9,1e-16) (15.2,1e-16) (15.5,1e-16) (15.8,1e-16) (16.1,1e-16) (16.4,1e-16) (16.7,1e-16) (17,2.22045e-16) (17.3,1e-16) (17.6,1e-16) (17.9,1e-16) (18.25,1e-16) (18.55,1e-16) (18.85,1e-16) (19.15,1e-16) (19.45,1e-16) (19.75,1e-16) (20.05,1e-16) (20.35,1e-16) (20.65,1e-16) (20.95,1e-16) (21.25,1e-16) (21.55,2.22045e-16) (21.85,1e-16) (22.2,1e-16) (22.5,1e-16) (22.8,1e-16) (23.1,1e-16) (23.4,1e-16) (23.7,1e-16) (24,1e-16)};
\addplot [color=mplcolor1, mark=none, line width=1.20pt]
  coordinates {(0.303797,2.81349e-06) (0.607595,2.62988e-06) (0.911392,1.39852e-05) (1.21519,6.32557e-05) (1.51899,7.44274e-05) (1.82278,3.69735e-06) (2.12658,5.11993e-05) (2.43038,1.62336e-05) (2.73418,3.32825e-05) (3.03797,9.85595e-06) (3.34177,3.99415e-05) (3.64557,1.19543e-05) (3.94937,2.39974e-05) (4.25316,5.50637e-06) (4.55696,3.34544e-05) (4.86076,1.77842e-05) (5.16456,2.14487e-05) (5.46835,6.66628e-06) (5.77215,4.61121e-05) (6.07595,6.16625e-05) (6.37975,1.92109e-05) (6.68354,1.72962e-05) (6.98734,1.48088e-05) (7.29114,6.70016e-05) (7.59494,3.52858e-05) (7.89873,4.42522e-06) (8.20253,1.47684e-05) (8.50633,4.45754e-06) (8.81013,7.56801e-06) (9.11392,8.31451e-06) (9.41772,4.03144e-05) (9.72152,9.30109e-05) (10.0253,7.64816e-05) (10.3291,1.00199e-05) (10.6329,8.10247e-06) (10.9367,5.77299e-05) (11.2405,0.000126041) (11.5443,0.000146173) (11.8481,0.000168603) (12.1519,0.000238725) (12.4557,0.000188013) (12.7595,0.000194661) (13.0633,0.000216644) (13.3671,0.000186652) (13.6709,0.000155877) (13.9747,0.0002059) (14.2785,0.000314546) (14.5823,0.000370204) (14.8861,0.00031902) (15.1899,0.00025515) (15.4937,0.000313745) (15.7975,0.000414014) (16.1013,0.000292041) (16.4051,0.000513593) (16.7089,0.000433011) (17.0127,0.000293582) (17.3165,0.000177024) (17.6203,0.00016831) (17.9241,0.000310118) (18.2278,0.000624204) (18.5316,0.00108988) (18.8354,0.0014583) (19.1392,0.00177987) (19.443,0.00167567) (19.7468,0.00167624) (20.0506,0.00121455) (20.3544,0.000853658) (20.6582,0.000891467) (20.962,0.000894269) (21.2658,0.000654353) (21.5696,0.0001615) (21.8734,0.00045752) (22.1772,0.000833859) (22.481,0.000816817) (22.7848,0.000681049) (23.0886,0.00084067) (23.3924,0.00143421) (23.6962,0.00164681) (24,0.000203628)};
\addplot [color=mplcolor2, mark=none, line width=1.20pt]
  coordinates {(0.303797,5.55112e-16) (0.607595,3.33067e-16) (0.911392,3.33067e-16) (1.21519,5.55112e-16) (1.51899,3.33067e-16) (1.82278,3.33067e-16) (2.12658,6.66134e-16) (2.43038,5.55112e-16) (2.73418,3.33067e-16) (3.03797,3.33067e-16) (3.34177,5.55112e-16) (3.64557,6.66134e-16) (3.94937,3.33067e-16) (4.25316,3.33067e-16) (4.55696,6.66134e-16) (4.86076,5.55112e-16) (5.16456,1e-16) (5.46835,3.33067e-16) (5.77215,8.88178e-16) (6.07595,3.33067e-16) (6.37975,3.33067e-16) (6.68354,5.55112e-16) (6.98734,3.33067e-16) (7.29114,5.55112e-16) (7.59494,3.33067e-16) (7.89873,9.99201e-16) (8.20253,1e-16) (8.50633,3.33067e-16) (8.81013,3.33067e-16) (9.11392,3.33067e-16) (9.41772,3.33067e-16) (9.72152,5.55112e-16) (10.0253,6.66134e-16) (10.3291,5.55112e-16) (10.6329,9.99201e-16) (10.9367,3.77476e-15) (11.2405,1.99396e-13) (11.5443,8.894e-13) (11.8481,7.77933e-13) (12.1519,1.34481e-12) (12.4557,1.46367e-11) (12.7595,7.8042e-12) (13.0633,7.56935e-11) (13.3671,2.68452e-10) (13.6709,2.88081e-10) (13.9747,1.14709e-09) (14.2785,1.60045e-09) (14.5823,6.38017e-09) (14.8861,1.2966e-08) (15.1899,1.80476e-08) (15.4937,2.91495e-09) (15.7975,6.90614e-08) (16.1013,2.78412e-07) (16.4051,8.28583e-07) (16.7089,1.4516e-06) (17.0127,1.44187e-06) (17.3165,1.24554e-06) (17.6203,7.99635e-07) (17.9241,2.00304e-06) (18.2278,1.18615e-05) (18.5316,3.07375e-05) (18.8354,5.36973e-05) (19.1392,6.942e-05) (19.443,5.9892e-05) (19.7468,9.42307e-06) (20.0506,0.000121185) (20.3544,0.000125491) (20.6582,4.52144e-05) (20.962,7.04347e-05) (21.2658,4.21343e-05) (21.5696,0.000180527) (21.8734,0.000278117) (22.1772,0.000371958) (22.481,0.000516379) (22.7848,0.000619551) (23.0886,0.000564097) (23.3924,0.000401651) (23.6962,0.000666744) (24,0.000892952)};
\addplot [color=mplcolor3, mark=none, line width=1.20pt]
  coordinates {(0.303797,0.00163686) (0.607595,0.00175532) (0.911392,0.00147006) (1.21519,0.00131966) (1.51899,0.00142562) (1.82278,0.00168932) (2.12658,0.0019564) (2.43038,0.00211655) (2.73418,0.00213871) (3.03797,0.00206077) (3.34177,0.00195784) (3.64557,0.00190963) (3.94937,0.00197858) (4.25316,0.00220119) (4.55696,0.00258835) (4.86076,0.00312969) (5.16456,0.00379816) (5.46835,0.00455427) (5.77215,0.00535049) (6.07595,0.00613611) (6.37975,0.00686244) (6.68354,0.00748741) (6.98734,0.00797893) (7.29114,0.00831667) (7.59494,0.00849237) (7.89873,0.00850902) (8.20253,0.00837921) (8.50633,0.0081232) (8.81013,0.00776684) (9.11392,0.00733957) (9.41772,0.00687254) (9.72152,0.00639688) (10.0253,0.00594219) (10.3291,0.00553513) (10.6329,0.00519835) (10.9367,0.00494962) (11.2405,0.00480135) (11.5443,0.00476042) (11.8481,0.00482845) (12.1519,0.00500224) (12.4557,0.00527447) (12.7595,0.00563446) (13.0633,0.00606896) (13.3671,0.00656285) (13.6709,0.00709979) (13.9747,0.00766288) (14.2785,0.00823516) (14.5823,0.00880013) (14.8861,0.00934225) (15.1899,0.00984737) (15.4937,0.0103031) (15.7975,0.0106992) (16.1013,0.0110277) (16.4051,0.0112833) (16.7089,0.011463) (17.0127,0.0115668) (17.3165,0.011597) (17.6203,0.0115581) (17.9241,0.0114569) (18.2278,0.011302) (18.5316,0.0111034) (18.8354,0.0108724) (19.1392,0.010621) (19.443,0.010362) (19.7468,0.0101084) (20.0506,0.009873) (20.3544,0.00966861) (20.6582,0.0095075) (20.962,0.00940132) (21.2658,0.00936098) (21.5696,0.00939648) (21.8734,0.00951686) (22.1772,0.00973007) (22.481,0.0100429) (22.7848,0.0104611) (23.0886,0.0109891) (23.3924,0.0116303) (23.6962,0.0123868) (24,0.0132597)};
\addlegendentry{SL}
\addlegendentry{MLP flow}
\addlegendentry{SympNet flow}
\addlegendentry{Eulerian PINN}
\end{axis}
\end{tikzpicture}%

    \end{minipage}%
    \begin{minipage}[t]{0.25\linewidth}\centering
        \setlength{\figurewidth}{1\linewidth}
        \setlength{\figureheight}{0.85\linewidth}
        %
    \definecolor{mplcolor0}{rgb}{0.000000,0.000000,0.000000}
\definecolor{mplcolor1}{rgb}{0.121570,0.466670,0.705880}
\definecolor{mplcolor2}{rgb}{0.839220,0.152940,0.156860}
\definecolor{mplcolor3}{rgb}{0.172550,0.627450,0.172550}
\begin{tikzpicture}
\begin{axis}[
  width=\figurewidth,
  height=\figureheight,
  ymode=log, log basis y=10,
  xlabel={$t$},
  xmin=-0.885, xmax=25.185,
  ymin=4.00706e-17, ymax=12.2182,
  xmajorgrids=true,
  ymajorgrids=true
]
\addplot [color=mplcolor0, mark=none, line width=1.20pt, forget plot]
  coordinates {(0.3,3.84665e-15) (0.6,7.84316e-15) (0.9,8.04299e-15) (1.2,6.84403e-15) (1.5,5.24543e-15) (1.8,6.06971e-15) (2.15,3.34708e-15) (2.45,7.49347e-16) (2.75,2.49782e-16) (3.05,8.4926e-16) (3.35,2.99739e-16) (3.65,5.49521e-15) (3.95,6.34447e-15) (4.25,7.44351e-15) (4.55,8.84229e-15) (4.85,9.44177e-15) (5.15,9.4043e-15) (5.45,9.84142e-15) (5.75,1.09405e-14) (6.1,1.20395e-14) (6.4,5.74499e-15) (6.7,2.49782e-16) (7,2.62271e-16) (7.3,3.447e-15) (7.6,8.54255e-15) (7.9,1.31385e-14) (8.2,1.56364e-14) (8.5,1.59486e-14) (8.8,1.51368e-14) (9.1,1.44374e-14) (9.4,1.19396e-14) (9.7,9.74151e-15) (10.05,8.1429e-15) (10.35,2.298e-15) (10.65,5.5077e-15) (10.95,1.50369e-14) (11.25,1.64856e-14) (11.55,8.09294e-15) (11.85,1.44874e-15) (12.15,8.44264e-15) (12.45,6.24456e-15) (12.75,1.50369e-14) (13.05,4.50108e-14) (13.35,7.25867e-14) (13.65,9.35684e-14) (13.95,9.85641e-14) (14.3,1.00163e-13) (14.6,8.85478e-14) (14.9,1.18946e-13) (15.2,1.43525e-13) (15.5,1.34676e-11) (15.8,5.44322e-11) (16.1,3.00498e-10) (16.4,3.55897e-09) (16.7,1.60616e-08) (17,4.85901e-08) (17.3,1.1761e-07) (17.6,2.01639e-07) (17.9,7.48501e-08) (18.25,4.60331e-08) (18.55,6.61031e-07) (18.85,3.66908e-06) (19.15,2.15122e-06) (19.45,6.88367e-06) (19.75,5.86611e-06) (20.05,1.16968e-05) (20.35,1.1058e-05) (20.65,7.84702e-06) (20.95,2.31566e-06) (21.25,2.80024e-06) (21.55,1.61908e-05) (21.85,1.37387e-05) (22.2,2.65718e-05) (22.5,6.70272e-05) (22.8,0.000176582) (23.1,0.00043328) (23.4,0.000840283) (23.7,0.00113822) (24,0.00100638)};
\addplot [color=mplcolor1, mark=none, line width=1.20pt, forget plot]
  coordinates {(0.303797,0.0015577) (0.607595,0.00133584) (0.911392,0.000570295) (1.21519,0.00194117) (1.51899,0.0014896) (1.82278,0.000453473) (2.12658,0.000512667) (2.43038,0.00303771) (2.73418,0.00317929) (3.03797,0.000624937) (3.34177,0.000987101) (3.64557,0.00255526) (3.94937,0.00335811) (4.25316,0.00194591) (4.55696,0.00212917) (4.86076,0.0023333) (5.16456,0.000845173) (5.46835,3.01237e-05) (5.77215,0.00252319) (6.07595,0.0053071) (6.37975,0.00526743) (6.68354,0.00451031) (6.98734,0.00457222) (7.29114,0.00450392) (7.59494,0.0055041) (7.89873,0.00381654) (8.20253,0.0119409) (8.50633,0.00959596) (8.81013,0.00538683) (9.11392,0.00436707) (9.41772,0.00456264) (9.72152,0.00364172) (10.0253,0.00197819) (10.3291,0.00218752) (10.6329,0.00633748) (10.9367,0.0123837) (11.2405,0.0154822) (11.5443,0.0147209) (11.8481,0.0140161) (12.1519,0.0229078) (12.4557,0.0175208) (12.7595,0.0108102) (13.0633,0.0139316) (13.3671,0.0209817) (13.6709,0.0222412) (13.9747,0.0157803) (14.2785,0.00813546) (14.5823,0.00692994) (14.8861,0.0124429) (15.1899,0.0163091) (15.4937,0.012691) (15.7975,0.0103324) (16.1013,0.0120894) (16.4051,0.00894229) (16.7089,0.0377837) (17.0127,0.0766035) (17.3165,0.090436) (17.6203,0.0822816) (17.9241,0.0769463) (18.2278,0.0868787) (18.5316,0.0930216) (18.8354,0.0743328) (19.1392,0.0310405) (19.443,0.00769777) (19.7468,0.0163234) (20.0506,0.0797233) (20.3544,0.0357439) (20.6582,0.0120041) (20.962,0.0158641) (21.2658,0.0112067) (21.5696,0.0516421) (21.8734,0.13183) (22.1772,0.192932) (22.481,0.202141) (22.7848,0.184476) (23.0886,0.179403) (23.3924,0.187038) (23.6962,0.180085) (24,0.178626)};
\addplot [color=mplcolor2, mark=none, line width=1.20pt, forget plot]
  coordinates {(0.303797,1.93322e-05) (0.607595,5.19962e-05) (0.911392,9.12051e-05) (1.21519,0.00011567) (1.51899,6.74597e-05) (1.82278,4.15047e-06) (2.12658,1.60745e-05) (2.43038,9.82904e-06) (2.73418,5.6974e-06) (3.03797,2.19886e-05) (3.34177,1.93855e-05) (3.64557,4.73242e-06) (3.94937,1.96882e-05) (4.25316,4.42945e-06) (4.55696,3.22446e-05) (4.86076,2.95767e-05) (5.16456,2.06121e-05) (5.46835,2.05054e-05) (5.77215,2.31308e-05) (6.07595,1.35116e-05) (6.37975,2.41077e-07) (6.68354,1.64384e-05) (6.98734,4.2986e-05) (7.29114,7.51265e-05) (7.59494,0.000117613) (7.89873,0.000136211) (8.20253,0.000134469) (8.50633,0.0001475) (8.81013,0.000106743) (9.11392,8.66142e-05) (9.41772,0.000105944) (9.72152,0.000131638) (10.0253,0.000108523) (10.3291,1.41289e-05) (10.6329,5.31558e-05) (10.9367,3.8144e-05) (11.2405,0.000236779) (11.5443,0.000357847) (11.8481,0.000291016) (12.1519,0.000101637) (12.4557,0.000100368) (12.7595,0.000203595) (13.0633,0.000129506) (13.3671,0.000111641) (13.6709,0.000143758) (13.9747,0.000134596) (14.2785,0.00011504) (14.5823,2.96413e-05) (14.8861,0.000447196) (15.1899,0.000620379) (15.4937,0.000293257) (15.7975,0.000476234) (16.1013,0.00110458) (16.4051,0.00147723) (16.7089,0.00109181) (17.0127,0.00201939) (17.3165,0.00282642) (17.6203,0.0036719) (17.9241,0.00426053) (18.2278,0.00382616) (18.5316,0.0051263) (18.8354,0.00797309) (19.1392,0.0123604) (19.443,0.0129225) (19.7468,0.0135657) (20.0506,0.0183885) (20.3544,0.0194508) (20.6582,0.0176201) (20.962,0.0176944) (21.2658,0.0190418) (21.5696,0.0187381) (21.8734,0.0191888) (22.1772,0.0192857) (22.481,0.0199946) (22.7848,0.0178853) (23.0886,0.0120385) (23.3924,0.0104983) (23.6962,0.0153959) (24,0.0162479)};
\addplot [color=mplcolor3, mark=none, line width=1.20pt, forget plot]
  coordinates {(0.303797,0.192854) (0.607595,0.314438) (0.911392,0.405439) (1.21519,0.482752) (1.51899,0.551456) (1.82278,0.612405) (2.12658,0.666048) (2.43038,0.713727) (2.73418,0.75747) (3.03797,0.799049) (3.34177,0.839067) (3.64557,0.8766) (3.94937,0.909579) (4.25316,0.93569) (4.55696,0.953311) (4.86076,0.962109) (5.16456,0.963128) (5.46835,0.958497) (5.77215,0.950961) (6.07595,0.943432) (6.37975,0.938651) (6.68354,0.938982) (6.98734,0.946311) (7.29114,0.962018) (7.59494,0.987001) (7.89873,1.02172) (8.20253,1.06621) (8.50633,1.12012) (8.81013,1.18265) (9.11392,1.25256) (9.41772,1.3281) (9.72152,1.4071) (10.0253,1.48709) (10.3291,1.56548) (10.6329,1.63981) (10.9367,1.708) (11.2405,1.76846) (11.5443,1.82025) (11.8481,1.86303) (12.1519,1.89698) (12.4557,1.92268) (12.7595,1.94092) (13.0633,1.95264) (13.3671,1.95874) (13.6709,1.96007) (13.9747,1.95734) (14.2785,1.95116) (14.5823,1.94199) (14.8861,1.93016) (15.1899,1.91592) (15.4937,1.89943) (15.7975,1.8808) (16.1013,1.86012) (16.4051,1.83747) (16.7089,1.81296) (17.0127,1.78674) (17.3165,1.75902) (17.6203,1.73006) (17.9241,1.7002) (18.2278,1.66982) (18.5316,1.63936) (18.8354,1.6093) (19.1392,1.5801) (19.443,1.55226) (19.7468,1.52622) (20.0506,1.50238) (20.3544,1.48106) (20.6582,1.46247) (20.962,1.44672) (21.2658,1.43377) (21.5696,1.42342) (21.8734,1.41529) (22.1772,1.40882) (22.481,1.40325) (22.7848,1.39761) (23.0886,1.39074) (23.3924,1.38123) (23.6962,1.36749) (24,1.34771)};
\end{axis}
\end{tikzpicture}%

    \end{minipage}%
    \begin{minipage}[t]{0.25\linewidth}\centering
        \setlength{\figurewidth}{1\linewidth}
        \setlength{\figureheight}{0.85\linewidth}
        %
    \definecolor{mplcolor0}{rgb}{0.000000,0.000000,0.000000}
\definecolor{mplcolor1}{rgb}{0.121570,0.466670,0.705880}
\definecolor{mplcolor2}{rgb}{0.839220,0.152940,0.156860}
\definecolor{mplcolor3}{rgb}{0.172550,0.627450,0.172550}
\begin{tikzpicture}
\begin{axis}[
  width=\figurewidth,
  height=\figureheight,
  ymode=log, log basis y=10,
  xlabel={$t$},
  xmin=-0.885, xmax=25.185,
  ymin=2.84144e-09, ymax=0.133307,
  xmajorgrids=true,
  ymajorgrids=true
]
\addplot [color=mplcolor0, mark=none, line width=1.20pt, forget plot]
  coordinates {(0.3,2.44727e-07) (0.6,7.68827e-07) (0.9,1.14377e-06) (1.2,1.1243e-06) (1.5,8.51813e-07) (1.8,6.70405e-07) (2.15,7.98673e-07) (2.45,1.01466e-06) (2.75,1.03611e-06) (3.05,7.96419e-07) (3.35,5.06593e-07) (3.65,4.27709e-07) (3.95,6.07451e-07) (4.25,8.57083e-07) (4.55,9.66034e-07) (4.85,9.10879e-07) (5.15,8.4277e-07) (5.45,8.85902e-07) (5.75,9.86261e-07) (6.1,9.54488e-07) (6.4,7.20942e-07) (6.7,4.33506e-07) (7,2.98124e-07) (7.3,3.95407e-07) (7.6,6.22204e-07) (7.9,8.12998e-07) (8.2,8.96482e-07) (8.5,9.23345e-07) (8.8,9.54301e-07) (9.1,9.5359e-07) (9.4,8.16146e-07) (9.7,4.96972e-07) (10.05,3.6654e-08) (10.35,2.16829e-07) (10.65,2.19444e-07) (10.95,1.68565e-08) (11.25,2.5128e-07) (11.55,4.66575e-07) (11.85,5.82216e-07) (12.15,5.78959e-07) (12.45,4.07708e-07) (12.75,6.34211e-09) (13.05,6.17917e-07) (13.35,1.33283e-06) (13.65,1.92673e-06) (13.95,2.23116e-06) (14.3,2.19311e-06) (14.6,1.96228e-06) (14.9,1.73516e-06) (15.2,1.66592e-06) (15.5,1.86618e-06) (15.8,2.40929e-06) (16.1,3.31105e-06) (16.4,4.48735e-06) (16.7,5.73882e-06) (17,6.80605e-06) (17.3,7.48042e-06) (17.6,7.70089e-06) (17.9,7.57333e-06) (18.25,7.26757e-06) (18.55,7.12213e-06) (18.85,7.25065e-06) (19.15,7.71378e-06) (19.45,8.48445e-06) (19.75,9.42718e-06) (20.05,1.03471e-05) (20.35,1.0999e-05) (20.65,1.1226e-05) (20.95,1.09691e-05) (21.25,1.03399e-05) (21.55,9.53142e-06) (21.85,8.77542e-06) (22.2,8.1551e-06) (22.5,7.90941e-06) (22.8,7.92451e-06) (23.1,8.08111e-06) (23.4,8.23834e-06) (23.7,8.29532e-06) (24,8.19618e-06)};
\addplot [color=mplcolor1, mark=none, line width=1.20pt, forget plot]
  coordinates {(0.303797,1.1184e-05) (0.607595,5.5922e-06) (0.911392,1.24414e-05) (1.21519,2.12999e-05) (1.51899,2.32245e-05) (1.82278,5.81798e-05) (2.12658,0.000105405) (2.43038,7.10949e-05) (2.73418,4.94933e-05) (3.03797,7.59955e-05) (3.34177,8.51217e-05) (3.64557,4.39508e-05) (3.94937,2.93837e-05) (4.25316,6.28303e-05) (4.55696,3.27262e-05) (4.86076,4.28995e-05) (5.16456,7.15404e-05) (5.46835,4.34873e-05) (5.77215,1.10122e-05) (6.07595,1.79542e-05) (6.37975,2.55693e-05) (6.68354,5.73954e-05) (6.98734,2.87639e-05) (7.29114,1.96605e-05) (7.59494,3.62198e-06) (7.89873,2.87036e-05) (8.20253,4.40052e-05) (8.50633,7.42941e-06) (8.81013,4.20188e-05) (9.11392,5.71513e-05) (9.41772,0.000110743) (9.72152,0.000190186) (10.0253,0.000206247) (10.3291,0.000161975) (10.6329,0.000185181) (10.9367,0.000321695) (11.2405,0.000428655) (11.5443,0.00043672) (11.8481,0.0004617) (12.1519,0.000630603) (12.4557,0.000427741) (12.7595,0.000350964) (13.0633,0.000450817) (13.3671,0.000505769) (13.6709,0.000416488) (13.9747,0.000294945) (14.2785,0.000281065) (14.5823,0.000349658) (14.8861,0.000358559) (15.1899,0.000278432) (15.4937,0.000268244) (15.7975,0.000339102) (16.1013,0.000209939) (16.4051,3.95031e-05) (16.7089,2.61119e-05) (17.0127,6.10481e-05) (17.3165,2.39464e-06) (17.6203,4.59146e-06) (17.9241,0.000210679) (18.2278,0.000610935) (18.5316,0.00111274) (18.8354,0.00150895) (19.1392,0.00176139) (19.443,0.00177531) (19.7468,0.00194553) (20.0506,0.000564747) (20.3544,0.000629235) (20.6582,0.000845424) (20.962,0.00079218) (21.2658,0.000562957) (21.5696,0.000279084) (21.8734,0.000175621) (22.1772,0.00073754) (22.481,0.00112611) (22.7848,0.0011649) (23.0886,0.00112088) (23.3924,0.00157184) (23.6962,0.00227157) (24,0.00118715)};
\addplot [color=mplcolor2, mark=none, line width=1.20pt, forget plot]
  coordinates {(0.303797,1.35508e-06) (0.607595,2.43803e-06) (0.911392,1.94941e-06) (1.21519,3.01328e-07) (1.51899,2.12789e-06) (1.82278,3.22978e-06) (2.12658,2.47728e-06) (2.43038,1.50719e-06) (2.73418,1.32501e-06) (3.03797,1.76985e-07) (3.34177,3.36433e-06) (3.64557,4.75197e-06) (3.94937,4.65392e-06) (4.25316,5.05479e-06) (4.55696,5.64626e-06) (4.86076,6.00749e-06) (5.16456,6.49931e-06) (5.46835,6.77782e-06) (5.77215,7.79399e-06) (6.07595,8.11605e-06) (6.37975,7.18302e-06) (6.68354,7.47096e-06) (6.98734,8.70922e-06) (7.29114,9.58795e-06) (7.59494,9.61043e-06) (7.89873,9.5322e-06) (8.20253,9.87627e-06) (8.50633,9.91656e-06) (8.81013,9.62419e-06) (9.11392,8.97981e-06) (9.41772,7.92923e-06) (9.72152,6.94094e-06) (10.0253,6.12806e-06) (10.3291,5.24579e-06) (10.6329,3.87128e-06) (10.9367,1.71539e-06) (11.2405,1.50755e-06) (11.5443,4.31408e-06) (11.8481,2.77379e-06) (12.1519,5.95969e-07) (12.4557,2.1932e-07) (12.7595,2.10599e-06) (13.0633,3.35123e-06) (13.3671,4.38284e-06) (13.6709,5.04653e-06) (13.9747,5.5834e-06) (14.2785,6.73435e-06) (14.5823,6.23252e-06) (14.8861,2.59822e-06) (15.1899,5.22509e-08) (15.4937,1.53821e-06) (15.7975,3.39028e-06) (16.1013,2.78268e-06) (16.4051,3.81247e-06) (16.7089,2.61416e-06) (17.0127,4.22088e-06) (17.3165,8.65393e-06) (17.6203,8.83776e-06) (17.9241,4.35461e-06) (18.2278,1.07257e-05) (18.5316,1.47637e-05) (18.8354,8.26466e-06) (19.1392,5.094e-05) (19.443,4.31146e-06) (19.7468,2.95627e-05) (20.0506,9.17939e-05) (20.3544,3.3868e-05) (20.6582,1.13052e-06) (20.962,4.06569e-05) (21.2658,2.57842e-06) (21.5696,0.000156382) (21.8734,0.000237782) (22.1772,0.000318122) (22.481,0.000431958) (22.7848,0.000517029) (23.0886,0.000497345) (23.3924,0.00041371) (23.6962,0.000600701) (24,0.000728871)};
\addplot [color=mplcolor3, mark=none, line width=1.20pt, forget plot]
  coordinates {(0.303797,0.00577642) (0.607595,0.0082642) (0.911392,0.00911112) (1.21519,0.00952224) (1.51899,0.0101717) (1.82278,0.0112739) (2.12658,0.0127417) (2.43038,0.0143547) (2.73418,0.0158898) (3.03797,0.017198) (3.34177,0.0182331) (3.64557,0.0190463) (3.94937,0.0197626) (4.25316,0.0205468) (4.55696,0.0215675) (4.86076,0.0229634) (5.16456,0.0248161) (5.46835,0.0271361) (5.77215,0.0298616) (6.07595,0.0328697) (6.37975,0.0359958) (6.68354,0.0390563) (6.98734,0.0418698) (7.29114,0.0442746) (7.59494,0.0461417) (7.89873,0.0473832) (8.20253,0.0479572) (8.50633,0.047869) (8.81013,0.0471708) (9.11392,0.0459583) (9.41772,0.0443637) (9.72152,0.0425467) (10.0253,0.0406809) (10.3291,0.0389388) (10.6329,0.0374765) (10.9367,0.0364202) (11.2405,0.0358558) (11.5443,0.0358253) (11.8481,0.0363278) (12.1519,0.0373261) (12.4557,0.0387557) (12.7595,0.0405351) (13.0633,0.0425742) (13.3671,0.0447814) (13.6709,0.0470681) (13.9747,0.0493511) (14.2785,0.0515536) (14.5823,0.0536053) (14.8861,0.055443) (15.1899,0.0570102) (15.4937,0.0582581) (15.7975,0.0591459) (16.1013,0.0596422) (16.4051,0.0597252) (16.7089,0.0593844) (17.0127,0.0586204) (17.3165,0.057445) (17.6203,0.0558814) (17.9241,0.0539628) (18.2278,0.0517318) (18.5316,0.0492391) (18.8354,0.0465419) (19.1392,0.0437023) (19.443,0.0407865) (19.7468,0.0378629) (20.0506,0.0350014) (20.3544,0.0322718) (20.6582,0.0297433) (20.962,0.0274833) (21.2658,0.0255568) (21.5696,0.0240254) (21.8734,0.0229464) (22.1772,0.0223724) (22.481,0.0223498) (22.7848,0.0229188) (23.0886,0.0241119) (23.3924,0.0259538) (23.6962,0.0284602) (24,0.0316375)};
\end{axis}
\end{tikzpicture}%

    \end{minipage}%
    \begin{minipage}[t]{0.25\linewidth}\centering
        \setlength{\figurewidth}{1\linewidth}
        \setlength{\figureheight}{0.85\linewidth}
        %
    \definecolor{mplcolor0}{rgb}{0.000000,0.000000,0.000000}
\definecolor{mplcolor1}{rgb}{0.121570,0.466670,0.705880}
\definecolor{mplcolor2}{rgb}{0.839220,0.152940,0.156860}
\definecolor{mplcolor3}{rgb}{0.172550,0.627450,0.172550}
\begin{tikzpicture}
\begin{axis}[
  width=\figurewidth,
  height=\figureheight,
  ymode=log, log basis y=10,
  xlabel={$t$},
  xmin=-0.885, xmax=25.185,
  ymin=1.77646e-17, ymax=0.574586,
  xmajorgrids=true,
  ymajorgrids=true
]
\addplot [color=mplcolor0, mark=none, line width=1.20pt, forget plot]
  coordinates {(0.3,6.66134e-16) (0.6,1.55431e-15) (0.9,2.22045e-15) (1.2,3.10862e-15) (1.5,3.55271e-15) (1.8,4.66294e-15) (2.15,5.55112e-15) (2.45,6.21725e-15) (2.75,6.88338e-15) (3.05,7.77156e-15) (3.35,8.43769e-15) (3.65,9.54792e-15) (3.95,9.99201e-15) (4.25,1.08802e-14) (4.55,1.17684e-14) (4.85,1.22125e-14) (5.15,1.31006e-14) (5.45,1.37668e-14) (5.75,1.46549e-14) (6.1,1.55431e-14) (6.4,1.62093e-14) (6.7,1.70974e-14) (7,1.77636e-14) (7.3,1.82077e-14) (7.6,1.93179e-14) (7.9,1.9984e-14) (8.2,2.08722e-14) (8.5,2.15383e-14) (8.8,2.24265e-14) (9.1,2.33147e-14) (9.4,2.42029e-14) (9.7,2.4869e-14) (10.05,2.57572e-14) (10.35,2.64233e-14) (10.65,2.75335e-14) (10.95,2.79776e-14) (11.25,2.88658e-14) (11.55,2.9754e-14) (11.85,3.04201e-14) (12.15,3.13083e-14) (12.45,3.21965e-14) (12.75,3.30846e-14) (13.05,3.39728e-14) (13.35,3.4861e-14) (13.65,3.59712e-14) (13.95,3.66374e-14) (14.3,3.79696e-14) (14.6,3.90799e-14) (14.9,3.9968e-14) (15.2,4.08562e-14) (15.5,3.53051e-14) (15.8,4.29656e-14) (16.1,6.48259e-13) (16.4,3.63776e-12) (16.7,1.28407e-11) (17,2.89214e-11) (17.3,4.18141e-11) (17.6,8.69229e-11) (17.9,5.75184e-10) (18.25,4.06309e-09) (18.55,1.47636e-08) (18.85,4.14755e-08) (19.15,9.65156e-08) (19.45,1.94058e-07) (19.75,3.4521e-07) (20.05,5.50375e-07) (20.35,7.91254e-07) (20.65,1.03591e-06) (20.95,1.25438e-06) (21.25,1.44284e-06) (21.55,1.63712e-06) (21.85,1.92448e-06) (22.2,2.6009e-06) (22.5,3.7468e-06) (22.8,5.66714e-06) (23.1,8.42595e-06) (23.4,1.18569e-05) (23.7,1.56125e-05) (24,1.92211e-05)};
\addplot [color=mplcolor1, mark=none, line width=1.20pt, forget plot]
  coordinates {(0.303797,9.51488e-06) (0.607595,9.52839e-07) (0.911392,2.10346e-05) (1.21519,7.58956e-05) (1.51899,9.00469e-05) (1.82278,4.53528e-06) (2.12658,0.000100428) (2.43038,8.2347e-05) (2.73418,2.65525e-05) (3.03797,6.20287e-05) (3.34177,0.000143126) (3.64557,0.000106184) (3.94937,2.25674e-05) (4.25316,6.80091e-05) (4.55696,1.59039e-05) (4.86076,2.84667e-05) (5.16456,7.10797e-05) (5.46835,4.53833e-05) (5.77215,2.79655e-05) (6.07595,4.69417e-05) (6.37975,1.80949e-05) (6.68354,7.50057e-05) (6.98734,2.91793e-05) (7.29114,5.07384e-05) (7.59494,8.36592e-06) (7.89873,4.10741e-05) (8.20253,4.90823e-05) (8.50633,6.48207e-05) (8.81013,6.37999e-05) (9.11392,3.18433e-05) (9.41772,5.8309e-05) (9.72152,0.000129629) (10.0253,8.92712e-05) (10.3291,2.5561e-05) (10.6329,7.00552e-05) (10.9367,2.35226e-05) (11.2405,0.000161987) (11.5443,0.000254116) (11.8481,0.000309757) (12.1519,0.000411032) (12.4557,0.00030235) (12.7595,0.000276968) (13.0633,0.00030743) (13.3671,0.000299481) (13.6709,0.000297676) (13.9747,0.000395154) (14.2785,0.00056405) (14.5823,0.000663756) (14.8861,0.000625915) (15.1899,0.000578806) (15.4937,0.000704613) (15.7975,0.00085126) (16.1013,0.000588331) (16.4051,0.000955928) (16.7089,0.000867639) (17.0127,0.000612634) (17.3165,0.000359847) (17.6203,0.000267931) (17.9241,0.000480228) (18.2278,0.000994547) (18.5316,0.00172975) (18.8354,0.00240031) (19.1392,0.0027717) (19.443,0.00269659) (19.7468,0.00234239) (20.0506,0.000903289) (20.3544,0.000638248) (20.6582,0.000528749) (20.962,0.000276852) (21.2658,0.000271257) (21.5696,0.00109371) (21.8734,0.00193785) (22.1772,0.00243923) (22.481,0.0024488) (22.7848,0.00240592) (23.0886,0.00304187) (23.3924,0.00445096) (23.6962,0.00529069) (24,0.00328708)};
\addplot [color=mplcolor2, mark=none, line width=1.20pt, forget plot]
  coordinates {(0.303797,8.88178e-16) (0.607595,8.88178e-16) (0.911392,1e-16) (1.21519,4.44089e-16) (1.51899,4.44089e-16) (1.82278,5.55112e-16) (2.12658,5.55112e-16) (2.43038,1e-16) (2.73418,1e-16) (3.03797,1e-16) (3.34177,9.99201e-16) (3.64557,4.44089e-16) (3.94937,4.44089e-16) (4.25316,4.44089e-16) (4.55696,5.55112e-16) (4.86076,8.88178e-16) (5.16456,9.99201e-16) (5.46835,4.44089e-16) (5.77215,8.88178e-16) (6.07595,4.44089e-16) (6.37975,4.44089e-16) (6.68354,8.88178e-16) (6.98734,9.99201e-16) (7.29114,4.44089e-16) (7.59494,5.55112e-16) (7.89873,4.44089e-16) (8.20253,4.44089e-16) (8.50633,4.44089e-16) (8.81013,5.55112e-16) (9.11392,8.88178e-16) (9.41772,4.44089e-16) (9.72152,3.33067e-16) (10.0253,8.88178e-16) (10.3291,1.22125e-15) (10.6329,1.33227e-15) (10.9367,4.68514e-14) (11.2405,6.01741e-13) (11.5443,1.41231e-12) (11.8481,3.28404e-13) (12.1519,7.69518e-12) (12.4557,1.96554e-11) (12.7595,2.4744e-11) (13.0633,4.17084e-11) (13.3671,3.46598e-10) (13.6709,4.59444e-10) (13.9747,5.90187e-10) (14.2785,4.67177e-09) (14.5823,6.51058e-10) (14.8861,1.43877e-09) (15.1899,8.49373e-09) (15.4937,4.57335e-08) (15.7975,1.30985e-07) (16.1013,2.85506e-07) (16.4051,3.70408e-07) (16.7089,1.55331e-07) (17.0127,6.42923e-07) (17.3165,2.48245e-06) (17.6203,4.68552e-06) (17.9241,4.72344e-06) (18.2278,3.03892e-06) (18.5316,2.91278e-06) (18.8354,4.13553e-06) (19.1392,8.07343e-06) (19.443,8.48613e-06) (19.7468,3.88588e-06) (20.0506,2.79325e-05) (20.3544,2.1902e-05) (20.6582,3.9006e-05) (20.962,4.66339e-05) (21.2658,3.46848e-05) (21.5696,2.16123e-07) (21.8734,2.23363e-06) (22.1772,7.90931e-05) (22.481,0.000239625) (22.7848,0.000374998) (23.0886,0.000459936) (23.3924,0.000400581) (23.6962,0.000456048) (24,0.000479413)};
\addplot [color=mplcolor3, mark=none, line width=1.20pt, forget plot]
  coordinates {(0.303797,0.00316405) (0.607595,0.00669632) (0.911392,0.0103155) (1.21519,0.0139694) (1.51899,0.0176521) (1.82278,0.0213246) (2.12658,0.0249088) (2.43038,0.0283133) (2.73418,0.0314623) (3.03797,0.0343177) (3.34177,0.0368869) (3.64557,0.0392174) (3.94937,0.0413802) (4.25316,0.0434486) (4.55696,0.0454815) (4.86076,0.0475146) (5.16456,0.0495619) (5.46835,0.0516221) (5.77215,0.0536868) (6.07595,0.055746) (6.37975,0.0577908) (6.68354,0.0598142) (6.98734,0.0618102) (7.29114,0.0637736) (7.59494,0.0656996) (7.89873,0.0675835) (8.20253,0.0694206) (8.50633,0.0712073) (8.81013,0.0729402) (9.11392,0.0746169) (9.41772,0.0762354) (9.72152,0.0777947) (10.0253,0.0792938) (10.3291,0.0807321) (10.6329,0.082109) (10.9367,0.0834238) (11.2405,0.0846758) (11.5443,0.0858636) (11.8481,0.0869862) (12.1519,0.0880421) (12.4557,0.08903) (12.7595,0.0899486) (13.0633,0.0907973) (13.3671,0.0915757) (13.6709,0.0922845) (13.9747,0.0929254) (14.2785,0.0935009) (14.5823,0.0940151) (14.8861,0.0944728) (15.1899,0.0948804) (15.4937,0.095245) (15.7975,0.0955745) (16.1013,0.0958773) (16.4051,0.096162) (16.7089,0.096437) (17.0127,0.09671) (17.3165,0.0969881) (17.6203,0.0972771) (17.9241,0.0975813) (18.2278,0.0979036) (18.5316,0.098245) (18.8354,0.0986049) (19.1392,0.0989807) (19.443,0.0993684) (19.7468,0.099762) (20.0506,0.100155) (20.3544,0.100538) (20.6582,0.100902) (20.962,0.101236) (21.2658,0.101531) (21.5696,0.101773) (21.8734,0.101952) (22.1772,0.102056) (22.481,0.102073) (22.7848,0.10199) (23.0886,0.101795) (23.3924,0.101477) (23.6962,0.101025) (24,0.100426)};
\end{axis}
\end{tikzpicture}%

    \end{minipage}
    \caption{Vlasov-Poisson tests from \cref{sec:vlasov_poisson}:
        drift of the invariants for Landau damping (top) and the two-stream
        instability (bottom), evaluated on a $256\times 256$ phase-space grid.
        Relative drift $|Q(t)/Q(0)-1|$ of the mass, total energy and squared
        $L^2$ norm, absolute drift $|P(t)-P(0)|$ of the momentum, which vanishes
        at $t = 0$.}
    \label{fig:vp-conservation}
\end{figure}

Both learned flows follow the reference and the theoretical damping rate up
to $T = 24$ (\cref{fig:vp-landau}); the SympNet flow is indistinguishable
from the reference, while the MLP flow deviates slightly at the end. The Eulerian PINN, with seven times
more parameters, follows only the first oscillations of the electric energy, smooths
out the filaments, and its error on $f$ is about eight times larger than that
of the MLP flow (not shown).
The mass is conserved up to $10^{-10}$ by both flows
(\cref{fig:vp-conservation}, top). For the SympNet flow, which is volume
preserving, only the quadrature error remains, and the same holds for the
$L^2$ norm. The MLP flow is not volume preserving and, as all flows in this
paper, is used without its Jacobian determinant (\cref{rem:absolute_value}),
so that its mass conservation is not guaranteed by construction. The energy
drift of the SympNet flow ($\sim 10^{-5}$) is comparable to that
of the reference and fifty times smaller than that of the MLP flow. The
Eulerian PINN drifts on every invariant.

On an NVIDIA A100 GPU, training takes about $63$ minutes for the six MLP
submaps, $214$ minutes for the SympNet ones, and $28$ minutes for the
Eulerian PINN. The cost of a submap grows linearly with its index, since
evaluating the field on the $k$-th subinterval requires the composition of
all previous submaps: from $286$\,s for the first MLP submap to $949$\,s for
the sixth.

\paragraph{Two-stream instability.}
The initial condition is the sum of two Maxwellians centered at $\pm v_0$,
\begin{equation*}
    f_0(x,v) = \big(1 + \varepsilon \cos(k x)\big)\,
    \frac{e^{-(v-v_0)^2/2} + e^{-(v+v_0)^2/2}}{2\sqrt{2\pi}},
    \qquad \varepsilon = 0.05, \quad k = 0.3, \quad v_0 = 3,
\end{equation*}
the computational domain is $[0,L] \times [-11,11]$ with  $L = 2\pi/k$.
This is an unstable equilibrium of the system and a small perturbation leads to an exponential increase of the electric energy: this is the linear regime. Eventually the electric energy saturates and the system enters the nonlinear regime. The goal of this test is to correctly reconstruct these dynamics.

We keep the numerical parameters of the Landau damping test: $T = 24$,
$N_t = 6$ subintervals of length $4$, the same MLP ($N_\theta = 1142$) and
SympNet ($N_\theta = 1144$) submaps trained for $400$ steps on $5000$
collocation and $2000$ initial-condition points, $128$ quadrature points in
$v$, $32$ points in $x$ and $N_\tau = 80$ for the field. The Eulerian PINN has
the same architecture and training as for Landau damping. The reference uses a $512 \times 512$ grid
and $\Delta t = 0.05$.

\begin{figure}[!ht]
    \centering
    \input{figures/vlasov_poisson/vp_colors}\ref*{vplegendtwostreamenergy}\\[2pt]
    \begin{minipage}[c]{0.55\linewidth}\centering
        \setlength{\figurewidth}{0.98\linewidth}
        \setlength{\figureheight}{0.78\linewidth}
        %
    \definecolor{mplcolor0}{rgb}{0.000000,0.000000,0.000000}
\definecolor{mplcolor1}{rgb}{0.121570,0.466670,0.705880}
\definecolor{mplcolor2}{rgb}{0.839220,0.152940,0.156860}
\definecolor{mplcolor3}{rgb}{0.172550,0.627450,0.172550}
\begin{tikzpicture}
\begin{axis}[
  width=\figurewidth,
  height=\figureheight,
  ymode=log, log basis y=10,
  xlabel={$t$},
  ylabel={$\frac{1}{2}\int E^2\,\mathrm{d}x$},
  xmin=-1.2, xmax=25.2,
  ymin=1.40203e-06, ymax=11.4554,
  xmajorgrids=true,
  ymajorgrids=true,
  legend to name=vplegendtwostreamenergy,
  legend columns=4
]
\addplot [color=mplcolor0, mark=none, line width=1.20pt]
  coordinates {(0,0.145444) (0.3,0.122141) (0.6,0.0687771) (0.9,0.0203336) (1.2,0.000467465) (1.5,0.00492484) (1.8,0.0122638) (2.15,0.00714801) (2.45,0.000116141) (2.75,0.00770374) (3.05,0.0327305) (3.35,0.0594484) (3.65,0.0682073) (3.95,0.0537277) (4.25,0.028349) (4.55,0.00868207) (4.85,0.000979043) (5.15,1.07583e-05) (5.45,0.000356881) (5.75,0.005603) (6.1,0.0293807) (6.4,0.065619) (6.7,0.102579) (7,0.122249) (7.3,0.116416) (7.6,0.0921146) (7.9,0.0644207) (8.2,0.0452609) (8.5,0.039251) (8.8,0.0482765) (9.1,0.0765462) (9.4,0.128296) (9.7,0.19996) (10.05,0.287681) (10.35,0.342393) (10.65,0.364598) (10.95,0.357735) (11.25,0.338056) (11.55,0.325694) (11.85,0.337849) (12.15,0.387254) (12.45,0.482254) (12.75,0.623943) (13.05,0.801334) (13.35,0.990466) (13.65,1.16183) (13.95,1.29347) (14.3,1.39299) (14.6,1.4521) (14.9,1.51891) (15.2,1.62534) (15.5,1.79634) (15.8,2.04418) (16.1,2.36404) (16.4,2.7321) (16.7,3.10897) (17,3.45048) (17.3,3.72253) (17.6,3.91329) (17.9,4.03651) (18.25,4.13883) (18.55,4.2349) (18.85,4.37055) (19.15,4.55833) (19.45,4.7901) (19.75,5.03851) (20.05,5.26359) (20.35,5.42404) (20.65,5.49015) (20.95,5.45382) (21.25,5.33106) (21.55,5.15579) (21.85,4.9683) (22.2,4.77948) (22.5,4.66501) (22.8,4.59513) (23.1,4.55434) (23.4,4.51704) (23.7,4.45548) (24,4.34853)};
\addplot [color=mplcolor1, mark=none, line width=1.20pt]
  coordinates {(0,0.145637) (0.303797,0.121895) (0.607595,0.0680893) (0.911392,0.0196644) (1.21519,0.000357785) (1.51899,0.0052351) (1.82278,0.012425) (2.12658,0.00824138) (2.43038,0.000372958) (2.73418,0.00644915) (3.03797,0.0309892) (3.34177,0.0583258) (3.64557,0.0681593) (3.94937,0.0539734) (4.25316,0.0283029) (4.55696,0.00842377) (4.86076,0.000862477) (5.16456,2.38216e-05) (5.46835,0.000396221) (5.77215,0.00607212) (6.07595,0.0260819) (6.37975,0.0615251) (6.68354,0.0993603) (6.98734,0.120625) (7.29114,0.115787) (7.59494,0.091701) (7.89873,0.0637039) (8.20253,0.0436616) (8.50633,0.0366992) (8.81013,0.0451153) (9.11392,0.0739267) (9.41772,0.128231) (9.72152,0.203967) (10.0253,0.283161) (10.3291,0.341811) (10.6329,0.365201) (10.9367,0.356859) (11.2405,0.334234) (11.5443,0.318039) (11.8481,0.325089) (12.1519,0.37027) (12.4557,0.465547) (12.7595,0.610672) (13.0633,0.794828) (13.3671,0.991247) (13.6709,1.16587) (13.9747,1.29445) (14.2785,1.37568) (14.5823,1.43081) (14.8861,1.49263) (15.1899,1.59367) (15.4937,1.75831) (15.7975,1.99709) (16.1013,2.29894) (16.4051,2.66736) (16.7089,3.06067) (17.0127,3.42053) (17.3165,3.70051) (17.6203,3.89194) (17.9241,4.02109) (18.2278,4.12369) (18.5316,4.2297) (18.8354,4.37282) (19.1392,4.55677) (19.443,4.79799) (19.7468,5.03231) (20.0506,5.19041) (20.3544,5.37398) (20.6582,5.49382) (20.962,5.51986) (21.2658,5.44216) (21.5696,5.28079) (21.8734,5.07393) (22.1772,4.86517) (22.481,4.68688) (22.7848,4.55447) (23.0886,4.45336) (23.3924,4.36497) (23.6962,4.28511) (24,4.21692)};
\addplot [color=mplcolor2, mark=none, line width=1.20pt]
  coordinates {(0,0.145444) (0.303797,0.121484) (0.607595,0.0671963) (0.911392,0.0189528) (1.21519,0.000254816) (1.51899,0.0055838) (1.82278,0.0125594) (2.12658,0.00797422) (2.43038,0.000283772) (2.73418,0.00671226) (3.03797,0.0312445) (3.34177,0.0582101) (3.64557,0.0674245) (3.94937,0.0530938) (4.25316,0.0276035) (4.55696,0.00812688) (4.86076,0.000804192) (5.16456,2.89036e-06) (5.46835,0.000381104) (5.77215,0.00608016) (6.07595,0.0260907) (6.37975,0.0613411) (6.68354,0.0986342) (6.98734,0.119272) (7.29114,0.114201) (7.59494,0.0902858) (7.89873,0.0626441) (8.20253,0.0434976) (8.50633,0.0376953) (8.81013,0.0469527) (9.11392,0.0755922) (9.41772,0.127811) (9.72152,0.19954) (10.0253,0.27454) (10.3291,0.331165) (10.6329,0.355536) (10.9367,0.349776) (11.2405,0.329443) (11.5443,0.314944) (11.8481,0.324981) (12.1519,0.374104) (12.4557,0.4699) (12.7595,0.613316) (13.0633,0.792183) (13.3671,0.981116) (13.6709,1.14984) (13.9747,1.27662) (14.2785,1.35902) (14.5823,1.41595) (14.8861,1.47976) (15.1899,1.58406) (15.4937,1.75411) (15.7975,2.00129) (16.1013,2.32253) (16.4051,2.69439) (16.7089,3.07671) (17.0127,3.42571) (17.3165,3.70349) (17.6203,3.89941) (17.9241,4.02601) (18.2278,4.11792) (18.5316,4.21621) (18.8354,4.36674) (19.1392,4.56054) (19.443,4.81422) (19.7468,5.06647) (20.0506,5.29055) (20.3544,5.47934) (20.6582,5.55669) (20.962,5.52779) (21.2658,5.41987) (21.5696,5.2756) (21.8734,5.10419) (22.1772,4.94404) (22.481,4.81528) (22.7848,4.7134) (23.0886,4.64247) (23.3924,4.58304) (23.6962,4.47924) (24,4.31606)};
\addplot [color=mplcolor3, mark=none, line width=1.20pt]
  coordinates {(0,0.143874) (0.303797,0.102444) (0.607595,0.0639575) (0.911392,0.0361549) (1.21519,0.0195837) (1.51899,0.011306) (1.82278,0.00812511) (2.12658,0.00815243) (2.43038,0.0110399) (2.73418,0.0175614) (3.03797,0.0290365) (3.34177,0.0468378) (3.64557,0.07206) (3.94937,0.105368) (4.25316,0.146998) (4.55696,0.196874) (4.86076,0.254757) (5.16456,0.320378) (5.46835,0.39349) (5.77215,0.473863) (6.07595,0.561228) (6.37975,0.655207) (6.68354,0.755258) (6.98734,0.860636) (7.29114,0.97038) (7.59494,1.08331) (7.89873,1.19805) (8.20253,1.31305) (8.50633,1.42662) (8.81013,1.537) (9.11392,1.64241) (9.41772,1.74107) (9.72152,1.83133) (10.0253,1.91167) (10.3291,1.98078) (10.6329,2.03757) (10.9367,2.08121) (11.2405,2.11114) (11.5443,2.12707) (11.8481,2.12893) (12.1519,2.11692) (12.4557,2.09144) (12.7595,2.05311) (13.0633,2.0027) (13.3671,1.94116) (13.6709,1.86956) (13.9747,1.78908) (14.2785,1.70098) (14.5823,1.60658) (14.8861,1.50723) (15.1899,1.40428) (15.4937,1.29909) (15.7975,1.19295) (16.1013,1.08712) (16.4051,0.982771) (16.7089,0.880989) (17.0127,0.782752) (17.3165,0.68892) (17.6203,0.600235) (17.9241,0.517312) (18.2278,0.440642) (18.5316,0.370595) (18.8354,0.307424) (19.1392,0.251281) (19.443,0.202219) (19.7468,0.160211) (20.0506,0.125159) (20.3544,0.0969049) (20.6582,0.0752455) (20.962,0.0599409) (21.2658,0.0507244) (21.5696,0.0473098) (21.8734,0.0493975) (22.1772,0.0566789) (22.481,0.0688391) (22.7848,0.085558) (23.0886,0.106511) (23.3924,0.131367) (23.6962,0.159791) (24,0.191438)};
\addlegendentry{SL}
\addlegendentry{MLP flow}
\addlegendentry{SympNet flow}
\addlegendentry{Eulerian PINN}
\end{axis}
\end{tikzpicture}%

    \end{minipage}\hfill
    \begin{minipage}[c]{0.43\linewidth}\centering
        \setlength{\figurewidth}{0.95\linewidth}
        \setlength{\figureheight}{0.5\linewidth}
        %
    \input{figures/vlasov_poisson/two_stream_final_rho.tex}%
\\[6pt]
        %
    \input{figures/vlasov_poisson/two_stream_final_E.tex}%

    \end{minipage}\\[6pt]
    {\phasespacebbox
        \setlength{\figurewidth}{0.97\linewidth}%
        \setlength{\figureheight}{0.2626\linewidth}%
        %
    {\pgfplotsset{compat=1.18}
\pgfplotsset{colormap={turbo}{rgb255(0.0000)=(48,18,59) rgb255(0.0159)=(53,30,88) rgb255(0.0317)=(57,42,115) rgb255(0.0476)=(61,53,139) rgb255(0.0635)=(64,64,162) rgb255(0.0794)=(66,75,181) rgb255(0.0952)=(68,86,199) rgb255(0.1111)=(70,97,214) rgb255(0.1270)=(70,107,227) rgb255(0.1429)=(71,118,238) rgb255(0.1587)=(70,128,246) rgb255(0.1746)=(69,138,252) rgb255(0.1905)=(66,148,255) rgb255(0.2063)=(61,158,254) rgb255(0.2222)=(55,168,250) rgb255(0.2381)=(47,178,244) rgb255(0.2540)=(39,190,233) rgb255(0.2698)=(32,199,223) rgb255(0.2857)=(27,208,213) rgb255(0.3016)=(24,215,202) rgb255(0.3175)=(24,222,192) rgb255(0.3333)=(26,228,182) rgb255(0.3492)=(32,234,172) rgb255(0.3651)=(42,239,161) rgb255(0.3810)=(53,243,148) rgb255(0.3968)=(67,247,135) rgb255(0.4127)=(82,250,122) rgb255(0.4286)=(97,252,108) rgb255(0.4444)=(113,254,95) rgb255(0.4603)=(128,255,83) rgb255(0.4762)=(143,255,73) rgb255(0.4921)=(156,254,64) rgb255(0.5079)=(169,251,57) rgb255(0.5238)=(180,248,54) rgb255(0.5397)=(190,244,52) rgb255(0.5556)=(200,239,52) rgb255(0.5714)=(210,233,53) rgb255(0.5873)=(219,226,54) rgb255(0.6032)=(227,219,56) rgb255(0.6190)=(235,211,57) rgb255(0.6349)=(241,203,58) rgb255(0.6508)=(246,195,58) rgb255(0.6667)=(250,186,57) rgb255(0.6825)=(252,177,54) rgb255(0.6984)=(254,167,50) rgb255(0.7143)=(254,155,45) rgb255(0.7302)=(254,144,41) rgb255(0.7460)=(252,132,35) rgb255(0.7619)=(249,117,29) rgb255(0.7778)=(245,105,24) rgb255(0.7937)=(241,93,19) rgb255(0.8095)=(236,83,15) rgb255(0.8254)=(231,73,12) rgb255(0.8413)=(225,65,9) rgb255(0.8571)=(218,57,7) rgb255(0.8730)=(210,49,5) rgb255(0.8889)=(202,42,4) rgb255(0.9048)=(193,35,2) rgb255(0.9206)=(183,29,2) rgb255(0.9365)=(172,23,1) rgb255(0.9524)=(161,18,1) rgb255(0.9683)=(149,13,1) rgb255(0.9841)=(136,8,2) rgb255(1.0000)=(122,4,3)}}
\begin{tikzpicture}
\begin{axis}[
  at={(0.0211\figurewidth, 0.0779\figureheight)},
  anchor=south west,
  width=0.2250\figurewidth,
  height=0.8312\figureheight,
  xlabel={$x$},
  ylabel={$v$},
  title={SL},
  xmin=0, xmax=20.944,
  ymin=-11, ymax=11,
  xtick=\empty,
  ytick=\empty,
  point meta min=0,
  point meta max=0.215762,
  colormap name=turbo,
  title style={font=\scriptsize, text height=1.6ex, text depth=0.25ex},
  label style={font=\scriptsize},
  scale only axis
]
\addplot graphics [
  xmin=0, xmax=20.944,
  ymin=-11, ymax=11
] {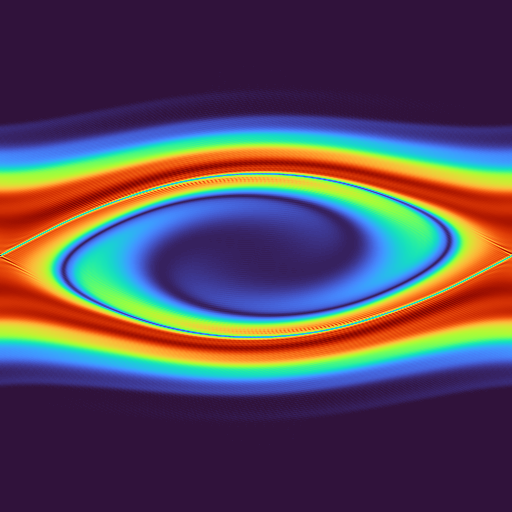};
\end{axis}

\begin{axis}[
  at={(0.2482\figurewidth, 0.0779\figureheight)},
  anchor=south west,
  width=0.2250\figurewidth,
  height=0.8312\figureheight,
  xlabel={$x$},
  title={MLP flow},
  xmin=0, xmax=20.944,
  ymin=-11, ymax=11,
  xtick=\empty,
  ytick=\empty,
  point meta min=0,
  point meta max=0.215762,
  colormap name=turbo,
  title style={font=\scriptsize, text height=1.6ex, text depth=0.25ex},
  label style={font=\scriptsize},
  scale only axis
]
\addplot graphics [
  xmin=0, xmax=20.944,
  ymin=-11, ymax=11
] {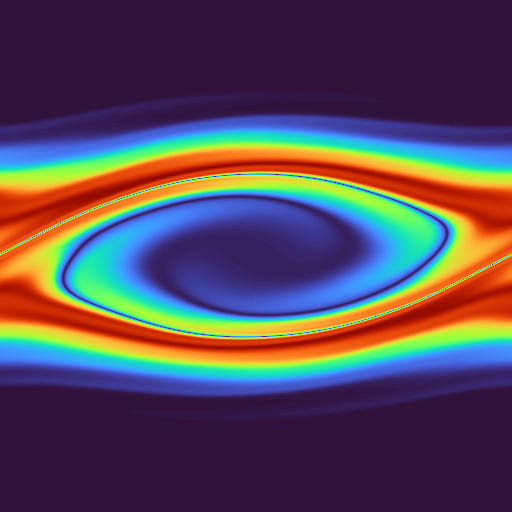};
\end{axis}

\begin{axis}[
  at={(0.4754\figurewidth, 0.0779\figureheight)},
  anchor=south west,
  width=0.2250\figurewidth,
  height=0.8312\figureheight,
  xlabel={$x$},
  title={SympNet flow},
  xmin=0, xmax=20.944,
  ymin=-11, ymax=11,
  xtick=\empty,
  ytick=\empty,
  point meta min=0,
  point meta max=0.215762,
  colormap name=turbo,
  title style={font=\scriptsize, text height=1.6ex, text depth=0.25ex},
  label style={font=\scriptsize},
  scale only axis
]
\addplot graphics [
  xmin=0, xmax=20.944,
  ymin=-11, ymax=11
] {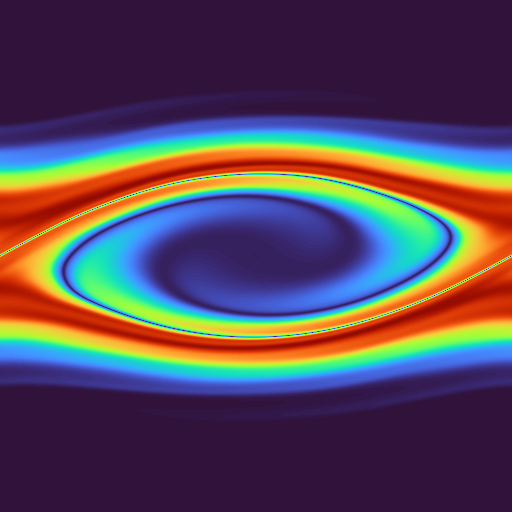};
\end{axis}

\begin{axis}[
  at={(0.7025\figurewidth, 0.0779\figureheight)},
  anchor=south west,
  width=0.2250\figurewidth,
  height=0.8312\figureheight,
  xlabel={$x$},
  title={Eulerian PINN},
  xmin=0, xmax=20.944,
  ymin=-11, ymax=11,
  xtick=\empty,
  ytick=\empty,
  point meta min=0,
  point meta max=0.215762,
  colormap name=turbo,
  colorbar,
  colorbar style={width=0.0127\figurewidth},
  title style={font=\scriptsize, text height=1.6ex, text depth=0.25ex},
  label style={font=\scriptsize},
  scale only axis
]
\addplot graphics [
  xmin=0, xmax=20.944,
  ymin=-11, ymax=11
] {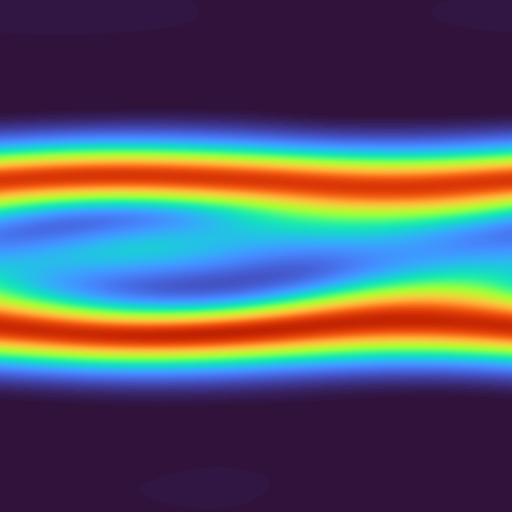};
\end{axis}
\end{tikzpicture}
}%
}
    \caption{Two-stream instability from \cref{sec:vlasov_poisson}.
        Top left: electric energy $\frac{1}{2}\int E^2\,\dx$ over time.
        Top right: density perturbation $\rho-1$ and electric field $E$
        at the final time $T = 24$.
        Bottom: distribution function $f$ at $T = 24$.}
    \label{fig:vp-ts}
\end{figure}

Both learned flows follow the electric energy through the growth, the saturation
and the subsequent slow decay (\cref{fig:vp-ts}), the SympNet flow more
closely than the MLP flow. The pointwise error, concentrated on the boundary of the
vortex, is not shown, since the reference solution is itself subject to numerical diffusion there. The
Eulerian PINN fails on this test: neither the electric energy nor $f$ is reproduced.
The drifts of the invariants (\cref{fig:vp-conservation}, bottom) follow the
Landau damping case until the saturation, then grow as the filaments become
finer. The SympNet flow remains about an order of magnitude more accurate
than the MLP flow on the energy and the momentum, and at $T = 24$ the drift of
the $L^2$ norm is $2\times 10^{-2}$ for the MLP flow, $5 \times 10^{-4}$ for
the SympNet flow and $2 \times 10^{-5}$ for the reference. Training takes
about $64$ minutes for the MLP submaps, $226$ minutes for the SympNet ones
and $28$ minutes for the Eulerian PINN.

\subsection{An axisymmetric drift-kinetic Vlasov-Poisson model in a tokamak}
\label{sec:gyrokinetic-test-case}

As a last test case, we present preliminary results on a
reduced, single-species, axisymmetric, electrostatic drift-kinetic model with adiabatic
electrons in a large-aspect-ratio circular tokamak.

The drift-kinetic electrostatic equations in axisymmetric geometry read~\cite{GyselaX-semi-Lagrangian2016,grandgirard2008computing}
\begin{align}
    \partial_t f
    + \dot{r}\,\partial_r f
    + \dot{\theta}\,\partial_\theta f
    + \dot v_\parallel\,\partial_{v_\parallel} f                          & = 0,                  \label{eq:vlasov} \\
    - \rho_*^2 a^2\, \Delta_{\perp} \phi + [ \phi - \langle \phi \rangle] & = \delta n, \label{eq:quasineutrality}
\end{align}
with $\Delta_{\perp} = \frac{1}{r}\,\partial_r\left( r\,\partial_r \right) + \frac{1}{r^{2}}\,\partial_\theta^{2}$.
The first equation is the Vlasov equation for the ion guiding-center
distribution function $f(t,r,\theta,v_\parallel,\mu)$. The characteristic curves are defined as
\begin{align}
    \dot r & = -\frac{1}{s\,rB}\left(\partial_\theta \phi + \frac{v_\parallel^{2} + \mu B}{B}\,\partial_\theta B\right),
    \qquad
    \dot\theta = \frac{1}{s\,rB}\left(\partial_r \phi + \frac{v_\parallel^{2} + \mu B}{B}\,\partial_r B\right)
    + \frac{v_\parallel}{s\,q R_0},
    \label{eq:characteristics-xy}                                                                                        \\
    \dot v_\parallel
           & = - \frac{1}{s\,q R_0}\left(\mu\,\partial_\theta B + \partial_\theta \phi\right)
    \;+\; v_\parallel\,\frac{\partial_r\phi\,\partial_\theta B - \partial_\theta\phi\,\partial_r B}{s\,rB^2},
    \label{eq:vpar-dot}
\end{align}
with
\begin{equation}
    B(r,\theta) = \frac{B_0 R_0\, s(r)}{R(r,\theta)}, \qquad
    s(r) = \sqrt{1 + \left(\frac{r}{qR_0}\right)^{2}}, \qquad
    R(r,\theta) = R_0 + r\cos\theta .
    \label{eq:field}
\end{equation}
The coefficients $\rho_*$, $a$, $B_0$, $R_0$ and $q$ are physical parameters given below.
In this model, $\mu$ is a parameter, meaning that $\dot \mu = 0$.
In the field equation, $\langle\cdot\rangle$ is the flux-surface average,
\begin{equation}
    \langle\phi\rangle(r) =
    \frac{\int_0^{2\pi} \phi(r,\theta)\,r\,R(r,\theta)\,d\theta}
    {\int_0^{2\pi} r\,R(r,\theta) \,d\theta},
    \label{eq:flux-surface-average}
\end{equation}
and $\delta n$ is the ion guiding-center density,
\begin{equation}
    \delta n(t,r,\theta) =
    \int_{\mathbb{R}_+} \!\! \int_{\mathbb{R}}
    2\pi\, B(r,\theta)\,\left[ f(t, r,\theta,v_\parallel,\mu)- f_{\mathrm{eq}}( r,\theta,v_\parallel,\mu) \right]
    \; dv_\parallel \, d\mu,
    \label{eq:density-moment}
\end{equation}
where the equilibrium distribution $f_{\mathrm{eq}}$ is
$$
    f_{\mathrm{eq}}( r,\theta,v_\parallel,\mu) =   \frac{1}{(2\pi)^{3/2}}
    \exp\!\left(-\frac{v_\parallel^{2}}{2} - \mu B(r,\theta)\right).
$$
The field
equation is solved with homogeneous Dirichlet conditions at both radial boundaries and
periodicity in $\theta$,
\begin{equation}
    \phi\big|_{r = r_{\min}} = \phi\big|_{r = r_{\max}} = 0.
    \label{eq:bcs}
\end{equation}

\paragraph{The Rosenbluth-Hinton test.} Following the Rosenbluth-Hinton test of~\cite{GyselaX-semi-Lagrangian2016,rosenbluth1998poloidal},
the initial condition is a Maxwellian with uniform density and temperature,
perturbed by a purely radial profile that vanishes at both radial
boundaries,
\begin{equation}
    f_0(r,\theta,v_\parallel,\mu) =
    f_{\mathrm{eq}}( r,\theta,v_\parallel,\mu)
    \left(1 + \varepsilon \sin\!\left(\pi\,\frac{r - r_{\min}}{r_{\max} - r_{\min}}\right)\right).
    \label{eq:initial-condition}
\end{equation}
The parameters are those
of~\cite{GyselaX-semi-Lagrangian2016}: major radius $R_0 = 1$, minor radius
$a = R_0/2.78$, safety factor $q = 1.9$, perturbation amplitude
$\varepsilon = 10^{-3}$ and radial domain $[r_{\min}, r_{\max}] = [0.2\,a, 0.8\,a]$.
The magnetic field is set by the normalized Larmor radius
$\rho_* = \rho_L/a = 10^{-2}$ through $B_0 = 1/(\rho_* a) \approx 278$.
The velocity variables are truncated to $|v_\parallel| \le 7$ and
$\mu \le 12/B_0$, and the final time is $T = 40$.

The time interval is split into $N_t = 80$ subintervals of length $0.5$.
Each submap is an MLP of width $w=64$ and depth $\ell = 3$,
returning $(r,\theta,v_\parallel)$, the magnetic moment $\mu$ entering as a
parameter. It is trained for $N_{\mathrm{epoch}} = 400$ natural-gradient
steps on $16000$
collocation points and $6000$ initial-condition points. The density
moment~\eqref{eq:density-moment} is computed with a uniform $128 \times 8$
quadrature in $(v_\parallel,\mu)$, and the field equation is solved with Q1
finite elements on a $32 \times 32$ grid in $(r,\theta)$ at $N_\tau = 25$
times per subinterval.

The quantity of interest is the potential averaged over each magnetic
surface~\eqref{eq:flux-surface-average}, taken at mid-radius $r^\star = a/2$
and normalized by its initial value. Starting from the radial
perturbation~\eqref{eq:initial-condition}, it oscillates at a well-defined
frequency $\omega_G$ (the so-called geodesic acoustic mode); the oscillation is damped by
a mechanism similar to Landau damping, and the average does not return to
zero but to a positive residual level $A_R$~\cite{rosenbluth1998poloidal}. As is
usual for this test, we summarize the curve by three numbers, obtained by a
least-squares fit of a damped cosine function
\begin{equation}
    \frac{\langle\phi\rangle(r^\star,t)}{\langle\phi\rangle(r^\star,0)}
    \simeq A_R + (1 - A_R)\cos(\omega_G t)\,e^{-\gamma_G t}:
    \label{eq:rh-fit}
\end{equation}
the residual level $A_R$, the frequency $\omega_G$ and the damping rate
$\gamma_G$. The value of $A_R$ is computed by an analytical formula, due to Rosenbluth and
Hinton~\cite{rosenbluth1998poloidal}, which gives $A_R = 0.068$ here.
The theoretical values of $\omega_G$ and $\gamma_G$ rely on asymptotic
approximations. We use the
values $\omega_G = 2.53$ and $\gamma_G = 0.0685$ measured with the
semi-Lagrangian code GYSELA~\cite{GyselaX-semi-Lagrangian2016}. It should be noted that GYSELA is a full gyrokinetic code, while we only consider a simplified axisymmetric drift-kinetic configuration here.

\begin{figure}[htbp]
    \centering
    \includegraphics[width=\linewidth]{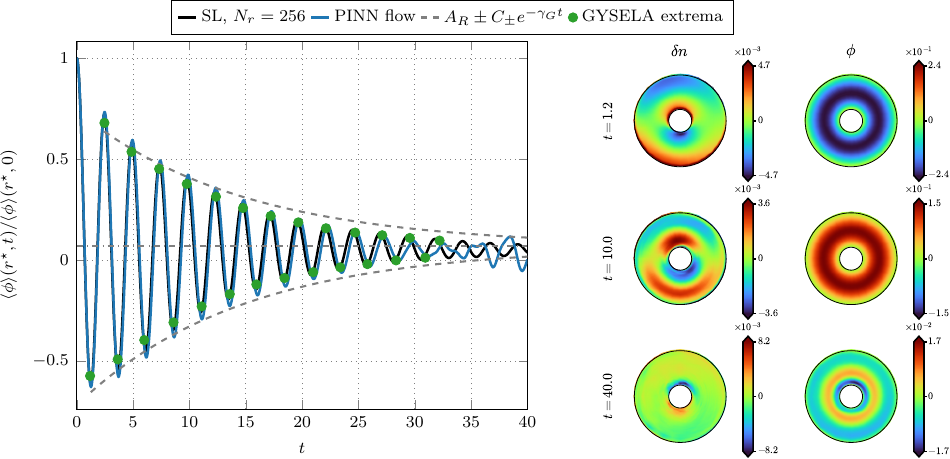}
    \caption{Rosenbluth-Hinton test over $t \in [0, 40]$ ($80$ submaps).
        Left: normalized averaged potential at $r^\star = a/2$, PINN, backward
        semi-Lagrangian reference (SL, $256 \times 64$ grid in $(r,\theta)$,
        $256 \times 8$ velocity quadrature, $\Delta t = 0.02$), GYSELA
        extrema~\cite{GyselaX-semi-Lagrangian2016}, residual level $A_R$
        (dash-dotted) and envelope
        $A_R \pm C_\pm e^{-\gamma_G t}$ with the Rosenbluth-Hinton $A_R$ and
        the GYSELA damping rate $\gamma_G$, its height $C_\pm$ set on the third
        maximum and minimum of the SL trace.
        Right: poloidal cuts of the PINN density perturbation $\delta n$ and
        potential $\phi$ at $t = 1.2$, $10$ and $40$, each with its own
        symmetric color scale (set on the interior radial nodes for
        $\delta n$); black circles mark $r_{\min}$ and $r_{\max}$.}
    \label{fig:rh-residual}
\end{figure}

\Cref{fig:rh-residual} shows that the PINN follows the oscillation
over $T = 40$, about $16$ periods. Fitting~\eqref{eq:rh-fit} on the whole
curve gives $\omega_G = 2.54$, within $0.5\%$ of GYSELA, and
$\gamma_G = 0.099$, about $45\%$ above the GYSELA value but in agreement with
our semi-Lagrangian reference computed on a much finer grid
($\omega_G = 2.57$, $\gamma_G = 0.102$). The same fit gives a residual level
$A_R = 0.052$ for the PINN and $A_R = 0.061$ for SL, below the
Rosenbluth-Hinton value $A_R = 0.068$, a deviation also observed with GYSELA
in this configuration.

\paragraph{Reconstruction of the distribution function.}
Write $z = (r,\theta,v_\parallel)$ and $E(z) = v_\parallel^2/2 + \mu B(r,\theta)$
for the energy. The initial condition~\eqref{eq:initial-condition} reads
$f_0(z) = f_{\mathrm{eq}}(z)\,\bigl(1 + \varepsilon \sin(\pi\,\frac{r - r_{\min}}{r_{\max} - r_{\min}})\bigr)$,
where $f_{\mathrm{eq}}(z) = (2\pi)^{-3/2} e^{-E(z)}$ is the equilibrium
appearing in~\eqref{eq:density-moment}. For $t \in [t_k, t_{k+1}]$, the points
$z_k, \dots, z_0$ are obtained by the backward
composition~\eqref{eq:discrete_composition}, and the pullback of $f_0$ gives
\begin{equation}
    f(t,z) = \frac{1}{(2\pi)^{3/2}}\, e^{-E_0}
    \Bigl(1 + \varepsilon \sin\Bigl(\pi\,\frac{r_0 - r_{\min}}{r_{\max} - r_{\min}}\Bigr)\Bigr),
    \label{eq:rh-reconstruction}
\end{equation}
where $r_0$ is the radial component of $z_0$ and $E_0$ is the energy at
$z_0$. This construction differs from the previous test cases in two
respects.

\begin{itemize}
    \item \emph{Radial boundary conditions.} In contrast with the previous
          test cases, where the domain was periodic, the
          characteristics~\eqref{eq:characteristics-xy} cross the boundaries
          $r = r_{\min}$ and $r = r_{\max}$, since $\dot r$ does not vanish
          there. The transport equation~\eqref{eq:vlasov} therefore requires
          inflow boundary conditions, which are not part of the model. In
          semi-Lagrangian codes, and in GYSELA in particular, a foot of a
          characteristic located outside the domain is given the value of $f$
          at the nearest boundary
          point~\cite{BourneMunschyGrandgirardMehrenbergerGhendrih2023}. The
          same question arises for the composition~\eqref{eq:discrete_composition}:
          each submap is trained on the domain only, but may receive a point
          that the previous submap has mapped outside of it. We impose the
          same condition: the radial component $r_i$ of each point $z_i$,
          $i = k, \dots, 0$, is replaced by
          $\bar r_i = \min\bigl(\max(r_i, r_{\min}), r_{\max}\bigr)$ before
          $z_i$ is used, either as the argument of the next submap or
          in~\eqref{eq:rh-reconstruction}. A backward characteristic may thus
          remain on the boundary over some subintervals and re-enter the
          domain afterwards. The projection is applied at the times $t_i$
          only, whereas semi-Lagrangian codes apply it at every time step.
    \item \emph{Computation of the density.} The density~\eqref{eq:density-moment}
          is of order $\varepsilon$, as the integral of the difference between
          $f$ and $f_{\mathrm{eq}}$, both of order one. The energy $E_0$ must
          therefore be computed with an accuracy small compared with
          $\varepsilon$. Evaluating $E$ at the computed point $z_0$ does not
          meet this requirement, since the approximation error of the
          characteristic is transferred to $E_0$ without any reduction.
          Instead, since $E(z)$ is known exactly, we only compute the
          variation of $E$ along the characteristic,
          \begin{equation*}
              \dot E = v_\parallel\,\dot v_\parallel
              + \mu\,\bigl(\dot r\,\partial_r B + \dot\theta\,\partial_\theta B\bigr),
          \end{equation*}
          with $\dot r$, $\dot\theta$ and $\dot v_\parallel$ given
          by~\eqref{eq:characteristics-xy}--\eqref{eq:vpar-dot}. Substituting
          these expressions shows that $\dot E$ vanishes for $\phi = 0$ and is
          linear in $\phi$, which is of order $\varepsilon$: an approximation
          error on the characteristic only results in an error multiplied by
          $\varepsilon$ on the variation of $E$. In addition, each projection
          keeps $v_\parallel$ but changes $r$, so that $E$ is discontinuous at
          the times $t_i$. Hence
          \begin{equation}
              E_0 = E(z) - \int_0^t \dot E\,ds
              + \sum_{i=0}^{k} \mu\,\bigl[B(\bar r_i,\theta_i) - B(r_i,\theta_i)\bigr],
              \label{eq:energy-foot}
          \end{equation}
          where the sum accounts for these jumps, which vanish when
          $r_i \in [r_{\min}, r_{\max}]$.
          On each subinterval $[t_i, t_{i+1}]$, the characteristic is
          approximated by the cubic Hermite interpolant of its end points and of
          the velocities~\eqref{eq:characteristics-xy}--\eqref{eq:vpar-dot} at
          these points, and the integral by the composite trapezoidal rule with
          $8$ sub-steps.

\end{itemize}

\section{Conclusion}
\label{sec:concl}
We have proposed a Lagrangian PINN for transport-dominated problems that
learns the flow map instead of the transported quantity and recovers the
solution by pulling back the exact initial condition. As in characteristic
mapping methods, the map is decomposed into submaps, each learned by its own
network on its own time interval, so that the growing complexity of the
solution is carried by their composition rather than by a single network.
Velocities given by a field equation are handled by a Picard iteration.

The results support three conclusions. First, learning the map rather than
the solution is beneficial: wherever both were run, the Lagrangian PINN clearly
outperforms the Eulerian PINN, by more than an order of magnitude on
advection, and it reproduces the two-stream instability on which the
Eulerian PINN failed. Second, the decomposition into submaps is equally
important: a single global map, as in~\cite{MontiMorozovskaPoggiShukla2026},
loses two to four orders of magnitude in enstrophy conservation on the Euler
tests, while the composition inherits the subgrid property of characteristic
mapping, resolving the spectrum of a $512^2$ reference from a $128^2$ field
grid. Third, structure-preserving networks pay off: at equal parameter
budget, symplectic networks are about an order of magnitude more accurate
than MLPs on Vlasov-Poisson and conserve all $L^p$ norms up to the quadrature error.
Preliminary drift-kinetic results capture the geodesic acoustic mode
frequency within $0.5\%$.

In the low-dimensional settings considered here, the method is not
competitive in cost with grid-based solvers, and the cost of each new submap
grows linearly with the number of previous ones when the velocity depends on
the solution.
Future work will address this by choosing the time partition adaptively.
It will also exploit the parametric
inputs of the networks for many-query problems and extend the drift-kinetic
test to the full five-dimensional gyrokinetic setting, where a meshfree
method is expected to pay off.
Furthermore, we plan to include source terms similar to \cite{YinKrahNaveSchneider2026,YinKrahLinNaveSchneider2026} to simulate plasmas with collisions or MHD turbulence in fusion devices.



\section*{Author Contribution Statement}

The following outlines the authors' contributions to this work.

{\small
\noindent
\begin{tabular}{@{}p{0.23\linewidth} p{0.7\linewidth}}

\textbf{M. Campos Pinto:} &Methodology, Investigation, Revision and Editing, Funding acquisition, Supervision \\
\textbf{V. Fournet:} & Investigation, Implementation, Software, Validation, Visualization, Writing original draft\\
\textbf{E. Franck:} &Initial Idea, Methodology, Implementation, Supervision\\
\textbf{P. Krah:} &Investigation, Implementation, Software, Validation, Visualization, Writing original draft \\
\textbf{V. Michel-Dansac:} &Methodology, Investigation, Implementation, Revision and Editing
\end{tabular}
}

\section*{Code Availability}
The code will be made publicly available upon publication or request.

\section*{Acknowledgments}
As part of the ``France 2030'' initiative, this work has benefited from a national grant managed by the French National
Research Agency (Agence Nationale de la Recherche) attributed to the Exa-MA and Exa-Dost projects of the NumPEx PEPR program, under the references
ANR-22-EXNU-0002, ANR-22-EXNU-0004.

\bibliographystyle{siamplain}
\bibliography{references}

\begin{thebibliography}{10}

\bibitem{ArberVann2002}
{\sc T.~D. Arber and R.~G.~L. Vann}, {\em A {{Critical Comparison}} of {{Eulerian-Grid-Based Vlasov Solvers}}}, J. Comput. Phys., 180 (2002), pp.~339--357, \url{https://doi.org/10.1006/jcph.2002.7098}.

\bibitem{Baty2024}
{\sc H.~Baty}, {\em A hands-on introduction to physics-informed neural networks for solving partial differential equations with benchmark tests taken from astrophysics and plasma physics}, arXiv preprint arXiv:2403.00599,  (2024).

\bibitem{BirdsallLangdon2018}
{\sc C.~K. Birdsall and A.~B. Langdon}, {\em Plasma physics via computer simulation}, CRC press, 2018.

\bibitem{BourneMunschyGrandgirardMehrenbergerGhendrih2023}
{\sc E.~Bourne, Y.~Munschy, V.~Grandgirard, M.~Mehrenberger, and P.~Ghendrih}, {\em Non-uniform splines for semi-{L}agrangian kinetic simulations of the plasma sheath}, J. Comput. Phys., 488 (2023), p.~112229, \url{https://doi.org/10.1016/j.jcp.2023.112229}.

\bibitem{Canosa1973}
{\sc J.~Canosa}, {\em Numerical solution of {{Landau}}'s dispersion equation}, J. Comput. Phys., 13 (1973), pp.~158--160, \url{https://doi.org/10.1016/0021-9991(73)90131-9}.

\bibitem{CanutoHussainiQuarteroniZang2007}
{\sc C.~Canuto, M.~Y. Hussaini, A.~Quarteroni, and T.~A. Zang}, {\em Spectral Methods: Evolution to Complex Geometries and Applications to Fluid Dynamics}, Scientific Computation, Springer, Berlin, Heidelberg, 2007, \url{https://doi.org/10.1007/978-3-540-30728-0}.

\bibitem{ChenKnorr1976}
{\sc C.-Z. Cheng and G.~Knorr}, {\em The integration of the {V}lasov equation in configuration space}, J. Comput. Phys., 22 (1976), pp.~330--351.

\bibitem{deng2023fluid}
{\sc Y.~Deng, H.-X. Yu, D.~Zhang, J.~Wu, and B.~Zhu}, {\em Fluid simulation on neural flow maps}, ACM Trans. Graph., 42 (2023), pp.~1--21.

\bibitem{DeVoreHaninPetrova2021}
{\sc R.~DeVore, B.~Hanin, and G.~Petrova}, {\em Neural network approximation}, Acta Numerica, 30 (2021), pp.~327--444, \url{https://doi.org/10.1017/S0962492921000052}.

\bibitem{FilbetSonnendrucker2003}
{\sc F.~Filbet and E.~Sonnendr{\"u}cker}, {\em Comparison of {E}ulerian {V}lasov solvers}, Comput. Phys. Commun., 150 (2003), pp.~247--266, \url{https://doi.org/10.1016/S0010-4655(02)00694-X}.

\bibitem{fournet2026nonlinearcontrolvariatedeltaf}
{\sc V.~Fournet, M.~{Campos Pinto}, E.~Franck, and V.~Michel-Dansac}, {\em Non-linear control variate in ${\delta}f$ particle-in-cell methods using symplectic neural networks}, 2026, \url{https://arxiv.org/abs/2606.30622}.

\bibitem{scimba}
{\sc E.~Franck, R.~Imbach, and V.~Michel-Dansac}, {\em Scimba python library}, July 2026, \url{https://gitlab.com/scimba/scimba/}.

\bibitem{FranckMichel-DansacNavoretVigon2026}
{\sc E.~Franck, V.~Michel-Dansac, L.~Navoret, and V.~Vigon}, {\em Neural semi-{L}agrangian method for high-dimensional advection-diffusion problems}, Computer Methods in Applied Mechanics and Engineering, 448 (2026), p.~118481.

\bibitem{Golse2013}
{\sc F.~Golse}, {\em Mean field kinetic equations}.
\newblock Lecture notes, École Polytechnique, Palaiseau, 2013.
\newblock Available at \url{https://www.cmls.polytechnique.fr/perso/golse/M2/PolyKinetic.pdf}.

\bibitem{GyselaX-semi-Lagrangian2016}
{\sc V.~Grandgirard, J.~Abiteboul, J.~Bigot, T.~Cartier-Michaud, N.~Crouseilles, G.~Dif-Pradalier, C.~Ehrlacher, D.~Esteve, X.~Garbet, P.~Ghendrih, et~al.}, {\em A {5D} gyrokinetic full-f global semi-{L}agrangian code for flux-driven ion turbulence simulations}, Comput. Phys. Commun., 207 (2016), pp.~35--68.

\bibitem{grandgirard2008computing}
{\sc V.~Grandgirard, Y.~Sarazin, X.~Garbet, G.~Dif-Pradalier, P.~Ghendrih, N.~Crouseilles, G.~Latu, E.~Sonnendr{\"u}cker, N.~Besse, and P.~Bertrand}, {\em Computing {ITG} turbulence with a full-f semi-{L}agrangian code}, Communications in nonlinear science and numerical simulation, 13 (2008), pp.~81--87.

\bibitem{GreifUrban2019}
{\sc C.~Greif and K.~Urban}, {\em Decay of the {{Kolmogorov}} ${N}$-width for wave problems}, arXiv:1903.08488 [math],  (2019), \url{https://arxiv.org/abs/1903.08488}.

\bibitem{HesthavenPeherstorferUnger2026}
{\sc J.~S. Hesthaven, B.~Peherstorfer, and B.~Unger}, {\em Nonlinear model reduction for transport-dominated problems}, Acta Numer., 35 (2026), pp.~173--272, \url{https://doi.org/10.1017/S0962492926100294}, \url{https://arxiv.org/abs/2602.01397}.

\bibitem{HornikStinchcombeWhite1989}
{\sc K.~Hornik, M.~Stinchcombe, and H.~White}, {\em Multilayer feedforward networks are universal approximators}, Neural Networks, 2 (1989), pp.~359--366.

\bibitem{jin2020sympnets}
{\sc P.~Jin, Z.~Zhang, A.~Zhu, Y.~Tang, and G.~E. Karniadakis}, {\em {SympNets}: Intrinsic structure-preserving symplectic networks for identifying {H}amiltonian systems}, Neural Networks, 132 (2020), pp.~166--179.

\bibitem{JniKiyShuUrbAhmMue2026}
{\sc A.~Jnini, E.~Kiyani, K.~Shukla, J.~F. Urbán, N.~Ahmadi~Daryakenari, J.~Müller, M.~Zeinhofer, and G.~E. Karniadakis}, {\em Curvature-aware optimization for high-accuracy physics-informed neural networks}, Comput. Method. Appl. M., 462 (2026), p.~119289, \url{https://doi.org/10.1016/j.cma.2026.119289}.

\bibitem{KirchhartWilhelm2024}
{\sc M.~Kirchhart and R.~P. Wilhelm}, {\em The numerical flow iteration for the {V}lasov--{P}oisson equation}, SIAM J. Sci. Comput., 46 (2024), pp.~A1972--A1997.

\bibitem{KrahLinWilhelmBacchiniNaveGrandgirardSchneider2026}
{\sc P.~Krah, Z.~Lin, R.~Wilhelm, F.~Bacchini, J.-C. Nave, V.~Grandgirard, K.~Schneider, et~al.}, {\em A hybrid semi-{L}agrangian flow mapping approach for {V}lasov systems: Combining iterative and compositional flow maps}, arXiv preprint arXiv:2601.21668,  (2026).

\bibitem{Krah_P_2024_PP_robust_spodpmdfmt}
{\sc P.~Krah, A.~Marmin, B.~Zorawski, J.~Reiss, and K.~Schneider}, {\em A robust shifted proper orthogonal decomposition: Proximal methods for decomposing flows with multiple transports}, SIAM J. Sci. Comput., 47 (2025), pp.~A633--A656, \url{https://doi.org/10.1137/24M164392X}, \url{https://arxiv.org/abs/2403.04313}.

\bibitem{KrahYinBergmannNaveSchneider2024}
{\sc P.~Krah, X.-Y. Yin, J.~Bergmann, J.-C. Nave, and K.~Schneider}, {\em A characteristic mapping method for {V}lasov–{P}oisson with extreme resolution properties}, Commun. Comput. Phys., 35 (2024), p.~905–937, \url{https://doi.org/10.4208/cicp.oa-2024-0012}.

\bibitem{krall1986principles}
{\sc N.~A. Krall and A.~W. Trivelpiece}, {\em Principles of Plasma Physics}, San Francisco Press, San Francisco, 1986.
\newblock Reprint of the 1973 McGraw-Hill edition.

\bibitem{Landau1946}
{\sc L.~D. Landau}, {\em On the vibrations of the electronic plasma}, J. Phys. (USSR), 10 (1946), pp.~25--34.

\bibitem{MajdaBertozzi2002}
{\sc A.~J. Majda and A.~L. Bertozzi}, {\em Vorticity and Incompressible Flow}, vol.~27 of Cambridge Texts in Applied Mathematics, Cambridge University Press, Cambridge, 2002.

\bibitem{MercierYinNave2019}
{\sc O.~Mercier, X.-Y. Yin, and J.-C. Nave}, {\em The characteristic mapping method for the linear advection of arbitrary sets}, SIAM J. Sci. Comput., 42 (2020), pp.~A1663--A1685, \url{https://doi.org/10.1137/18M1234424}.

\bibitem{Mojgani_R_2023_j-comput-meth-appl-mech-eng_kolmogorov_wlpinnccmcdpde}
{\sc R.~Mojgani, M.~Balajewicz, and P.~Hassanzadeh}, {\em {K}olmogorov n-width and {L}agrangian physics-informed neural networks: A causality-conforming manifold for convection-dominated {PDE}s}, Computer Methods in Applied Mechanics and Engineering, 404 (2023), p.~115810, \url{https://doi.org/10.1016/j.cma.2022.115810}.

\bibitem{MontiMorozovskaPoggiShukla2026}
{\sc E.~Monti, K.~Morozovska, A.~Poggi, and K.~Shukla}, {\em {Conservative Semi-Lagrangian Physics-Informed Neural Networks (c{SL}-{PINN}s) for solving hyperbolic Partial Differential Equations ({PDE}s)}}, in AI{\&}PDE: ICLR 2026 Workshop on AI and Partial Differential Equations, 2026, \url{https://openreview.net/forum?id=LYU0dchvlN}.

\bibitem{MueZei2023}
{\sc J.~M\"{u}ller and M.~Zeinhofer}, {\em {Achieving High Accuracy with PINNs via Energy Natural Gradient Descent}}, in Proceedings of the 40th International Conference on Machine Learning, A.~Krause, E.~Brunskill, K.~Cho, B.~Engelhardt, S.~Sabato, and J.~Scarlett, eds., vol.~202 of Proceedings of Machine Learning Research, PMLR, 2023, pp.~25471--25485, \url{https://proceedings.mlr.press/v202/muller23b.html}.

\bibitem{Orszag1971}
{\sc S.~A. Orszag}, {\em On the elimination of aliasing in finite-difference schemes by filtering high-wavenumber components}, Journal of the Atmospheric Sciences, 28 (1971), p.~1074, \url{https://doi.org/10.1175/1520-0469(1971)028<1074:OTEOAI>2.0.CO;2}.

\bibitem{RaissiPerdikarisKarniadakis2019}
{\sc M.~Raissi, P.~Perdikaris, and G.~E. Karniadakis}, {\em Physics-informed neural networks: A deep learning framework for solving forward and inverse problems involving nonlinear partial differential equations}, J. Comput. Phys., 378 (2019), pp.~686--707.

\bibitem{ReissSchulzeSesterhennMehrmann2018}
{\sc J.~Reiss, P.~Schulze, J.~Sesterhenn, and V.~Mehrmann}, {\em The {{Shifted Proper Orthogonal Decomposition}}: {{A Mode Decomposition}} for {{Multiple Transport Phenomena}}}, SIAM J. Sci. Comput., 40 (2018), pp.~A1322--A1344, \url{https://doi.org/10.1137/17M1140571}.

\bibitem{rosenbluth1998poloidal}
{\sc M.~Rosenbluth and F.~Hinton}, {\em Poloidal flow driven by ion-temperature-gradient turbulence in tokamaks}, Phys. Rev. Lett., 80 (1998), p.~724.

\bibitem{SchFur2025}
{\sc N.~Schwencke and C.~Furtlehner}, {\em {AN}a{GRAM}: A natural gradient relative to adapted model for efficient {PINN}s learning}, in The Thirteenth International Conference on Learning Representations, 2025, \url{https://openreview.net/forum?id=o1IiiNIoaA}.

\bibitem{SonnendrueckerRocheBertrandGhizzo1999}
{\sc E.~Sonnendr{\"u}cker, J.~Roche, P.~Bertrand, and A.~Ghizzo}, {\em The semi-{L}agrangian method for the numerical resolution of the {V}lasov equation}, J. Comput. Phys., 149 (1999), pp.~201--220.

\bibitem{vlasov1938vibrational}
{\sc A.~A. Vlasov}, {\em On the vibrational properties of an electron gas}, Zhurnal Eksperimental'noi i Teoreticheskoi Fiziki, 8 (1938), p.~291.

\bibitem{WilhelmEifertTorrilhonOrland2025}
{\sc R.-P. Wilhelm, J.~Eifert, M.~Torrilhon, and F.~Orland}, {\em High fidelity simulations of the multi-species {V}lasov equation in the electro-static, collisional-less limit}, Plasma Phys. Controlled Fusion, 67 (2025), p.~025011, \url{https://doi.org/10.1088/1361-6587/ad9fdb}.

\bibitem{WilhelmKormann2025}
{\sc R.-P. Wilhelm and K.~Kormann}, {\em Restarting the numerical flow iteration through low rank tensor approximations}, arXiv preprint arXiv:2509.08474,  (2025).

\bibitem{WilhelmKrahSchneiderBacchiniGrandgirard2026}
{\sc R.-P. Wilhelm, P.~Krah, K.~Schneider, F.~Bacchini, and V.~Grandgirard}, {\em Revisiting kinetic electrostatic electron non-linear ({KEEN}) waves in the presence of dynamical ions}, 2026, \url{https://arxiv.org/abs/2609.04941}.

\bibitem{YinKrahLinNaveSchneider2026}
{\sc X.-Y. Yin, P.~Krah, Z.~Lin, J.-C. Nave, and K.~Schneider}, {\em Characteristic mapping method for {V}lasov-{P}oisson with {BGK}-collisions}, 2026, \url{https://arxiv.org/abs/2609.16328}.

\bibitem{YinKrahNaveSchneider2026}
{\sc X.-Y. Yin, P.~Krah, J.-C. Nave, and K.~Schneider}, {\em A characteristic mapping method with source terms: Applications to ideal magnetohydrodynamics}, Journal of Computational Physics,  (2026), p.~114766.

\bibitem{YinMercierYadavSchneiderNave2021}
{\sc X.-Y. Yin, O.~Mercier, B.~Yadav, K.~Schneider, and J.-C. Nave}, {\em {A Characteristic Mapping Method for the Two-Dimensional Incompressible {{Euler}} Equations}}, J. Comput. Phys., 424 (2021), p.~109781, \url{https://doi.org/10.1016/j.jcp.2020.109781}.

\bibitem{YinSchneiderNave2023}
{\sc X.-Y. Yin, K.~Schneider, and J.-C. Nave}, {\em A characteristic mapping method for the three-dimensional incompressible {E}uler equations}, J. Comput. Phys., 477 (2023), p.~111876.

\end{thebibliography}




\end{document}